
\input amstex
\input amssym
 
\documentstyle{amsppt}
 
\NoRunningHeads
\magnification=1100
\baselineskip=10pt
\parskip 9pt
\pagewidth{5.2in}
\pageheight{7.2in}

\def\aast{A_{\ast}}
\let\alp=\alpha

\def\a{\Cal A}
\def\aut{\operatorname{Aut}}
\def\abul{\a^{\bullet}}

\def\bl{\bigl(}
\def\br{\bigr)}
\def\bbul{\Cal B^{\bullet}}
\def\bdot{\Cal B^{\bullet}}
\def\bb{\Cal B}

\def\d{\Cal D}
\def\dual{^{\vee}}

\def\cc{\Cal C}
\def\cbul{\Cal C^{\bullet}}
\def\codim{\operatorname{codim}}
\def\calb{{\Cal B}}

\def\Del{\Delta}

\def\e{\Cal E}
\def\Ext{\operatorname{Ext}}

\def\edot{\e^{\bullet}}

\def\ebul{\Cal E^{\bullet}}
\def\ext{\Cal Ext}

\def\fs{f}
\def\f{\Cal F}

\def\fgna{\frak F_{\alpha,g,n}^X}

\def\FF{\frak F}

\def\fdgn{\FF^X_{d,g,n}}

\def\g{\Cal G}
\def\gbul{\g^{\bullet}}

\def\hbul{\Cal H^{\bullet}}
\let\hdot=\hbul
\def\HH{\frak H}
\def\Hbul{\HH^{\bullet}}
\def\Hom{\operatorname{Hom}}
\def\h{\Cal H}

\def\hom{\Cal Hom}

\def\i{\Cal I}

\def\k{\Cal K}
\def\kbul{\k^{\bullet}}

\def\l{\Cal L}

\let\lra=\longrightarrow
\def\lsta{_{\ast}}

\def\lgnax{_{\alpha,g,n}^X}

\def\lagnx{_{\alpha,g,n}^X}
\def\lagnmx{_{\alpha,g,n-1}^X}

\def\mgna{\m_{\alpha,g,n}^X}
\def\m{\Cal M}
\def\mapright#1{\,\smash{\mathop{\lra}\limits^{#1}}\,}
\def\mapto#1{\,\smash{\mathop{\to}\limits^{#1}}\,}
\def\M{{\Cal M}}
\def\mh{\!:\!}
\def\mgn{\M_{g,n}}
\def\mgny{\M_{\alp,g,n}^X}
\def\mgnma{\Cal M_{\alpha,g,n-1}^X}
\def\mdgn{\M_{d,g,n}^X}

\def\mods{\frak{Mod}_S}

\def\mm{\frak m}

\def\mgnx{\m_{\alpha,g,n}^X}
\def\mgnao{\M_{\alpha_1,g_1,n_1+1}^X}
\def\mgnat{\M_{\alpha_2,g_2,n_2+1}^X}
\def\mgnai{\M_{\alpha_i,g_i,n_i+1}^X}
\def\mgno{\M_{g_1,n_1+1}}
\def\mgnt{\M_{g_2,n_2+1}}
\def\mgni{\M_{g_i,n_i+1}}

\let\mgnaot=\mgna

\def\n{\Cal N}

\def\nopo{n_1+1}

\def\o{\Cal O}
\def\om{\Omega}

\def\o{\Cal O}

\def\ob{\operatorname{ob}}

\def\os{{\o_S}}

\def\pri{^{\prime}}

\def\p{\Cal P}

\def\qgna{\Cal Q}
\def\q{\Cal Q}
\def\r{\Cal R}

\def\rdot{\Cal R^{\bullet}}

\def\sub{\subset}
\def\spec{\operatorname{Spec}}
\def\sta{^{\ast}}

\def\s{\Cal S}

\def\sha{^{!}}

\def\simalg{\sim_{\text{alg}}}

\def\t{\Cal T}
\def\tone{\Cal T^1}
\def\ttwo{\Cal T^2}
\def\tdot{\t^{\bullet}}

\def\uvir{^{\text{vir}}}
\def\upmo{^{-1}}

\def\vir{^{\text{vir}}}

\def\v{{\Cal V}}

\def\w{\Cal W}

\def\x{\Cal X}

\def\xs{\x}
\def\xsn{\x^N}

\def\sn{S_N}

\def\y{\Cal Y}

\def\ZZ{\bold Z}
\def\QQ{\bold Q}

\def\CC{\bold C}

\let\pro=\proclaim
\let\endpro=\endproclaim

\def\spec{\operatorname{Spec}}
\def\ker{\operatorname{ker}}
\def\coker{\operatorname{coker}}
\def\rank{\operatorname{rank}}

\def\Ext{\operatorname{Ext}}

\def\newrel#1#2{
        \buildrel
                {\smash{\mathop {\scriptstyle#1}}}   \over 
                {\smash {\raise -0.6ex \hbox {$\textstyle#2$}}}  }

\def\bbpri{\beta\beta\pri}

\def\lb{_{\beta}}
\def\lab{_{\alpha,\beta}}
\def\la{_{\alpha}}
\def\aa{\alpha_1\alpha_2}

\def\ub{U_{\beta}}
\def\wb{W_{\beta}}
\def\wbb{W\lbb}

\def\unm{^{m}}

\def\ual{^{\alpha}}
\def\ubb{U_{\beta\beta\pri}}
\def\aota{\alpha_1+\alpha_2=\alpha}

\def\lkk{_{K}}
\def\lbb{_{\beta\beta\pri}}
\def\lubb{_{U\lbb}}
\def\mgnm{\M_{g,m}}
\def\lb{_{\beta}}
\def\lbp{_{\beta\pri}}
\def\lub{_{\ub}}
\def\lkkb{_{K,\beta}}
\def\lbp{_{\beta\pri}}
\def\lkkbb{_{K,\beta\beta\pri}}
\def\lkkbp{_{K,\beta\pri}}
\def\lubp{_{U\lbp}}
\def\laa{_{\alpha_1\alpha_2}}

\def\and{\text{and}}
\def\lwb{_{\wb}}
\def\lwbb{_{W\lbb}}
\def\ubp{{U\lbp}}
\def\lwb{_{\wb}}

\def\laab{_{\alpha_1\alpha_2,\beta}}
\def\laabb{_{\aa,\beta\beta\pri}}
\def\laabp{_{\alpha_1\alpha_2,\beta\pri}}
\def\lbo{_{\beta,1}}
\def\lfb{_{F\lb}}
\def\labb{_{\alpha,\bbpri}}
\def\lwbp{_{W\lbp}}
\def\lnb{_{n,\beta}}
\def\bp{\beta\pri}
\def\ua{^{\alpha}}
\def\lnmb{_{m,\beta}}
\def\lggnn{_{g_{\bullet}n_{\bullet}}}
\def\lvb{_{V\lb}}

\def\lalp{_{\alpha}}
\def\hhdot{\frak h^{\bullet}}
\def\vt{\text{Vect}}
\let\ecomplex=\edot
\def\oae{\otimes_A E_2}
\def\jko{J_{k-1}}
\def\mko{M_{k+1}}
\def\iko{I_{k-1}}
\def\mo{M_1}

\def\homem{\Hom_A(E_2,M\otimes_AE_2)}
\let\homae=\homem
\def\ide{\text{1}_{E_2}}
\def\oae{\otimes_A E_2}

\def\hh{\frak h}
\def\bob{B\otimes_k B}
\def\hatbob{\hat B}
\def\fdot{\Cal F^{\bullet}}

\let\ohp=\oph
\def\hatfdot{{\hat F}^{\bullet}}
\def\hatedot{{\hat E}^{\bullet}}

\def\hatsym{\hat{\text{Sym}^{\bullet}}}
\def\ce{C^{\edot}}
\def\vt{\text{Vect}}
\def\ao{\bold A^1}
\def\piao{\pi_{\ao}}
\def\qed{\ \ $\square$\enddemo}

\topmatter
 
\title
Virtual moduli cycles and Gromov-Witten invariants
of algebraic varieties
\endtitle
 
\author
Jun Li and Gang Tian
\endauthor
 
\address
Stanford University
\endaddress
 
\email
jli\@gauss.stanford.edu
\endemail
 
\thanks
Supported in part by NSF grant, A. Sloan fellowship and Terman
fellowship.
\endthanks
 
\address
Massachusetts Institute of technology
\endaddress
 
\email
tian\@math.mit.edu
\endemail
 
\thanks
Supported in part by NSF grant.
\endthanks
 


\endtopmatter

\document

\head Introduction
\endhead

The study of moduli spaces plays a
fundamental role in our understanding
the geometry and topology of manifolds.
One example is the Donaldson theory (and more recently
the Seiberg-Witten invariants), which provides a set of
differential invariants of 4-manifolds [Do]. When the 
underlying manifolds are smooth algebraic surfaces, then they are
the intersection theories on the moduli spaces
of vector bundles over these surfaces [Li, Mo].
Another example is the mathematical theory, inspired
by the sigma model theory in mathematical physics
([W1], [W2]), of quantum cohomology. The quantum cohomology
uses the GW-invariants, which are the intersection numbers of
certain induced homology classes on the moduli spaces of rational
curves in a given symplectic manifold. This is a generalization of
the classical enumerative invariant which counts
the number of algebraic curves with appropriate constraints in a variety.
The first mathematical foundation of quantum cohomology was established
by Ruan and the second named author in [RT1] for
semi-positive symplectic manifolds,
which include all algebraic manifolds of
complex dimension less than $4$,
all Fano manifolds and Calabi-Yau spaces.
In [RT2], general GW-invariants of higher
genus are constructed to establish a mathematical theory of
the sigma model theory coupled with gravity on any 
semi-positive symplectic manifolds (also see [Ru] for the special cases). 
There are some related works we would like to 
mention. In [KM], Kontsevich and Manin
proposed an axiomatic approach to GW-invariants, 
and in [Ko2], Kontsevich introduced the notion of stable maps to
study GW-invariants. There are also works dealing with special 
classes of Fano varieties, such as homogeneous manifolds. (cf. [BDW; Ber;
Ci; CM; LT].)

Now let us discuss the new issue in intersection theory raised from
studying GW-invariants, and more generally the Donaldson type invariants.
The core of an intersection theory is the fundamental class. 
For a manifold (or a variety), the ordinary cup product with 
the fundamental class given by the underlying manifold provides a
satisfactory intersection theory.
However, for the GW-invariants, which should be an
intersection theory on the moduli space of stable maps, we can not
take the fundamental class of the whole moduli space  
directly. This is because the relative
moduli space (i.e., the family version) in general do not form a flat family 
over the parameter space. One guiding principle of our search of
a ``good'' intersection theory is that such a theory should
be invariant under deformation of the
underlying manifolds. In [Do, RT1, RT2], they employed
analytic methods to construct ``good'' intersection theory
using generic moduli spaces (they are almost always non algebraic).

Abiding with algebraic methods, we have no luxury of having a ``generic
moduli space''. Instead, we will construct directly a cycle in
the moduli space, called the virtual moduli cycle, and define an intersection
theory by using this cycle as the fundamental class. Such 
a construction commutes with Gysin maps.
In this paper, we will construct such
a cycle by first constructing a 
cone cycle inside a vector bundle, which functions as a normal cone, 
and then intersecting this cone cycle with the 
zero section of the vector bundle. To make this construction
sufficiently general, we shall carry it out  
based on the moduli functor solely. The data we need is a choice of
tangent-obstruction complex of the moduli functor,
which is a global obstruction 
theory of the moduli
problem. The virtual moduli cycle depends on the choice of
such a complex, so does the virtual intersection theory defined.
The so constructed intersection theory will have 
the following invariance property.
Given a family of moduli functors, namely, a relative moduli functor,
if we assume that the tangent-obstruction complex of the 
relative moduli functor
and that of the specialized moduli functor are compatible, then the
specialization of the virtual intersection theory on the relative moduli
space is the same as the virtual intersection theory of the 
specialized moduli space. 
Applying to the moduli space of stable maps from $n$-pointed
nodal curves into a smooth projective variety $X$, we can define
the GW invariants of $X$ purely algebraically. 

We now describe briefly the key idea to our construction.
When we are working with a moduli space, usually we
can compute its virtual dimension. 
However, the virtual dimension may not coincide with
the actual dimension of the moduli space. One may view this as
if the moduli space is a subspace of an ``ambient'' space cut out
by a set of ``equations'' whose vanishing loci do not meet
properly. Such a situation is well understood in the following
setting: let
$$\CD Z @>>> X\\
      @VVV @VV f V\\
      Y @> g >> W
\endCD
$$
be a fiber square, where $X$, $Y$ and $W$ are
smooth varieties and subvarieties. Then $[X]\cdot [Y]$,
the intersection of the cycle $[X]$ and $[Y]$, is a
cycle in $A\lsta W$ of dimension $\dim X+\dim Y-\dim W$.
When $\dim Z=\dim X+\dim Y-\dim W$, then $[Z]=[X]\cdot [Y]$.
Otherwise, $[Z]$ may not be $ [X]\cdot[Y]$. The excess intersection
theory tells us that we can find a cycle in $A\lsta Z$ 
so that it is $[X]\cdot[Y]$. We may view this cycle as the virtual
cycle of $Z$ representing $[X]\cdot[Y]$. 
Following Fulton-MacPherson's normal cone
construction, this cycle is the image of
the cycle of the normal cone to $Z$\ in $X$, denoted
by $C_{Z/X}$, under the Gysin homomorphism 
$s\sta\mh A\lsta\bl C_{Y/W}\times
_YZ\br \to A\lsta Z$, where $s\mh Z\to C_{Y/W}\times_YZ$\ is the
zero section. This theory does not apply directly to moduli
schemes, since, except for some isolated cases, it is impossible
to find pairs $X\to W$\ and $Y\to W$\ so that $X\times_WY$\
is the moduli space and $[X]\cdot[Y]$\ so defined is the
virtual moduli cycle we need.

The strategy to our approach is that rather than
trying to find an embedding of the moduli space
into some ambient space, we will
construct a cone in a vector bundle directly, say $C\sub V$,
over the moduli space and then define the virtual moduli 
cycle to be $s\sta[C]$, where $s$ is the zero section of $V$.
The pair $C\sub V$ will be constructed
based on a choice of the tangent-obstruction complex
of the moduli functor (The definition of
tangent-obstruction complex is given in section 1).

Let $\m$ be a moduli space. We first construct
its tangent-obstruction complex, which usually comes from
studying obstruction
theory of the moduli problem. 
For a large class of moduli problems, 
their tangent-obstruction complexes
are the sheaf cohomologies of complexes of locally free 
sheaves
$$\edot=[\e_1\lra \e_2].
$$
Assume that $\m$ belongs to this class of moduli problems.
Then at each closed point $w\in\m$, $T_1=
\hh^1(\edot\otimes_{\o_{\m}} k(w))$ is the tangent space $T_w\m$ and
$T_2=\hh^2(\edot\otimes_{\o_{\m}}k(w))$ is the obstruction space to
deformations of $w$ in $\m$. 
Here $\hh^i(\edot)$ is the $i$-th sheaf cohomology of the complex.
There is an ``intrinsic''  set of 
defining equations of the germ of $\m$ at $w$, namely, the Kuranishi
map 
$$f: \hatsym(T_2\dual)\lra \hatsym(T_1\dual):=
\lim_{\longleftarrow}\oplus_{l=0}^n S^l(T_1\dual).
$$
Note that if we denote by $\hat w$ be the formal completion of
$\m$ along $w$, then [La]
$$\hat w\cong \spec \hatsym(T_1\dual)\otimes_{\hatsym(T_2\dual)}k.
$$
The normal cone to $\hat w$ in $\spec \hatsym(T_1\dual)$
is canonically a subcone in $\hat w\times_k T_2$.
We denote this cone by $C_w$.
The virtual normal cone we seek will be a cycle
$[C]$ in $Z\lsta\vt(\e_2)$, where $\vt(\e_2)$ is the vector bundle over
$\m$ so that its sheaf of sections is $\e_2$. The $[C]$ is
uniquely determined by the following criterion.
At each $w\in\m$, there is a surjective vector bundle homomorphism
$$\vt(\e_2)\times_{\m}\hat w\lra \vt(T_2)\times\hat w,
$$
where $T_2=\hh^2(\edot\otimes_{\o_{\m}}k(w))$,
that extends the given homomorphism
$\e_2\to\hh^2(\edot\otimes_{\o_{\m}}k(w))$ such that the 
restriction of $[C]$ to $\vt(\e_2)\times_{\m}\hat w$ is
the pull back of $C_w$. In short, the virtual normal cone
is the result of patching these local normal cones defined by the
Kuranishi maps of the moduli space.
The virtual moduli cycle $[\m]\vir$ is 
then defined to be the image 
of $[C]$ in $A\lsta\m$\ via the Gysin homomorphism 
$$s\sta:A\lsta \vt(\e_2)\to A\lsta {\Cal M}.
$$
The GW-invariants is defined by
applying this construction to the moduli spaces of 
stable morphisms from nodal curves to $X$.

\pro{Theorem}
For any smooth projective variety $X$, and any choice of integers $n$\ and
$g$\ and
$\alp\in A_1X/\sim_{\text{alg}}$, there is a virtual
moduli cycle
$[\mgny]\vir\in A_k X\otimes_{\ZZ}\QQ$,
where $k$ is the virtual dimension of $\mgna$. Using this cycle, we can
define the GW-invariants
$$\Psi^X_{\alpha,g,n}: (A\sta X)^{\times n}\times A\sta\mgn \lra
\aast\mgny\otimes_{\ZZ}\QQ\,
$$
in the usual way. Let
$\psi_{\alpha,g,n}^X$ be
the composite of $\Psi^X_{\alp,g,n}$\ with the degree homomorphism
$\aast\mgnx\to\QQ$. Then
$\psi^X_{\alpha,g,n}$ are invariant under deformations of $X$ and satisfy
all the expected properties of the GW-invariant, including the composition law.
\endpro

The construction of the virtual moduli cycle $[\mgny]\uvir$\
is the main purpose of this paper. The proof of the composition
law is almost straightforward, following a similar process as in [RT1].
Since maps in $\mgny$\ may have non-trivial automorphisms.
our classes may have rational coefficients.
The approach to this problem is
the usual descent argument. In the end, we obtain a
cycle supported on an effective cone over $\mgna$ inside a
$\QQ$-vector bundle. The virtual moduli cycle is 
then the image of the
Gysin map as described before. 
The resulting class is of rational coefficients.


This construction of virtual cycles
was finished in early 1995. 
During AMS summer school held at Santa Cruz, July, 1995,
the first named author reported this work.
In his talk, he described the ideas of our construction of virtual
moduli cycles and the definition of GW-invariants. After the talk, 
S. Katz kindly informed the first named author that he 
had studied the problem of constructing
virtual moduli cycles and obtained some partial results in special cases [Kz1].
During the preparation of this paper, we learned that K. Behrend and B.
Fantechi had gave an alternative construction of 
virtual moduli cycles [BF]. Also
as an application, K. Behrend defined GW-invariants and proved 
the basic property of these invariants [Bh].
A similar idea can be applied to constructing symplectic invariants
for general symplectic manifolds.

The layout of this paper is as follows. In section one, we 
introduce the notion of tangent-obstruction complexes of functors, which is a 
global obstruction theory of the moduli functors. We then describe the
tangent-obstruction complex of the 
moduli functor of stable morphisms.
The next two sections are devoted to
construct and investigate the virtual normal cone of 
any tangent-obstruction complex.
In section 4 and 5, we will construct the GW-invariants and
prove some basic properties of these invariants, including their
deformation invariance and the composition laws. 

The first named author thanks
W. Fulton and D. Gieseker for many stimulating discussions.
Part of this work
was done when the second named author
visited Department of Mathematics, Stanford University
in the winter quarter of 1994.
He would like to thank his colleagues there
for providing a stimulating atmosphere.
We thank the referee for many comments and suggestions.

\head 
1. Tangent--Obstruction complex
\endhead

In this section we will introduce the notion of 
tangent-obstruction complexes of moduli functors. 
Such a notion was implicit in many earlier works and
should be viewed as another way of presenting deformation 
theory.

In this paper, we will fix an algebraically closed field $k$ of
characteristic 0 and will only consider schemes over $k$.

We first define the functor of tangent spaces. Let $\s$ be the category of
all schemes and let $\FF\mh \s\to(\text{sets})$ be a (contravariant)
moduli functor. Here we call $\FF$ a moduli functor if for any
$S\in\s$ the object $\FF(S)$ is the set of {\it isomorphism classes} of
flat families of objects (to be parameterized) over $S$. For our
purpose, we will introduce an associated functor, called the
pre-moduli functor of $\FF$ and denoted by
$\FF^{\text{pre}}$. For any $S\in\s$,
$\FF^{\text{pre}}(S)$ is the set of all flat families of objects 
(to be parameterized) over $S$. Note that we
do not take isomorphism classes in this case. Following the convention, 
for $\xi_1,\xi_2\in\FF(S)$ we
will denote by $\xi_1\cong\xi_2$ if $\xi_1$ and $\xi_2$ are isomorphic
and denote by $\sim$ the equivalence relation induced by 
$\cong$. Hence $\FF(S)=\FF^{\text{pre}}(S)/\sim$. 
In this paper, for $S\in\s$ we let $\mods$ be the category of
sheaves of $\o_S$-modules. For $\n\in\mods$, we denote by $\sn$
the trivial extension of $S$ by the sheaf $\n$\ {}
\footnote{By this we mean
$\sn=\spec(\Gamma(\o_S)\ast\Gamma(\n))$, where $\Gamma(\o_S)\ast\Gamma(\n)$
is the trivial ring extension of $\Gamma(\o_S)$ by $\Gamma(\n)$. Note that
there is an inclusion $S\to S_N$ and projection $S_N\to S$
so that $S\to S_N\to S$ is the identity. (See [Ma,p191]).} 
and denote by 
$\pi_{\n}\mh\FF^{\text{pre}}(\sn)\to\FF^{\text{pre}}(S)$
the restriction morphism induced by the obvious inclusions $S\sub\sn$.
Given $\eta_0\in\FF^{\text{pre}}(S)$ and 
$\zeta_1,\zeta_2\in\pi_{\n}\upmo(\eta_0)$, we say 
$\zeta_1\cong_{\eta_0}\zeta_2$ if there is an isomorphism $\rho\mh\zeta_1\cong
\zeta_2$ so that its restriction $\pi_{\n}(\rho)\mh\eta_0\cong\eta_0$ 
is the identity isomorphism. Now we define the functor of tangent spaces.
For any affine $S\in\s$ and $\eta_0\in\FF^{\text{pre}}(S)$,
we let $\t\FF(\eta_0)\mh\mods\to(sets)^0$ be the functor
that assigns $\n$ to the set $\pi_{\n}\upmo(\eta_0)$ modulo the
equivalence relation induced by the isomorphism $\cong_{\eta_0}$.
In short, $\t\FF(\eta_0)$ consists of all isomorphism classes of
$\tilde \eta\in\FF^{\text{pre}}(\sn)$ whose restrictions to
$S$ is $\eta_0$. In case $\eta_0\cong\eta_0\pri$,
then there is a canonical isomorphism of sets $\t\FF(\eta_0)$ and
$\t\FF(\eta_0\pri)$. This way, for $\eta\in\FF(S)$ we obtain an
isomorphism class of functors $\f\FF(\eta_0)\mh\mods\to(sets)^0$.
It is clear that if $\rho\mh S_1\to S_2$ is a morphism between
affine schemes and if
$\n_1\in{\frak Mod}_{S_1}$ and $\n_2\in {\frak Mod}_{S_2}$
are two sheaves with $\o_{S_1}$-homomorphism 
$\o_{S_1}\otimes_{\o_{S_2}}\n_2\to\n_1$, then for any $\eta_2\in\FF(S)$ 
with $\eta_1=\FF(\rho)(\eta_2)$ the induced object in $\FF(S_1)$, there is
a canonical morphism 
$$\t\FF(\eta_2)(\n_2)\lra\t\FF(\eta_1)(\n_1),
$$
satisfying the base change property. Note that when
$\FF$ is represented by a scheme $Y$ and $\eta\in\FF(S)$ is
represented by a morphism $f\mh S\to Y$, then
$$\t\FF(\eta)(\n)=\Gamma(\hom_S(f\sta\Omega_Y,\n)).
$$
In this case, $\t\FF(\eta)$ is a functor $\mods\to\mods$.

\pro{Assumption}
In this paper, we will only consider moduli functor $\FF$ such that
$\t\FF$ is  induced by a
sheaf-valued functor over fibered category of modules over
schemes over $\FF$. Namely, for any affine $S\in\mods$,
$\eta\in\FF(S)$ and $\n\in\mods$ the set $\t\FF(\eta)(\n)$ is
canonically isomorphic to the set of all sections of a sheaf of
$\o_S$-modules, denoted by $\tone\FF(\eta)(\n)$,
and the arrows (above) in the base change property are 
induced by sheaf homomorphisms
$$\tone\FF(\eta_2)(\n_2)\otimes_{\o_{S_2}}\o_{S_1}
\lra\tone\FF(\eta_1)(\n_1).
$$
\endpro

In the following, we will call $\tone\FF$ the functor of tangent spaces.
We remark that we have not exhausted the literatures to see how
restrictive is this assumption. Nevertheless, the moduli functors that
will be discussed in this paper all satisfy this condition.

Next, we recall the definition of an 
obstruction theory. An obstruction theory to 
deformations of $p\in\FF(\spec k)$ with values in a vector space $O$ 
is an assignment as follows. Given a
pair $(\eta\in\FF(\spec B/I), I\sub B)$, where
$B$ is an Artin ring with the residue field $k$, $I\sub B$ 
is an ideal annihilated by the maximal ideal of $B$ 
and $\eta\otimes_{B/I} k=p$,
the obstruction theory assigns a natural obstruction class
$ob(\eta,B/I,B)\in O$ whose vanishing is the necessary and sufficient 
condition for $\eta$ to be extendible to $\tilde \eta\in\FF(\spec B)$. 
We now introduce its relative analogue.

\pro{Definition 1.1} Let $\Cal K=\{\k_{\eta}\}$ be a collection of 
sheaves of $\o_S$-modules $\k_{\eta}$ indexed by $\eta$ for 
$S\in\s$ and $\eta\in\FF(S)$. We say that $\k$ is a sheaf over $\FF$ 
if for any morphism $f\mh T\to S$
of schemes, there is an isomorphism 
$f\sta\k_{\eta}\cong\k_{f\sta\eta}$ canonical under base change.
\endpro

\pro{Definition 1.2}
An obstruction theory of the moduli functor $\FF$ with values in a sheaf 
$\ob_{\bullet}$ over $\FF$ 
consists of the following data: let $S$ be any affine scheme,
let $S\to Y_0\to Y$ be schemes and embedding morphisms over $S$. Namely,
$Y_0\to Y$ is an embedding and $i\mh S\to Y_0$ is a section of $Y_0\to S$.  
Let $\mm$ be the ideal sheaf of $S\sub Y$ and let $\Cal I\sub\o_Y$ be the 
ideal sheaf of $Y_0\sub Y$. Assume that $\Cal I\cdot\mm=0$. Then 
for any $\eta\in\FF(Y_0)$, there is an obstruction class
$$ob(\eta,Y_0,Y)\in\Gamma_S(\ob_{\eta_0}\otimes_{\o_S}\Cal I),
$$ 
where $\eta_0=i\sta(\eta)\in\FF(S)$,
whose vanishing is the necessary and sufficient condition for $\eta$
to be extendible to $\tilde\eta\in\FF(Y)$. We call $ob(\eta,Y_0,Y)$ the
obstruction class to extending $\eta$ to $Y$. The obstruction class
is canonical under base change: let $\sigma\mh S\pri\to S$ be another morphism, 
$Y\pri$ a scheme over $S\pri$ with a section $i\pri$ and $f\mh Y\pri\to Y$ 
a morphism such that
$$\CD Y\pri @>f>> Y\\ @VVV @VVV\\ S\pri @>>> S
\endCD
\qquad\text{and}\qquad
\CD Y\pri @>f>> Y\\ @A{i\pri}AA @AiAA\\ S\pri @>{\sigma}>> S
\endCD
$$
are commutative. Let $Y_0\pri=Y_0\times_Y Y\pri$ and let
$\eta\pri\in\FF(S\pri)$ be the pull back of $\eta$. Let $g\mh 
\sigma\sta(\ob_{\eta_0}\otimes_{\o_S}
\Cal I_{Y_0\sub Y})\to \ob_{\eta\pri_0}\otimes_{\o_{S\pri}}
\Cal I_{Y_0\pri\sub Y\pri}$
be the obvious homomorphism. Then
$$ob(\eta\pri,Y_0\pri,Y\pri)=g (ob(\eta,Y_0,Y)).
$$
\endpro

\noindent
{\bf Example}.
Let $X\sub \bold A^n$ be a subscheme defined by the ideal $I=(f_1,\ldots,f_m)$.
Let $\FF_X$ be the functor ${\bold Mor}(-,X)$
and let $\cbul$ be the complex $[\o_X(T\bold A^n)\mapright{\sigma}\o_X^{\oplus
m}]$ where $\sigma=(df_1,\ldots,df_m)$. Then for any affine $S$, morphism
$\eta\mh S\to X$ and sheaf $\n\in\mods$, $\tone\FF_X(\eta)(\n)$ is the first
sheaf cohomology of $\eta\sta\cbul\otimes_{\o_S}\n$. The defining sections
$f_1,\ldots,f_m$ defines an obstruction theory of the functor $\FF_X$
with values in $\coker(\sigma)$ (See section 2 for explicit description).

\noindent
{\bf Example} ([Al]).
This example concerns the moduli of {\it stable} sheaves $\e$ on a
smooth algebraic surface $X$ of fixed Poincare polynomial $\chi$.
Here we implicitly fix an ample divisor on $X$. We denote 
the corresponding moduli functor by $\FF_{\chi}$. For any affine $S\in\s$
and $\eta\in\FF_{\chi}(S)$ representing the sheaf $\e$ of 
$\o_{X\times S}$-modules, then 
$$\tone\FF_{\chi}(\eta)(\n)=\ext_{S\times X/S}^1(\e,\e\otimes\pi_S\sta\n)^0,
$$
where the superscript means the traceless part
of the extension sheaf. The canonical obstruction theory of 
$\FF_{\chi}$ takes values in the sheaf $\ext_{S\times X/S}^2(\e,\e)^0$.

Let $\text{ob}_{\bullet}$ be the sheaf in Definition 1.2. 
For simplicity, we will use the convention $\ttwo\FF(\eta)(\n)=
\text{ob}_{\eta}\otimes_{\o_S}\n$ and $\tdot\FF=[\tone\FF\to
\ttwo\FF]$, where the arrow is the zero homomorphism. 
By the assumption of this section,
$\tone\FF(\eta)(\n)$ is a two-term complex of
sheaves of $\o_S$-modules connected by the zero arrow.

\pro{Definition 1.3}
Let $\FF$ be as before. A tangent-obstruction complex
of $\FF$ is a complex $\tdot\FF=[\tone\FF\to\ttwo\FF]$,
where the arrow is the zero arrow,
such that $\tone\FF$ is the functor of the tangent spaces 
of $\FF$ and that there is an obstruction theory of $\FF$
taking values in $\ttwo\FF$. The tangent-obstruction complex
$\tdot\FF$ is said to be a perfect
tangent-obstruction complex if for any affine $S$ and $\eta\in
\FF(S)$, there is an affine covering $\{ S_{\alpha}\}$ of $S$
such that there are two-term complexes of locally free sheaves 
$\edot_{\alpha}$ such that for any $\n\in\frak{Mod}_{S\la}$,
$\t^i\FF(\eta\la)(\n)$ is the $i$-th sheaf cohomology of 
$\edot_{\eta\la}\otimes_{\o_{S\la}}\n$,
where $\eta\la\in\FF(S\la)$ is the induced object of $\eta$
via $S\la\to S$. In case the complex $\edot$ is
explicitly given, we will write $\tdot\FF=\hhdot(\edot)$.
\endpro

We emphasize that the obstruction theory is part of the data
making up the tangent-obstruction complex. This notion is a
convenient way to group the data of tangent spaces and obstruction
theory. It will be clear later that when the tangent-obstruction
complex is perfect then the classical construction of Kuranishi
maps can be adopted to construct the relative Kuranishi families,
which is the heart of the construction of virtual moduli cycles.

\noindent
{\bf Remark.}
The assumption that $\FF$ admits a perfect tangent-obstruction
complex is a strong requirement. For instance, moduli functors of
stable vector bundles over threefolds other than Calabi-Yau manifolds may
not have perfect tangent-obstruction complexes.

In this paper, our main interest is in the moduli spaces of stable
morphisms from marked curves to
smooth projective varieties. Let $X$ be a fixed smooth project
variety. We first recall the notion of 
stable morphisms introduced by Kontsevich [Ko1]. 
An $n$-pointed nodal curve is a nodal curve $C$ and $n$
ordered marked points $D \subset C$ away from the singular
locus of $C$ (We will use $D$ to denote the $n$-ordered marked points
on $C$ in this paper). A morphism $f: D\sub C \to X$ is said to be stable
if $D\sub C$\ is an $n$-pointed connected nodal curve and
$f\mh C\to X$\ is a morphism such that 
$$
\Hom_C (\Omega_C (D), \Cal{O}_C) \to \Hom_C (f^* \Omega_X, \Cal{O}_C)
$$
is surjective, where
$f^* \Omega_X \to \Omega_C (D)$ is induced by $f^* \Omega_X \to
\Omega_C$.  
We will call $f$ stable relative to $D$ or simply stable
if the marked points $D\sub C$\ are understood.

>From now on, we fix a class $\alpha \in A_1 X /\simalg$ 
and two integers $n$ and $g$.
We let $\fgna\mh\s\to(sets)^0$ be the functor that
assigns any $S\in\s$ to the set of
all isomorphism classes of flat families over $S$ of stable morphisms
$$f:D\sub \x\lra X
$$
from $n$-pointed connected nodal curves
$ D\sub \x$ of arithmetic genus $g$ to $X$ such that 
$f $\ sends the fundamental classes of closed fibers of 
$\x$ over $S$\ to $\alpha$.
Since $X$ is a smooth projective
variety, by the work of [Al], $\fgna$\ is coarsely represented by
a projective scheme. We denote this scheme by $\mgna$.  
It is also known that $\fgna$ is represented by a Deligne-Mumford
stack [FP].

In the following, we will determine the natural
tangent-obstruction complex of $\fgna$.  
We fix an affine scheme $S$ and a sheaf $\n\in\mods$. Let
$\xi\in\fgna(S)$ be represented by
$f\mh \x\to X$ with marked sections $D\sub\x$.
Let $\Cal{X}^N$ be a flat family of nodal
curves over $S_N$, where $S_N$ is the trivial extension of $S$ by $\n$,
that extends the family $\x$.
Then we have a commutative diagram of exact sequences
$$
\CD
0 @>>> \Cal{O}_{\x} \otimes_{\o_S} \n          @>>>    
       \Omega_{\x^N/S} \otimes_{\o_{\x^N}} \o_{\x}   @>>> 
       \Omega_{\x/S}                       @>>>     0  \\ 
@.     @|     @AA{d_N}A     @AAdA                           \\
0 @>>> \Cal{O}_{\x} \otimes_\os \n          @>>>    
       \o_{\x^N}     @>>> 
       \o_{\x}                       @>>>    \ \, 0  
\; ,
\endCD
$$
where $d$ and $d_N$ are the differentials.
The upper sequence is exact because $\x^N$ is a family of
nodal curves flat over $S_N$, following Theorem 25.2 in [Ma].
Conversely, given any exact sequence of sheaves of 
$\Cal{O}_{\x}$-modules
$$
\CD
0 @>>> \o_{\xs} \otimes_\os \n @>>> \cc @>>> 
\Omega_{\Cal{X}/S} @>>> 0\,,
\endCD
$$
we obtain a pull-back from the $S$-homomorphism $d :
\Cal{O}_{\Cal{X}} \to \Omega_{\Cal{X}/S}$:
$$
\CD
0 @>>> \Cal{O}_{\Cal{X}} \otimes_S \n   @>>>   \cc  @>>> 
      \Omega_{\Cal{X}/S}                @>>>   0           \\ 
@.     @|     @AAA        @A{d}AA                              \\
0 @>>> \Cal{O}_{\Cal{X}} \otimes_S \n   @>>>   \Cal B   @>>> 
       \Cal{O}_{\Cal{X}}               @>>>     \ \, 0  \; .
\endCD
$$
One checks that there is a canonical way to give $\Cal B$ a structure
of sheaf of $\o_{S_N}$-algebras,
which is flat over $S_N$ automatically.  Thus we
obtain a flat family $\Cal{X}^N$ over $S_N$. One checks also that this 
correspondence is one-to-one and onto. This is the 
one-one correspondence between the
space of flat extensions of $\xs\to S$\ to $\xsn\to\sn$
and the module
$$
\Ext^1_{\xs/S}
\left(
  \Omega_{\Cal{X}/S}, \Cal{O}_{\Cal{X}} \otimes_{\o_S} \n
\right)  \; .
$$

Next we investigate when such an extension $\Cal{X}^N$
admits a morphism $f^N\mh \Cal{X}^N \to X$ extending $f\mh 
\x \to X$.  We claim that such an $f ^N$ comes from the
existence of an $\o_{\x}$-linear lifting $f^* \Omega_X \to
\Omega_{\Cal{X}^N/S}\otimes_{\o_{\xsn}}\o_{\x}$ 
of the obvious $f \sta\om_X\to\om_{\x/S}$.
Indeed, given $f^N$ restricting to $f $, we certainly have 
such a lifting from
$$(f ^N)\sta\om_X\lra \om_{\x^N/S}
\quad \text{and}\quad
(f^N)\sta \Omega_X\otimes_{\o_{\x^N}}\o_{\x}=\fs\sta\om_X.
$$
Conversely, given any diagram
$$
\CD
@. @.  f^* \Omega_X   @=     f^* \Omega_X       \\ 
@. @. @V\beta VV         @VVV       \\
0  @>>>  \Cal{O}_{\Cal{X}} \otimes_S \n@>>>\Cal B@>>> 
            \Omega_{\Cal{X}/S  }              @>>>     0
\; ,
\endCD
\tag 1.1
$$
we first obtain a flat extension $\xsn$\ of $\xs$ and 
an isomorphism $\Cal B\cong\om_{\xsn/S}\otimes_{\o_{\x^N}}\o_{\x}$  
based on the bottom exact sequence.
Observing that 
$$\CD
\o_X@>{\gamma}>> f\lsta(\o_{\x})\\
@VV{f\lsta(\beta)\circ d}V @VV{f\lsta(d)}V\\
f\lsta(\Omega_{\Cal{X}^N/S} \otimes_{\o_{\x^N}}\o_{\x}) 
 @>>>\,   f\lsta( \Omega_{\Cal{X}/S} )
\endCD
$$
is commutative, we can factor $f\lsta(\beta)\circ d$ and $\gamma$ through
$\o_X\to f\lsta(\o_{\x^N})$, because
$$
\CD
f\lsta(\Cal{O}_{\Cal{X}^N}  )   @>>>  f\lsta( \Cal{O}_{\Cal{X}}   )  \\
@VV{f\lsta(d_N)}V                             @VV{f\lsta(d)}V                      \\
f\lsta(\Omega_{\Cal{X}^N/S} \otimes_{\o_{\x^N}}\o_{\x}) 
 @>>>\,   f\lsta( \Omega_{\Cal{X}/S} )
\endCD
$$
is a pull-back diagram.
One checks directly that $\o_X\to f\lsta(\o_{\x^N})$ is a
homomorphism of sheaves of $S$-algebras. Therefore, it
defines a morphism $f ^N: \xsn\to X$
that is an extension of $f\mh\x\to X$.

In conclusion, we have shown that for any affine
$S\in\s$ and $\xi\in\FF_{\alpha,g,0}^X(S)$ that corresponds to the family
$f: \xs \to X$, the tangent $\t\FF_{\alpha,g,0}^X$
at $\xi$\ takes $\n\in\mods$ to the set of all commutative diagrams of
$\o_{\x}$-modules (1.1).
This set is naturally the first extension module
$$
\Ext^1_{\xs}
\left(
  [f^* \Omega_X \to \Omega_{\Cal{X}/S}], 
  \Cal{O}_{\Cal{X}} \otimes_{\os}\n
\right) \,,
$$
where $[f^* \Omega_X \to \Omega_{\Cal{X}/S}]$ is a complex 
indexed at $-1$ and 0. 
We should point out that in
[Ra], Z. Ran has identified the deformation space and the obstruction space of
this moduli problem to the above diagram and has expressed them in
terms of extension modules over non-commutative rings. He actually treated 
more general cases. Based on [Ra] and the above 
reasoning, we only need to check that the set of diagrams above is 
canonically isomorphic to the above extension module. 
This can be checked by using hypercohomology of a 
double complex based on a covering of
$\Cal{X}_S$ and a locally free resolution of
$\Omega_{\Cal{X}_S}$, similar to the example in [GH].  We leave
the details to readers.

Now we give the tangent of the functor $\fgna$.

\proclaim {Proposition 1.4}
Let $S$ be any affine scheme and let
$\xi\in\fgna(S)$\ be represented by
$f \mh\x \to X$ with marked points $D\sub\x$.
Then for any $\n\in\mods$,
$$
\tone\fgna(\xi)(\n)=
\ext^1_{\xs/S}
\bl
  [f^* \Omega_X \to \Omega_{\Cal{X}/S} (D)], 
  \Cal{O}_{\Cal{X}} \otimes_\os \n
\br \,.
$$
\endproclaim

\demo {Proof}
We will sketch the proof of one direction
and leave the other to the readers.
Given any section in the above sheaf, we can associate to it a diagram
$$
\CD
@.    @.    f^* \Omega_X  @=  f^* \Omega_X        \\
@.    @.    @VVV                     @VVV    \\
0     @>>>  \Cal{O}_{\Cal{X}} \otimes_\os \n  @>>> 
      \Cal{A}    @>>>    \Omega_{\Cal{X}/S} (D)   @>>>  0    \\
@.    @VVV       @VVV                     @|      \\
0     @>>>  \Cal{O}_{\Cal{X}} (D)\otimes_\os \n  @>>> 
      \Cal B    @>>>    \Omega_{\Cal{X}/S} (D)   @>>>  \ \, 0
\, ,
\endCD
$$
where the lower left square is the push-forward of sheaves of
$\o_{\x}$-modules. The last line (tensored by
$\Cal{O}_{\Cal{X}} (-D)$) defines an extension $\xsn$.
Since $f^* \Omega_X \to \Omega_{\Cal{X}/S}
(D)$ factors through $\Omega_{\Cal{X}/S}\sub
\Omega_{\Cal{X}/S}(D)$, 
$f^* \Omega_X \to \Cal B $ factors through 
$\Cal B(-D) \subset \Cal B $, and thus defines a morphism $f^N
: \Cal{X}^N \to X$.  The immersion $D^N
\to \Cal{X}^N$ extending $D \to \x$ is determined by the
data $\coker\{ \Cal{A}\to\Cal B \}$.  In this way, we have constructed an
extension 
$$f^N : D^N\sub \Cal{X}^N\lra X
\quad\text{of}\quad
f  : D\sub \x\lra X.
$$
It is routine to check that this
correspondence is one-to-one and onto, and satisfies the
required base change property. This proves the proposition.
\qed

We now describe the standard choice of the obstruction theory of $\fgna$.

\pro{Proposition 1.5}
For any $S\in\s$ and $\eta\in\fgna(S)$ corresponding to the family
$f\mh \x\to X$ over $S$ with the marked sections $D$, 
we define $\ob_{\eta}$ to be the sheaf 
$$\ext^2_{\x/S}([f\sta\Omega_X\to\Omega_{\x/S}(D)],\o_{\x}).
$$
We let $\ob_{\bullet}$ be the collection $\{\ob_{\eta}\}$ indexed by 
$\eta\in\FF(S)$ for 
$S\in\s$. Then $\text{ob}_{\bullet}$ 
forms a sheaf over $\fgna$. Furthermore, there is an obstruction
theory of $\fgna$ taking values in the sheaf $\ob_{\bullet}$.
\endpro

\demo {Proof}
It is clear that $\{\text{ob}_{\eta}\}$ is a sheaf over $\fgna$. 
Now we describe the 
obstruction theory taking values in this sheaf. Let $S$ be any affine scheme 
and let $S\to Y_0\to Y$ be a tuple of $S$-schemes described in Definition 1.2 
(Namely, $S\to Y_0\to Y$ are embeddings and the ideal sheaf of $Y_0\sub Y$
is annihilated by the ideal sheaf of $S\sub Y$).
Let $\eta\in\fgna(Y_0)$ be any object corresponding to a 
family $f_0\mh\x_0\to X$ 
with marked points $D_0\sub\x_0$ understood. We let 
$\bar\x=\x_0\times_{Y_0}S$,
$\bar D=D_0\times _{Y_0}S$ and $\bar f=f_0|_{\bar\x}$. We let $\bar\pi\mh\bar\x\to S$
be the projection. Since $\i:=\i_{Y_0\sub Y}$ is annihilated by the ideal 
sheaf $\i_{S\sub Y}$, $\i$ is a sheaf of $\o_S$-modules. 
Clearly, $D_0\sub\x_0$ can
be extended to a family over $Y$, say $D\sub \x$. 
Since $\x\to Y$ is a flat family of nodal curves, we obtain an exact sequence
$$0\lra\i_{\x_0\sub\x}\lra\Omega_{\x/S}\otimes_{\o_{\x}}\o_{\x_0}
\lra\Omega_{\x_0/S}\lra0.
\tag 1.2
$$
Let $\tau(\x)\in\Hom(f_0\sta\Omega_X,\Omega_{\x_0/S})$ be the obvious 
homomorphism and let 
$$\bar\tau(\x)\in\Ext^1_{\bar\x_0}({\bar f}\sta\Omega_X,{\bar \pi}\sta\i)
$$
be the image of $\tau(\x)$ under the connecting homomorphism
$$\Hom(f_0\sta\Omega_X,\Omega_{\x_0/S})\mapright{\delta}
\Ext^1_{\x_0}(f_0\sta\Omega_X,\i_{\x_0\sub\x})\
=\Ext^1_{\bar\x}({\bar f}\sta\Omega_X,{\bar\pi}\sta\i).
$$
Here the second identity holds because $f_0\sta\Omega_X$ is locally free.
It follows that $\bar\tau(\x)=0$ if and only if 
$f_0\sta\Omega_X\to\Omega_{\x_0/S}$ lifts to
$f_0\sta\Omega_X\to \Omega_{\x/S}\otimes_{\o_{\x}}\o_{\x_0}$.
This can be shown by similar arguments in studying $\t\fgna(\eta)$ 
that it
is the necessary and sufficient condition for $f_0\mh \x_0\to X$ to
be extendible to $f\mh\x\to X$. We let
$$ob(\eta,Y_0,Y)\in\Ext^2_{\bar\x}([{\bar f}\sta\Omega_X
\to\Omega_{\bar\x/S}(\bar D)],{\bar\pi}\sta\i)
$$
be the image of $\bar\tau(\x)$ under the obvious homomorphism
$$\Ext^1_{\bar\x}({\bar f}\sta\Omega_X,{\bar\pi}\sta\i)
\to 
\Ext^2_{\bar\x}([{\bar f}\sta\Omega_X
\to\Omega_{\bar\x/S}],{\bar\pi}\sta\i)\mapto{\cong}
\Ext^2_{\bar\x}([{\bar f}\sta\Omega_X
\to\Omega_{\bar\x/S}(\bar D)],{\bar\pi}\sta\i).
$$

To complete the proof, we need to check that
the definition of $ob(\eta_0,Y_0,Y)$ is
independent of the choice of extension $D\sub \x$, 
$ob(\eta_0,Y_0,Y)$ has the required
base change property and is the obstruction to extending $f$ to $Y$.
Since the choice of the
marked points of the nodal curve is irrelevant to extending
$f_0$ to $f$ and since the definition of $ob(\eta_0,Y_0,Y)$ is independent
of the choice of the marked points, 
to study the obstruction problem, we suffice to look at the situation where
$D=\emptyset$. We will assume this in the rest of this section.
We now check that $ob(\eta_0,Y_0,Y)$
is independent of the choice of the extension $\x$.
Indeed, let $\x\pri$ be another extension over $Y$. Then by the deformation
theory of nodal curves, there is an extension class 
$v\in \Ext^1_{\bar\x}(\Omega_{\bar\x/S},{\bar\pi}\sta\i)$ 
defining the exact sequence
$$0\lra {\bar\pi}\sta\i \lra\Cal A\lra\Omega_{\bar\x/S}\lra 0,
\tag 1.3
$$
of which the following holds. Let 
$$0\lra ({\bar\pi}\sta\i)^{\oplus 2} \lra
\Omega_{\x/S}\otimes_{\o_{\x}}\o_{\x_0}\oplus \Cal A\pri
\lra\Omega_{\x_0/S}\lra 0
\tag 1.4
$$
be the exact sequence induced by (1.2), (1.3) and
$\Omega_{\x_0/S}\to\Omega_{\bar\x/S}$. Then the bottom
exact sequence in 
$$\CD
0@>>> ({\bar\pi}\sta\i)^{\oplus 2} @>>>
\Omega_{\x/S}\otimes_{\o_{\x}}\o_{\x_0}\oplus \Cal A\pri
@>>>\Omega_{\x_0/S}@>>> 0\\
@.@V{(1,1)}VV @VVV @|\\
0@>>> {\bar\pi}\sta\i @>>>
\Cal B @>>>\Omega_{\x_0/S}@>>> 0,\\
\endCD
\tag 1.5
$$
where the left square is the push forward of sheaves,
is isomorphic to the exact sequence
$$0\lra {\bar\pi}\sta\i \lra
\Omega_{\x\pri/S}\otimes_{\o_{\x\pri}}\o_{\x_0}
\lra\Omega_{\x_0/S}\lra 0.
\tag 1.6
$$
Therefore, $\bar\tau(\x\pri)=\bar\tau(\x)+\delta_0(v)$, where
$$\delta_0\mh \Hom({\bar f}\sta\Omega_X,\Omega_{\bar\x/S})\to \Ext^1(
{\bar f}\sta\Omega_X,{\bar\pi}\sta\i)
$$ 
is the obvious connecting homomorphism, and hence the images
of $\bar\tau(\x)$ and $\bar\tau(\x\pri)$ in
$\Ext^2_{\bar\x}([{\bar f}\sta\Omega_X\to\Omega_{\bar\x/S}],{\bar\pi}\sta\i)
$
coincide. This proves that $ob(\eta_0,Y_0,Y)$ is well-defined.
For the same reason, the class $ob(\cdot,\cdot,\cdot)$ satisfies
the required base change property.

It remains to show that $ob(\eta_0,Y_0,Y)$
is the obstruction to extending $\eta$ to $Y$.
Obviously, if $\eta_0$ can be extended to $\eta\in\fgna(Y)$,
say $f\mh \x\to X$, we can take $\x$ to be the extension
of $\x_0$ and then $\bar\tau(\x)=0$ by construction.
Hence $ob(\eta_0,Y_0,Y)=0$. Now assume $ob(\eta_0,Y_0,Y)=0$. 
Because of the exact sequence
$$\Ext^1_{\bar\x}(\Omega_{\bar\x/S}, {\bar\pi}\sta\i)
\mapright{\beta}
\Ext^1_{\bar\x}({\bar f}\sta\Omega_X,{\bar\pi}\sta\i)
\lra
\Ext^2_{\bar\x}([{\bar f}\sta\Omega_X\to\Omega_{\bar\x/S}],{\bar\pi}\sta\i)
\lra 0,
$$
$\bar\tau(\x)$ is $\beta(-v)$ for some $v\in \Ext^1_{\bar\x}(\Omega_{\bar\x/S},
{\bar\pi}\sta\i)$. It follows from the deformation theory of nodal curves
that we can find an extension $\x\pri$ over $Y$ (of $\x_0$) such that the
diagrams of exact sequence (1.3)-(1.6) hold. Hence $\bar\tau(\x\pri)
=\bar\tau(\x)+\beta(v)=0$, which implies that $f$ extends to $f\pri\mh\x\pri\to X$.
This proves that $ob(\eta_0,Y_0,T)$ is the obstruction class to extending
$\eta_0$ to $Y$.
\qed

In section 4, we will show that $\tdot\fgna$ is a perfect 
tangent-obstruction complex of $\fgna$.


\head
2. Relative Kuranishi Families
\endhead

In this section, we will construct the relative Kuranishi families
of a perfect tangent-obstruction complex. We will show
that any two such families are equivalent under an explicit transformation.
This will be used to construct virtual normal cones and cycles of
moduli spaces in the next section.

We begin with the notion of relative tangent-obstruction complex and
the observation on how defining equation induces relative tangent-obstruction
complex. In this section, we assume
that $S$ is an affine scheme and $Z$ is a formal $S$-scheme
with a section $i\mh S\to Z$ so that as sets
$\text{Supp}(Z)=\text{Supp}(i(S))$.  
We further assume that there is a finite rank locally free sheaf $\f$ of
$\o_S$-modules such that $Z$ is embedded in $\spec\hatsym(\f)$
(Recall $\hatsym(\f)=\lim_{\leftarrow}\oplus_{l=0}^n S^n(\f)$).

\pro{Definition 2.1} Let $Z/S$ be as before.
A perfect relative tangent-obstruction complex is a two term complex
$[\e_1\mapto{\sigma}\e_2]$ of locally free sheaves of $\o_S$-modules 
of which the following holds.
\roster
\item
The cokernel of $\e_2\dual\mapto{\sigma\dual}\e_1\dual$ is isomorphic
to $\Omega_{Z/S}\otimes_{\o_Z}\o_S$.
\item
Let $O=\coker{\sigma}$. Then there is a relative obstruction 
theory to extending $S$-morphisms to $Z$ with values in $O$.
\endroster
Here, an obstruction theory to extending $S$-morphisms
to $Z$ is an assignment that, to each tuple of
$S$-schemes $S\to Y_0\to Y$ described in Definition 1.2 
and any $S$-morphism $\varphi_0\mh Y_0\to Z$, 
assigns a canonical obstruction class
$$ob(\varphi_0,Y_0,Y)\in\Gamma_S(O\otimes_{\o_S}\i_{Y_0\sub Y})
$$ 
to extending $\varphi_0$ to $Y\to Z$. 
\endpro

In this paper, 
for a scheme $W$, we often need to consider the formal completion of 
$W\times W$ along its diagonal. We denote this completion by 
$\hat W$. We will view $\hat W$ as a $W$-scheme
where $\pi\mh\hat W\to W$ is induced by the first projection of
$W\times W$. We will denote by $p_W\mh\hat W\to W$ the projection
induced by the second projection of $W\times W$. Note that there is a
canonical section $W\to\hat W$ of $\pi\mh\hat W\to W$
induced by the diagonal embedding $W\to W\times W$.

\pro{Lemma 2.2}
Let $W$ be a quasi-projective scheme.
Assume that $W$ admits a perfect tangent-obstruction complex 
$\tdot_W=\hhdot(\ecomplex)$. Let $S\sub W$ be
a locally closed subscheme. Then $\tdot_W$ canonically induces a 
perfect relative tangent-obstruction complex of $\hat W\times_W S/S$,
denoted by $\tdot_{\hat W\times_WS/S}$.
\endpro

\demo {Proof}
Since $\hat W$ is the formal completion of $W\times W$ along its diagonal,
$$\Omega_{\hat W\times_W S/S}\otimes_{\o_{\hat W\times_W S}}\o_S\cong
\Omega_W\otimes_{\o_W}\o_S.
$$
Now let $S\to Y_0\to Y$ 
be a tuple of $S$-schemes as before. Assume that 
$\varphi_0\mh Y_0\to \hat W\times_WS$ is an 
$S$-morphism. Then $p_W \circ\varphi_0$ is a morphism from $Y_0$ to $W$.
Clearly, $\varphi_0$ extends to an $S$-morphism 
$\varphi\mh Y\to\hat W\times_WS$ if and only if 
$p_W\circ\varphi_0$ extends to $Y\to W$, 
which is possible if and only if the obstruction
$$ob(p_W\circ\varphi_0,Y_0,Y)\in\Gamma_S(\hh^2(\edot\otimes_{\o_W}\o_S)
\otimes_{\o_S}\i_{Y_0\sub Y})
$$
vanishes. Hence $ob(p_W\circ\varphi_0,Y_0,Y)$ is the obstruction to extending
$\varphi_0$ to $Y$. 
\qed

To relate a Kuranishi family to an obstruction theory, 
we need to investigate how
defining equations induce a perfect relative tangent-obstruction complex.
Before we proceed, let us introduce the convention that will be used 
throughout this section. In this section, $S$ will always be an affine
scheme. Let $\e_1$ and $\e_2$ be two locally free sheaves of $\o_S$-modules.
We will assume throughout this section that $\Gamma(\e_i)$ are
free $\Gamma(\o_S)$-modules.
We will denote by  $A$ the ring $\Gamma(\o_S)$ and by
$E_i$ the free $A$-module $\Gamma(\e_i)$.
Given an $A$-module $N$, we will denote by $\hat{\text{Sym}^{\bullet}}(N)$
the inverse limit $\lim_{\leftarrow}\oplus_{l=0}^n S^l(N)$ of the direct
sum of the symmetric products of $N$. In this section, we will always use
$M$ to denote $\hat{\text{Sym}^{\bullet}}(E_1\dual)$. We denote by 
$M_1\sub M$ be the ideal generated by $E_1\dual\sub M$ and denote by $M_k$
the ideal $M_1^k$. For any $A$-homomorphism
$F\mh E_2\dual\to M$,
sometimes denoted $F\in M\otimes_A E_2$,
we will use $(F)$ to denote the ideal of $M$ generated by the
components of $F$. We now fix an $F\mh E_2\dual\to M$. 
We assume $(F)\sub M_1$. Let $\sigma\mh E_1\to E_2$ be the
dual of $E_2\dual\to M_1/M_2\equiv E_1\dual$
which is induced by $F$. We let $O=\coker(\sigma)$.

We now describe how $F$ induces a relative tangent-obstruction
complex to deformations of $S$-morphisms to $Z$.
Let $S\to Y_0\to Y$ be a tuple of $S$-schemes
as before and let $\varphi_0
\mh M/(F)\to\Gamma(\o_{Y_0})$ be an $A$-homomorphism. 
Let $g\mh M\to\Gamma(\o_Y)$ be a lift of $M\to M/(F)\to\Gamma(\o_{Y_0})$.
Clearly, $E_2\dual\mapto{F} M\mapto{g} \Gamma(\o_Y)$ 
factors through $o\mh E_2\dual\to\Gamma(\i_{Y_0\sub Y})$. Let
$$ob(\varphi_0,Y_0,Y)\in O\otimes_A\Gamma(\i_{Y_0\sub Y})
$$ 
be the image of $o$ under the obvious
$E_2\otimes_A\Gamma(\i_{Y_0\sub Y})\to O\otimes_A\Gamma(\i_{Y_0\sub Y})$.
We claim that this
is the obstruction to extending $Y_0\to Z$ to $Y\to Z$.
Assume $ob(\varphi_0,Y_0,Y)=0$. Then $o$ lifts to
an $h\mh E_1\dual\to\Gamma(\i_{Y_0\sub Y})$. Let $\hat h\mh M\to\Gamma
(\o_Y)$ be the induced homomorphism. It follows that
$g-\hat h\mh M\to\Gamma(\o_Y)$ factors through
$M/(F)\to\Gamma(\o_Y)$. Thus
$\varphi_0$ extends. The other direction is clear. 
We leave it to readers to check that such an
assignment of obstruction
class is canonical under base change.

\pro{Definition 2.3}
Let $Z/S$ be as before. Assume that $\tdot_{Z/S}=\hhdot(\ecomplex)$
is a perfect relative tangent-obstruction
complex of $Z/S$. A relative Kuranishi family
of $\tdot_{Z/S}=\hhdot(\ecomplex)$ is a pair
$(F,\Phi)$, where
$$F: E_2\dual\lra M
\qquad\and\qquad
\Phi:\spec_A M/(F)\lra Z,
$$
of which the following holds.
\roster
\item
$\Phi$ is an $S$-isomorphism.
\item
The complex $\ecomplex$ is identical to the sheafification
 of the complex $E_1\mapright{\sigma} E_2$ induced by $F$.
\item
The induced relative tangent-obstruction complex (from $F$)
is identical to the relative tangent-obstruction complex
$\tdot_{Z/S}=\hhdot(\ecomplex)$.
\endroster
\endpro

If the choice of the complex $\tdot_{Z/S}=\hhdot(\edot)$ is understood
from the context, we will simply call $(F,\Phi)$ a Kuranishi family and
call $F$ a Kuranishi map.

The relative Kuranishi families of $\tdot_{Z/S}=\hhdot(\ecomplex)$,
if exist, are not unique. Let 
$(\xi,\eta)\in\operatorname{Aut}_A(M)\times\homem$ be a pair such that
$$\xi\equiv 1_M\mod M_2\qquad\and\qquad
\eta\equiv\ide\mod M_1.
\tag 2.1
$$
We will show momentarily that if $(F,\Phi)$ is a relative 
Kuranishi family, then the pair
$(F\pri,\Phi\pri)$ defined by
$$F\pri=((1_M\otimes\eta)\circ(\xi\otimes\ide))(F)
\qquad\and\qquad
\Phi\pri=\Phi\circ\bar\xi
$$
is also a relative Kuranishi family.
Here, $\xi\otimes\ide$ and $1_M\otimes\eta$ are maps from 
$M\otimes_AE_2$
to $M\otimes_AE_2$ and 
$\bar\xi$ is the induced morphism
$$\spec_A M/(F\pri)\mapright{\bar\xi}\spec_A M/(\eta(F))=\spec_AM/(F).
$$
We will denote the pair $(F\pri,\Phi\pri)$ above by $(\xi,\eta)(F,\Phi)$.
We will call those $(\xi,\eta)$ satisfying (2.1) transformations. 
Given two transformations
$(\xi,\eta)$ and $(\xi\pri,\eta\pri)$, we define
$$(\xi,\eta)\cdot(\xi\pri,\eta\pri)=(\xi\circ\xi\pri,(1_M\otimes\eta)
\circ(\xi\otimes\ide)\circ(1_M\otimes\eta\pri)).
$$ 
It follows that 
$$\bl(\xi,\eta)\cdot(\xi\pri,\eta\pri)\br(F,\Phi)=
(\xi,\eta)\bl (\xi\pri,\eta\pri)(F,\Phi)\br.
$$
Let $\Cal K$ be the set of all relative Kuranishi families 
and let $\Cal H$ be the set of all transformations. 
It follows that $\Cal H$ is a group acting on $\Cal K$.

\proclaim{Proposition 2.4}
Let $Z/S$ be as before and let $\tdot_{Z/S}=\hhdot(\ecomplex)$ be its perfect
relative tangent-obstruction complex. Then the set $\Cal K$ 
of all relative Kuranishi families
is non-empty and the group $\Cal H$ acts 
transitively on $\Cal K$.
\endproclaim

We let $F_1\mh E_2\dual\to M$ be the map 
$E_2\dual\to E_1\dual\sub M$ induced by 
$\sigma\mh \e_1\to\e_2$ in the complex $\ecomplex$. As before, we denote by $(F_1)$
the ideal in $M$ generated by the components of $F_1$. 
Let $J_1=(F_1)+M_2$. Since
$\coker\{\sigma\dual\}\cong\Omega_{Z/S}$, there is a canonical $S$-morphism
$$\Phi_1: \spec_A M/J_1\lra Z
$$
such that the isomorphism between $\Omega_{(\spec M/J_1)/S}
\otimes_{\o_{\spec M/J_1}}\o_S$
and $\Omega_{Z/S}\otimes_{\o_Z}\o_S$ induced by $\Phi_1$ is the identity
map.

The existence part of Proposition 2.4 follows from the following Lemma.
Recall $O=\coker\{E_1\to E_2\}$.

\pro{Lemma 2.5}
Let the notation be as before. 
Then there are sequences $F_k\in M\oae$ and $\Phi_k\mh\spec
M/((F_{k})+M_k)\to Z$, where
$k=1,\cdots$, of which the following holds.
\roster
\item
$F_1\in M\otimes_A E_2$ and
$\Phi_1\mh\spec_A M/J_1\to Z$ are given as before.
\item
$F_k-F_{k-1}\in M_k\otimes E_2$.
\item
Let $J_k=(F_k)+M_{k+1}$. Then the image of 
$F_k$ in $(\jko/(\mko+\jko\cdot\mo))\otimes_AO$ (using quotient $E_2\to O$)
is the obstruction class
$$o_k=ob(\Phi_{k-1},M/\jko,M/(\mko+\jko\cdot\mo))
$$
to extending $\Phi_{k-1}$ to $\spec M/(M_{k+1}+J_{k-1}\cdot M_1)\to Z$.
\endroster
\endpro

\demo {Proof}
We prove the lemma by induction. Assume that we have constructed a sequence
$F_1,\ldots,F_{k-1}$ satisfying the property of the lemma. We let
$$\iko=\mko+\jko\cdot\mo\sub M.
$$
Note that $M/\jko$ is
a quotient ring of $M/\iko$ and its kernel $\jko/\iko$ is annihilated by the
ideal $M_1$. We let $f_{k-1}$ be the residue class of $F_{k-1}$ in
$(\jko/(M_k+\iko))\oae$. 
We claim that the sequence 
$$\align
(\jko/  \iko)\oae  \mapright{(h_1,h_2)}
&(\jko/(M_k+\iko))\oae\oplus 
(\jko/\iko)\otimes_AO\\
&\lra(\jko/(\M_k+\iko))\otimes_AO\lra 0
\tag 2.2
\endalign
$$
is exact. Indeed, since $h_2$ is surjective and has kernel
$(\jko/\iko)\otimes_A\text{Im}\{E_1\to E_2\}$, the cokernel of
$(h_1,h_2)$ is the cokernel of
$$(\jko/\iko)\otimes_A\text{Im}\{E_1\to E_2\}\lra (\jko/(M_k+\iko))\otimes_A E_2,
$$
which is the last non-zero term in (2.2).
Now we consider $(f_{k-1},o_k)$ in the middle group of 
the above exact sequence.
By the induction hypothesis and the base change property of 
obstruction class, the
images of $f_{k-1}$ and $o_k$ in 
$(\jko/(M_k+\iko))\otimes_AO$ are the obstruction class 
$$ob(\Phi_{k-1},M/\jko,M/(M_k+\iko)),
$$
hence they coincide. It follows that there is an
$\bar f_k\in (\jko/\iko)\oae$
such that its image under $(h_1,h_2)$ is $(f_{k-1},o_k)$.
Now we choose $F_k$. We first select an $F_k\pri\in\jko\oae$ 
so that its residue class is
$\bar f_k$. Since $\bar f_k\equiv f_{k-1}\mod M_k$, which by 
definition is the residue of $F_{k-1}$ in $(M/(M_k+\iko))\oae$, 
it follows that
$$F_k\pri-F_{k-1}\in M_k+\jko\cdot M_1.
$$
Therefore, we can find an $F_k$ so that $F_k-F_{k-1}\in M_k$ and
$F_k-F_k\pri\in \jko\cdot I_{k-1}$. 

Let $J_k=(F_k)+M_{k+1}$.
It remains to show that $\Phi_{k-1}$ extends to an $S$-morphism
$$\Phi_k:\spec_A M/J_k\lra Z.
$$
Because of the obstruction theory, it suffices to show that under
$h_2$ the residue class of $F_k$ in $(\jko/\iko)\oae$
is mapped to the obstruction class $o_k$, because then the
obstruction class
$$ob(\Phi_{k-1},M/\jko,M/J_k)
$$
will be the image of $F_k$ in $(\jko/J_k)\otimes_AO$, which will be zero.
But this is exactly the condition imposed on $F_k$ in our selection.
This proves that $\Phi_{k-1}$ lifts to $\Phi_k$ as desired.
Finally, it follows from $F_{k+1}-F_k\in M_k$ that $\lim F_k=F\in M\oae$
exists and $F-F_k\in M_k$. Also, since $\Phi_{k+1}$ is an 
extension of $\Phi_k$, the limit
$\lim\Phi_k=\Phi: \spec_AM/(F)\to Z$ 
is an $S$-morphism.
This proves the lemma.
\qed

To complete the existence part of Proposition 2.4, it remains 
to show that $\Phi$ is 
an isomorphism and the induced perfect tangent-obstruction complex
from $(F,\Phi)$ is identical to $\tdot_{Z/S}=\hhdot(\edot)$.
We now show that $\Phi$ is an isomorphism. By our technical assumption,
$Z$ embeds in $\spec\hatsym(H)$ for some finitely generated free $A$ module $H$.
Without loss of generality, we can assume that
$\rank H=\rank E_1$. Let $N=\spec\hatsym(H)$, let $N_1\sub N$ be the
ideal generated by $H\sub N$ and let $N_k=N_1^k$.
We let $K\sub N$ be the ideal of $Z\sub \spec N$. Let $\Phi\sta\mh N/K\to
M/(F)$ be the ring homomorphism induced by $\Phi$. Because
$\Phi$ induces an isomorphism between $\Omega_{M/(F)}\otimes_{M/(F)}A$
and $\Omega_{N/K}\otimes_{N/K}A$, $\Phi\sta$ induces an isomorphism
$N/(K+N_2)\cong M/((F)+M_2)$. It follows that we can find an
$A$-isomorphism $\phi\mh N\to M$ such that $\phi(K)\sub (F)$ and
$\phi/K\mh N/K\to M/(F)$ is $\Phi\sta$.
We now show that $\phi(K)=(F)$. 
Let $k$ be the least integer so that
$\phi(K)+M_{k+1}
\ne (F)+M_{k+1}$. Clearly, $k$ must be at least
2. Let $J=\phi(K)+M_{k+1}$.
Since $J+M_k=(F)+M_k$, $((F)+M_k)/J$ is annihilated by $M_1$.
Now let $o$ be the obstruction class to extending
$\Gamma(\o_Z)=N/K\lra M/((F)+M_k)$ to $\Gamma(\o_Z)\to M/J$.
Because such an extension does exist, we have $o=0$. On the
other hand, since $J\sub (F)\cdot M_1+M_{k+1}$, by the definition of $F$
the obstruction $o$ is the residue of $F$ in
$\bl((F)+M_k)/J\br\otimes_A O$. Hence
$$F\in \bl (F)+M_k\br\otimes_A\text{Im}\{E_1\to E_2\}+J
\sub \bl(F)\cdot M_1+M_{k+1}\br+J=J.
$$
This implies that $(F)+M_{k+1}=J$, contradicting to our assumption 
that $J\ne (F)+M_{k+1}$. This proves that $\phi(K)+M_k
=(F)+M_k$ for all $k$, and hence $\Phi$ is an isomorphism.
The proof of that the tangent-obstruction complex of
$(F,\Phi)$ is 
$\tdot_{Z/S}=\hhdot(\edot)$ is straightforward, 
and will be omitted. The existence part of Proposition 2.4 is proved.

Now we study the group action on $\Cal K$. 
We first check that if $T=(\xi,\eta)\in
\Cal H$ and $(F,\Phi)\in\Cal K$, then $T(F,\Phi)\in\Cal K$.
Let $(F\pri,\Phi\pri)=T(F,\Phi)$. Because of the base change property of
the obstruction class, 
it suffices to show that the the sequence $(F\pri_k,\Phi\pri_k):=
(F\pri,\Phi\pri)$ satisfies the four properties listed in Lemma 2.5. 
Indeed, in case $T=(\xi,\ide)$, 
property 1, 2 and 4 in Lemma 2.5 are obviously satisfied. 
Property 3 also holds because $\xi\mh M\to M$ induces an isomorphism
between $(F)+M_k$ and $(F\pri)+M_k$ for all $k$. 
This shows that $(\xi,\ide)(F,\Phi)\in\Cal K$.
Now we consider $(1_M,\eta)\in\h$. Let $F\pri=\eta(F)$. Then since 
$(F)=(F\pri)\sub M$, property 1 to 4 in Lemma 2.5 hold for $(F\pri,\Phi\pri)$
as well. This proves that $(\xi,\eta)=(1,\eta)\cdot(\xi,1)$ acts on $\Cal K$.

\pro{Lemma 2.6}
${\Cal H}$ acts transitively on $\Cal K$.
\endpro

\demo {Proof}
Let $(F,\Phi)$ and $(G,\Psi)$ be any two elements in $\Cal K$.
By definition, we know that $(F,\Phi)\equiv(G,\Psi)\mod M_2$. Now assume 
that there is a $k\geq 2$ so that $(F,\Phi)\equiv(G,\Psi)\mod M_k$.
We will show that there is a transformation $(\xi,\eta)\in\h$ satisfying
$\xi\equiv 1_M\mod M_k$ such that $(\xi,\eta)(F,\Phi)\equiv
(G,\Psi)\mod M_{k+1}$. 

Let $\jko=(F)+M_k$, which is $(G)+M_k$ by the assumption.
Let $\iko=M_{k+1}+\jko\cdot M_1$
and let $f_k$ and $g_k$ be the residue classes of $F$ and $G$ in
$(\jko/\iko)\oae$ respectively. Let 
$$\beta\mh
(\jko/\iko)\oae\to(\jko/\iko)\otimes_AO
$$
be the obvious homomorphism. 
By definition, $\beta(f_k)$ and $\beta(g_k)$ are the obstruction classes to
extending $\Phi_{k-1}=\Psi_{k-1}\mh\spec_AM/\jko\to Z$ to $\spec_A M/\iko$. It
follows that $\beta(f_k)=\beta(g_k)$, and hence 
$$f_k-g_k\in (\jko/\iko)\otimes_A\text{Im}\{E_1\to E_2\}.
$$
Let $t\in M\otimes_A E_2$ be a lift of $f_k-g_k$. Then
$F_k-G_k-t\in \iko\oae$. On the other hand, since 
$M\otimes_A  \text{Im}\{E_1\to E_2\} \sub J_1\oae\sub M_1\otimes_A E_2$,
$$t\in (J_1\cdot \jko)\oae\sub\iko\oae.
$$
Therefore $f_k-g_k\in\iko\oae$. This implies that for some $\eta\in\homem$ 
satisfying property (2.1), we have $G-\eta(F)\equiv0\mod M_{k+1}$.
Hence $J_k=(F)+M_{k+1}=(G)+M_{k+1}$.

Next we analyze $\Phi_k$ and $\Psi_k$. Let 
$\varphi_k,\,\psi_k:\Gamma(\o_Z)\lra M/J_k$
be the homomorphisms of rings induced by
$\Phi_k$ and $\Psi_k$ respectively. Then since $\varphi_k\equiv\psi_k\mod
M_k$, there is a $D\in \text{Der}_A(\Gamma(\o_Z),\jko/J_k)$ 
such that $\psi_k=\varphi_k+D$
(see [Ma, p191]). Since $(\jko/J_k)\cdot M_1=0$ and since $\Phi_k$ induces 
an isomorphism between $\Omega_{Z/S}\otimes_{\o_Z}\o_S$ and
$\Omega_{Y/S}\otimes_{\o_Y}\o_S$,
where $Y=\spec_A M/J_k$, there is a $D_0\in\text{Der}_A(M/J_k,\jko/J_k)$
so that $\psi_k=(\text{id}+D_0)\circ\phi_k$. Since $\text{id}+D_0$ is an isomorphism of
$M/J_k$ that is the identity modulo $M_k$, there is an isomorphism
$\xi\mh M\to M$ so that $\xi\equiv 1_M\mod M_k$, $\xi(J_k)=J_k$ and the 
induced homomorphism $M/J_k\to M/J_k$ is exactly $\text{id}+D_0$.

The transformation 
$(\xi,\eta)$ is not quite what we want, since it satisfies the
relation 
$$G\equiv \eta(F)\mod M_{k+1}
\qquad\and\qquad
\Psi\equiv \Phi\circ\bar \xi\mod M_{k+1}.
$$
To obtain $\eta$ so that $G=\bl(1_M\otimes\eta)\circ(\xi\otimes\ide)\br(F)$, we instead
look at the relative Kuranishi family $(F\pri,\Phi\pri)=(\xi,\ide)(F,\Phi)$.
Then since
$$(F\pri,\Phi\pri)\equiv (F,\Phi)\equiv (G,\Psi)\mod M_k,
$$
by the previous argument we can find an $\eta\in\homem$
satisfying (2.1) so that $G\equiv\eta(F\pri)\mod M_{k+1}$.

Now we apply induction on $k$. The previous argument 
shows that there is a sequence of 
transformations $T_k\in \Cal K$ so that if we let 
$S_k=T_k\circ\cdots\circ T_2$ then
$S_k(F,\Phi)\equiv(G,\Psi)\mod M_{k+1}$. Let $T_k$ be $(\xi_k,\eta_k)$.
Since $\xi_k\equiv 1_M\mod M_k$, $\xi_k\circ\cdots\circ\xi_2$
converges to an automorphism $\xi_{\infty}\mh M\to M$. $\xi_{\infty}$ satisfies
the property (2.1). Now let $(F\pri,\Phi\pri)=(\xi_{\infty},\ide)(F,\Phi)$.
For the same reason, there is a sequence 
$\eta_k\in\Hom_A(E_2,M_1\otimes_AE_2)$ 
so that if we let $H_1=F\pri$ and $H_{k+1}=(\ide+\eta_k)(H_k)$,
then $G\equiv H_{k+1}\mod M_{k+1}$. Applying the Artin-Rees lemma
to the ideal $L=(F\pri)\sub M$, we can find an integer $c$ so that
$$LM_1\cap M_1^n=M_1^{n-c}(LM_1\cap M_1^c)
$$
for any $n>c$. Let $L_k=(H_k)$. Since $L$ is isomorphic to $L_k$ under 
an isomorphism of $M$, the same identity holds with $L$ replaced by $L_k$.
Then since 
$$\eta_k(H_k)=H_{k+1}-H_k\in\Hom_A(E_2,E_2)\otimes_A(L_kM_1\cap M_1^{k}),
$$
for $k\geq c$, we can assume that $\eta_k$ have already been chosen so that
$$\eta_k\in\Hom_A(E_2,E_2)\otimes_AM_1^{k-c}.
$$
With this choice of $\eta_k$, the composite
$(\ide+\eta_k)\cdots(\ide+\eta_1)$ will converge to an $\eta_{\infty}\in
\homae$ such that $\eta_{\infty}\equiv\ide\mod M_1$ 
and $G=\eta_{\infty}(F\pri)$.
Therefore, $(\xi_{\infty},\eta_{\infty})(F,\Phi)=(G,\Psi)$. 
This proves Lemma 2.6 and Proposition 2.4.
\qed

\pro{Corollary 2.7}
Let $Z/S$ and $\tdot_{Z/S}=\hhdot(\edot)$ be as in the situation of
Proposition 2.4. Let $T\sub S$ be a closed subscheme.
Let $W=Z\times_ST$ and $B=\Gamma(\o_T)$.
Then $\tdot_{Z/S}$ canonically induces a relative perfect
tangent-obstruction complex $\tdot_{W/T}=\hhdot(\edot\otimes_{\o_S}\o_T)$.
Further, if $(F,\Phi)$ is a relative Kuranishi family of $\tdot_{Z/S}$,
then $F\pri=F\otimes_AB$ 
and the restriction of $\Phi$ to $\spec N/(F\pri)$, where
$N=M\otimes_AB$,
is a relative Kuranishi family of $\tdot_{W/S}=\hhdot(\edot\otimes_{\o_S}\o_T)$.
\endpro

\demo {Proof}
This is obvious from the proof of Lemma 2.5.
\qed

Before we close this section, we will point out the relation
between the relative Kuranishi families and the Kuranishi families
in the usual sense. Let $S$ and $\tdot_S=\hhdot(\edot)$ be an affine
scheme and a perfect tangent-obstruction complex of $S$. Let
$Z$ be the formal completion of $S\times S$ along its
diagonal and let
$\tdot_{Z/S}=\hhdot(\edot)$ be its 
induced perfect relative tangent-obstruction complex.
Let $(F,\Phi)$ be a relative
Kuranishi family of this complex. 
In the following, we will localize $(F,\Phi)$
and compare it with the usual Kuranishi maps.

Let $q\in S$ be any closed point and let $\mm\sub A$ be the maximal ideal
of $q\in S$. Let $\hat A=\lim A/\mm^n$, let $\hat E_i=E_i\otimes_A\hat A$
and let $\hat M=M\otimes_A \hat A$ be their respective
formal completions. We denote by $\hat Z$ the formal completion of
$Z$ along $Z\times_S\{q\}$ and by $\hat S$ the formal completion of
$S$ along $q$. Then $\tdot_{Z/S}=\hhdot(\edot)$ canonically induces a
perfect relative
tangent-obstruction complex $\tdot_{\hat Z/\hat S}=\hhdot(\hatedot)$.
Obviously, $(F,\Phi)$ induces a Kuranishi family
$$\hat F\in \hat M\otimes_{\hat A}\hat E_2
\qquad\and\quad
\hat\Phi: \spec \hat M/(\hat F)\mapright{\cong} \hat Z
$$
of $\tdot_{\hat Z/\hat S}=\hhdot(\hatedot)$.

Now we turn to the usual Kuranishi families.
For simplicity, we assume $\e_2\otimes_{\o_S}k(q)$ is isomorphic to
$\hh^2(\edot)\otimes_{\o_S}k(q)$.
Let $T_i=\e_i\otimes_{\o_S}k(q)$. Then $T_1$ is the tangent space $T_qS$.
The complex $\tdot_S$ induces an obstruction theory to deformations
of $q$ in $S$ taking values in $T_2$.
Now let $B=\hatsym(T_1\dual)$ and let
$$f\in B\otimes_k T_2\dual\qquad\and\qquad
\varphi:\spec B/(f)\mapright{\cong} \hat S
$$
be a Kuranishi family (cf. [La]). In the following,
we will construct a pair
$(\hat f,\hat\varphi)$ from $(f,\varphi)$ analogous to
$(\hat F,\hat\Phi)$.
Let $I\sub \bob$ be the ideal generated by $a\otimes 1-1\otimes a$ and
let $\hatbob=\lim\bob/I^n$. Let $p_1,p_2\mh
B\to\hatbob$ be the homomorphisms defined by $p_1(a)=a\otimes 1$ and
$p_2(a)=1\otimes a$. For $f\in B\otimes_kT_2$ given before, we denote by
$p_1(f)$ the image of $f$ under $p_1\otimes 1_{T_2}\mh B\otimes T_2\to
\hatbob\otimes_kT_2$. As before, we denote by $(p_1(f))\sub\hatbob$
the ideal generated by the components of $p_1(f)$.
We let $C=\hatbob/(p_1(f))$. $C$ is an $\hat A$-algebra via
$\hat A= B/(f)\to C$ induced by $p_1$. Let $\hat f=p_2(f)$.
It follows that $\spec C/(\hat f)$ is the
formal completion of $\hat S\times
\hat S$ along its diagonal. Because $\hat Z$ is the formal
completion of
$\hat S\times \hat S$ along its diagonal, we obtain a
canonical $\hat S$-isomorphism
$$\hat\varphi: \spec C/(\hat f)\mapright{\cong} \hat Z.
\tag 2.3
$$
Now let $I_C\sub C$ be the ideal generated by the images of $I\sub\bob$.
Note that $C/I_C\cong \hat A$.
We consider the complex of $\hat A$-modules
$$\hatfdot=[ (I_C/I_C^2)\dual\mapright{} T_2\otimes_k\hat A]
$$
index at 1 and 2, where the arrow is the dual of $d\hat f\mh T_2\dual\otimes
\hat A\to I_C/I_C^2$.
Since $H^{-1}(\Hom_{\hat A}(\hatfdot,\ohp))=\Omega_{B/(f)}$,
there is an isomorphism of complexes
$$\Hom_{\ohp}(\hatedot,\ohp)\mapright{\cong}
\Hom_{\ohp}(\hatfdot,\ohp)
\tag 2.4
$$
so that their induced isomorphism on $H^{-1}$ is the
canonical isomorphism between $\Omega_{\hat S}$ and $\Omega_{B/(f)}$.
Namely, we have the commutative diagram
$$\CD
H\upmo(\Hom_{\hat A}(\hatedot,\hat A))@= \Omega_{\hat S}\\
@V{\cong}VV @VV{\cong}V\\
H\upmo(\Hom_{\hat A}(\hatfdot,\hat A))@>>> \Omega_{B/(f)}
\endCD
$$
Let $(r_1,r_2):  \hatedot\to \hatfdot$, where $r_i\mh \hat E_i\to\hat F_i$
are the corresponding
isomorphisms. Since $\hatfdot\otimes_{\hat A}k$ and
$\hatedot\otimes_{\hat A}k$ are
$T_1\mapto{\times0} T_2$, we can choose $r_2$ so that its tensoring with
$k$ is the identity of $T_2$.
We now compare the pairs $(\hat f, \hat\varphi)$ and $(\hat F,\hat\Phi)$.
Let $\xi_1\mh C\to \hat M$ be the $\hat A$-isomorphism induced
by the dual of $r_1\mh \hat E_1\to (I_C/I_C^2)\dual$
and let $\eta_1\mh T_2\otimes_k\hat A\to \hat E_2$
be $r_2\upmo$. Then our choice of
$r_1$ and $r_2$ guarantees that
$$\bl(1_{\hat M}\otimes\eta_1)\circ(\xi_1\otimes 1_{T_2})\br(\hat f)
\equiv \hat F\mod \hat M_2
\qquad\and\qquad
\hat\varphi\circ\bar\xi_1\equiv\hat\Phi\mod\hat M_2,
$$
where $\bar\xi_1\mh \spec \hat M/(\hat F)\to\spec C/(\hat f)$ is the
isomorphism induced by $\xi_1$.

\pro{Lemma 2.8}
There is an $\ohp$-isomorphism $\xi\mh C\to \hat M$ and
$\eta\in \Hom_{\ohp}(T_2\otimes_k\hat A,\hat M\otimes_{\hat A}\hat E_2)$
satisfying
$$\xi\equiv\xi_1\mod \hat M_2
\qquad\and\qquad
\eta\equiv\eta_1\mod\hat M_1
\tag 2.5
$$
such that
$$((1_{\hat M}\otimes\eta)\circ(\xi\otimes 1_{T_2}))(\hat f)
=\hat F\qquad\and\qquad
\hat\varphi\circ\bar\xi=\hat\Phi,
$$
where, as usual, $\bar\xi$ is the isomorphism induced by $\xi$.
\endpro

\demo{Proof}
The proof is parallel to that of Lemma 2.6. The difference is that
in this case we can only compare the obstruction 
classes when they lie in $\hh^2
(\hatedot\otimes_{\hat A}k)$ or in $\hh^2(\hatfdot\otimes_{\hat A}k)$,
because the identification
(2.4) is canonical only after tensoring $k$. We proceed as follows.
Let $\mm_0\sub C$ and $\mm_0\pri\sub \hat M$ be their maximal ideals,
and let $J_k=I_C^2\cdot\mm_0^{k-1}\sub C$ and $J_k\pri
=\hat M_2\cdot{\mm_0\pri}^{k-1}\sub\hat M$.
Assume that there are $\xi_{k-1}\mh C\to \hat M$
and $\eta_{k-1}\in\Hom_{\ohp}(T_2\otimes_k\ohp,
\hat M\otimes_{\hat A}\hat E_2)$ satisfying (2.5) such that
$$((1_{\hat M}\otimes\eta_{k-1})\circ(\xi_{k-1}\otimes
1_{T_2}))(\hat f)\equiv \hat F\mod J_k\pri
\qquad\and\qquad
\hat\varphi\circ\bar\xi_{k-1}\equiv\hat\Phi\mod J_k\pri.
$$
We consider the residue classes
$$\CD
o_k=\text{residue of $ \hat f$ in }\ \bl ((\hat f)+J_k)/((\hat f)\cdot\mm_0+J_{k+1})\br
\otimes_k T_2;\\
o_k\pri=\text{residue of $ \hat F$ in }\ \bl ((\hat F)+J_k\pri)/((\hat F)\cdot\mm_0\pri
+J_{k+1}\pri)\br\otimes_k T_2.\\
\endCD
$$
Because $f$ and $F$ are the (relative) Kuranishi maps,
they are the obstruction class to lifting
$$\ohp\mapright{p_2}\o_{\hat Z}\lra C/((\hat f)+J_k)\qquad\text{to}\qquad
\ohp\to C/((\hat f)\cdot\mm_0+J_{k+1})
$$
and the obstruction class to lifting
$$\ohp\mapright{p_2}\o_{\hat Z}\lra\hat M/((\hat F)+J_k\pri)
\qquad\text{to}\qquad
\ohp\to \hat M/((\hat F)\cdot\mm_0\pri+J_{k+1}\pri),
$$
where $p_2$ is induced by
$Z\to S\times S\mapright{\text{pr}_2} S$.
It follows that they must coincide under the isomorphism
$$\bl ((\hat f)+J_k)/((\hat f)\cdot\mm_0+J_{k+1})\br
\otimes_k T_2\cong \bl ((\hat F)+J_k\pri)/((\hat F)\cdot\mm_0\pri
+J_{k+1}\pri)\br\otimes_k T_2
$$
induced by $\xi_{k-1}$. The remainder argument is a repetition of the
proof of Lemma 2.6, and will be omitted.
This proves the lemma.
\qed

\head
3. Virtual normal cones
\endhead

In the first part of this section, for any quasi-projective scheme $W$ and
a perfect tangent-obstruction complex $\tdot_W=\hhdot(\edot)$, we will
construct a virtual normal cone $C^{\edot}\sub
\vt(\e_2)$. Here $\vt(\e_2)$ is the vector bundle
on $W$ so that its sheaf of sections is $\e_2$.
By abuse of notation, we will not distinguish a vector bundle with the
scheme of its total space.
The cone $C^{\edot}$ will be the restriction to $W$ of the normal cone
to the zero locus of a relative Kuranishi map in its graph.  
Based on the property
of the relative Kuranishi families, we will show that $\ce$ is
unique (as scheme).
Because of this, this construction can be applied to the moduli functors
represented by Deligne-Mumford moduli stacks.

We begin with an affine scheme $S$ and a perfect tangent-obstruction complex
$\tdot_S=\hhdot(\edot)$. Let $Z$ be the formal completion of $S\times S$
along its diagonal. We continue to use the convention adopted in the
previous section. Namely, $A=\Gamma(\o_S)$, $E_i=\Gamma(\e_i)$,
which are assumed to be free $A$-modules, and
$M=\hatsym(E_1\dual)$. We let
$N=\hatsym(E_2\dual)$. Let $(F,\Phi)$, where $F\in M\oae$, be a
relative Kuranishi family of $\tdot_{Z/S}=\hhdot(\edot)$.
It is clear that
$F$ extends to an $A$-homomorphism $N\to M$ of $A$-algebras. We let
$\Gamma_F\sub \spec N\otimes_AM$ be its graph.
We let $j\mh S\to\spec N$ be the obvious section and let
$$\iota=j\times_S 1: \spec M=S\times_S\spec M\to
\spec N\otimes_AM.
$$
We view $\iota$ as the
0-section of $\spec N\otimes_AM\to\spec M$.
In the following, we will view $Z$ as a subscheme of $\spec M$
via the isomorphism $\spec M/(F)\cong Z$.
It follows that
$$\iota(Z)=\Gamma_F\times_{\spec N\otimes_AM}\iota(\spec M).
$$
We let $\n^F$ be the normal cone to $\iota(Z)$ in 
$\Gamma_F$. $\n^F$ is canonically embedded as a closed
subcone in $\vt(\e_2)\times_SZ$, which is
the normal bundle to $\iota(\spec M)$ in $\spec N\otimes_AM$. 
Finally, we let $C^{\edot}$ be the
restriction of $\n^F$ to $S$:
$$C^{\edot}=\n^F\times_ZS
$$

\noindent
{\bf Remark}.
The cone $\n^F$ is the normal cone to $Z$ in $\spec M$, denoted
$C_{Z/\spec M}$, and $\ce= C_{Z/\spec M}\times_ZS$.
However, using the graph description it is clear how these cones are
canonically embedded in the vector bundles $\vt(\e_2)\times_SZ$
and $\vt(\e_2)$, respectively.

\pro{Lemma 3.1}
let $(F,\Phi)$ and $(G,\Psi)$ be two relative Kuranishi families of
$\tdot_{Z/S}=\hhdot(\edot)$. Then as subschemes of $\vt(\e_2)\times_SZ$,
$$\n^F\times_ZS=\n^G\times_ZS.
$$
In particular, the cone
$C^{\edot}\sub\vt(\e_2)$ does not depend on the choice of the relative
Kuranishi families.
\endpro

\demo{Proof}
By Proposition 2.4,
there is a transformation $(\xi,\eta)\in\Cal H$ so that $(\xi,\eta)(F,\Phi)
=(G,\Psi)$. Let $\theta$ be the automorphism of $N\otimes_AM$
defined by
$$a\otimes 1\mapsto \eta^{inv}(a)
\qquad\and\qquad
1\otimes b\mapsto 1\otimes\xi(b).
$$
Here, $\eta^{inv}\mh N\to N\otimes_AM$ is the homomorphism induced by
$\eta\upmo\in\homem$ such that $\eta\upmo\circ\eta=\ide$.
Let $\bar\theta\mh \spec N\otimes_AM\to\spec N\otimes_AM$ be the
induced isomorphism. Clearly, $\bar\theta$ preserves $\iota(\spec M)$
and induces an isomorphism
between $\Gamma_F$ and $\Gamma_G$.
Hence $\bar\theta$ induces an isomorphism,  denoted
$\bar\theta\lsta$, of the normal bundle to
$\iota(\spec M)$ in $\spec N\otimes_AM$ with itself. It follows that
$\bar\theta\lsta$ induces an isomorphism between $\n^F$ and $\n^G$.
Finally, because $\eta\equiv\ide\mod M_1$ and
$\xi\equiv 1_M\mod M_2$, the restriction of $\bar\theta\lsta$
to $\vt(\e_2)\sub
\vt(\e_2)\times_SZ$ is the identity homomorphism. Therefore,
$\n^F\times_ZS=\n^G\times_ZS$.
This proves the lemma.
\qed

We will call $\ce\sub\vt(\e_2)$ the virtual normal cone of the
tangent-obstruction complex $\tdot_S=\hhdot(\edot)$.

\pro{Lemma 3.2}
Let the notation be as before. Assume $\fdot
=[\f_1\to\f_2]$ is another complex  of locally free sheaves so that
$\tdot_S=\hhdot(\fdot)$. Assume further that there is a surjective homomorphism
of complexes $\fdot\to\edot$ such that the induced isomorphism of
their sheaf cohomologies
$\hhdot(\fdot)\cong\hhdot(\edot)$ is the identity, using the isomorphisms
$\hhdot(\fdot)=\tdot_S=\hhdot(\edot)$. Let $\varphi_2\mh \f_2\to\e_2$ be one of the
homomorphism and let $C(\varphi_2)\mh \vt(\f_2)\to\vt(\e_2)$ be the
induced submersive morphism. Then
$$C(\varphi_2)\upmo(\ce)=C^{\fdot}.
$$
Here by $C(\varphi_2)\upmo(\ce)$ we mean
$\ce\times_{\vt(\e_2)}\vt(\f_2)$.
\endpro

\demo{Proof}
This is a local problem. By shrinking $S$ if necessary, we
can assume that there is an
isomorphism $\fdot\cong\edot\oplus [\o_S^{\oplus a}\mapto{\text{id}}
\o_S^{\oplus a}]$ so that the given $\fdot\to\edot$ is the obvious projection.
Let $F_1\mh E_2\dual\to M$
be a relative Kuranishi family of $\tdot_{Z/S}=\hhdot(\edot)$.
Let $M\pri=\hatsym(A^{\oplus a})$ and let $
F_2\mh A^{\oplus a}\to M\pri$ be induced by $\text{id}\mh A^{\oplus a}\to A^{\oplus a}$.
Then
$$F_1\otimes 1+1\otimes F_2: E_2\dual\oplus A^{\oplus a}\lra M\otimes M\pri
$$
is a relative Kuranishi family of $\tdot_{Z/S}=\hhdot(\fdot)$.
A direct computation on normal cones shows that $C^{\fdot}$ is the pull back of
$\ce$ under the obvious projection $\vt(\f_2)\to\vt(\e_2)$.
This proves that $C(\varphi_2)\upmo(\ce)=C^{\fdot}$.
\qed

Let $q\in S$ be any closed point, let $T_1=T_qS$ and
$T_2=\hh^2(\edot)\otimes_{\o_S}k(q)$. Then $\tdot_S$ provides an obstruction
theory to deformations of $q$ in $S$. As in Section 2, 
we let $\hat S$ be the formal
completion of $S$ at $q$, let $f\mh T_2\dual\to B$,
where $B=\hatsym(T_1\dual)$,
be a Kuranishi map. Let $C^f$ be the normal cone to $\spec B/(f)$
in $\spec B$. It follows from the Remark before that $C^f$ is canonically
embedded in $\vt(T_2)\times_k\hat S$.

\pro{Lemma 3.3}
Let the notation be as above. Then there is a quotient vector bundle map
$j\mh \vt(\e_2)\times_S\hat S\to \vt(T_2)\times\hat S$ extending the given quotient map
$$\vt(\e_2)\times_S\{q\}= \vt(\e_2\otimes_{\o_S}k(q))\to\vt( T_2)
$$
so that as subschemes in $\vt(\e_2)\times_S\hat S$,
$$(\vt(\e_2)\times_S\hat S)\times_{\vt(T_2)\times\hat S}C^f=\hat S\times_SC^{\edot}
$$
\endpro

\demo{Proof}
We first consider the case where $\dim T_2=\rank \e_2$.
Then we are in the situation of Lemma 2.8 and its proof.
We continue to use the notations introduced there.
Let $R_1=\spec \hat M$, let $i_1\mh \hat S\to R_1$ be the obvious
section induced by $S\to\spec M$
and let $\hat F\in \hat M\otimes_{\hat A}\hat E_2$ be the image of $F$
under $M\otimes_A E_2\to \hat M\otimes_{\hat A}\hat E_2$.
We let $R_2=\spec C\sub\spec\hatbob$, $i_2\mh \hat S\to R_2$ be the
section induced by $a\otimes b\mapsto ab$ and let
$\hat f\in C\otimes_k T_2$ be the image of the Kuranishi map
$f\mh T_2\dual\to B$ under
$$B\otimes_kT_2\mapright{p_2\times 1_{T_2}} \bob\otimes_kT_2\lra C\otimes_kT_2,
$$
where $p_2(a)=1\otimes a$.
Let $V_1=\spec\bl\hatsym({\hat E_2}\dual)\br$ and
$V_2=\spec\bl\hatsym(T_2\dual\otimes\ohp)\br$.
Let $0_{V_i}$ be the 0-section of $V_i\to \hat S$.
In the proof of Lemma 2.8, we have shown that there is an isomorphism
$$K: V_1\times_{\hat S}R_1\lra V_2\times_{\hat S}R_2
$$
of which the following holds.
First, it induces an isomorphism between $0_{V_1}\times_{\hat S}R_1$
and $0_{V_2}\times_{\hat S}R_2$,
and induces an isomorphism between the graphs $\Gamma_{\hat F}$
and $\Gamma_{\hat f}$.
Secondly, let
$$\varphi : (\vt(\e_2)\times_S\hat S)\times_{\hat S}R_1\lra
\vt(T_2)\times R_2
$$
be the induced isomorphism between the normal bundle to $0_{V_1}\times
_{\hat S}R_1$ in $V_1\times_{\hat S}R_1$ and the normal bundle to
$0_{V_2}\times_{\hat S}R_2$ in $V_2\times_{\hat S}R_2$.
Then the
restriction of $\varphi$ to the fiber over the closed point of $R_1$ is
the identity homomorphism between $\vt(\e_2)
\times_S\{q\}=\vt(\e_2\otimes_{\o_S} k(q))$ and $\vt(T_2)$.

Now let $\n_1$ be the normal cone to
$\spec \hat M/(\hat F)$ in $R_1$ and let $\n_2$ be the normal cone to
$\spec C/(\hat f)$ in $R_2$.
Note that $\n_1$ and $\n_2$ are canonically embedded
in $V_1\times_{\hat S}\spec \hat M/(\hat F)$ and
in $V_2\times_{\hat S}\spec C/(\hat f)$, respectively.
Let $\bar\varphi$ be the restriction of $\varphi$ to
$\vt(\e_2)\times_S \spec \hat M/(\hat F)$.
$\bar\varphi$ is an isomorphism between
$\vt(\e_2)\times_{\hat S}\spec \hat M/(\hat F)$ and
$\vt(T_2)\times\spec C/(\hat f)$.
Since $K$ preserves the 0-sections and the graphs,
$\bar\varphi(\n_1)=\n_2$ and hence
$$\bar\varphi\bl \n_1\times_{\spec\hat M}\hat S\br
=\n_2\times_{\spec C}\hat S.
$$
Since the term inside the parentheses on the
left hand side is $\ce\times_S\hat S$, to prove the lemma, we suffice to
show that the right hand is $C^f$.
Let $\pi\mh \spec C\to\spec B$ be the morphism induced by
$a\mapsto 1\otimes a$. Clearly, $\pi$ is flat and
$$\spec C\times_{\spec B}
\spec B/(f)=\spec C/(\hat f).
$$
It follows from [Vi, p639] that $C^f\times_{\spec B}\spec C=\n_2$.
Therefore
$$\n_2\times_{\spec C}\hat S=
C^f\times_{\spec B}\spec C\times_{\spec C}\hat S=C^f.
$$
This proves the lemma in case $\rank\e_2=\dim T_2$.
In general, by shrinking $S$ if necessary
we can find a complex $\fdot$ and an isomorphism
of complexes
$\edot\cong \fdot\oplus[\ohp^{\oplus a}\mapto{\text{id}}
\ohp^{\oplus a}]$ such that $a=\rank\e_2-\dim T_2$.
Then the general case follows from Lemma 3.2 and the situation just proved. This
proves the lemma.
\qed

Since $\dim\Gamma_f=\rank\e_1$, we have $\dim C^f=\rank \e_1$. This proves

\pro{Corollary 3.4}
The cone $C^{\edot}$ is equidimensional and has dimension $\rank \e_1$.
\endpro

\pro{Corollary 3.5}
Assume that we have two complexes $\edot$ and $\fdot$ so that
$\hhdot(\edot)=\tdot_S$ and $\hhdot(\fdot)=\tdot_S$.
Assume further that there is a surjective homomorphism
$\varphi\mh \f_2\to\e_2$ so that the following diagram is commutative:
$$\CD
\f_2 @>>> \Cal T^2_S(\o_S)\\
@V{\varphi}VV @|\\
\e_2 @>>> \Cal T^2_S(\o_S).
\endCD
$$
Then as cycles, we have
$$C(\varphi)\sta([\ce])=[C^{\fdot}]\in Z\lsta(\vt(\f_2)).
$$
\endpro

Here $C(\varphi)\mh\vt(\f_2)\to\vt(\e_2)$ is the induced morphism on vector bundles.

\demo{Proof}
Let $q\in S$ be any closed point. Lemma 3.3 says that there are
quotient vector bundle homomorphisms
$$j_1: \vt(\e_2)\times_S\hat S\lra \vt(T_2)\times \hat S
\quad \and\quad
j_2: \vt(\f_2)\times_S\hat S\lra \vt(T_2)\times\hat S,
$$
extending the given $\vt(\e_2)\times_S\{q\}\to\vt(T_2)$ and
$\vt(\f_2)\times_S\{q\}\to \vt(T_2)$
respectively, such that
$j_1\upmo(C^f)=\ce\times_S\hat S$ and
$j_2\upmo(C^f)=C^{\fdot}\times_S\hat S$.
It follows that there is a vector bundle quotient homomorphism
$j\mh\vt(\f_2)\times_S\hat S\to\vt(\e_2)\times_S\hat S$ extending
$\vt(\f_2)\times_S\{q\}\to\vt(\e_2)\times_S\{q\}$ such that
$j\upmo(\ce\times_S\hat S)=C^{\fdot}\times_S\hat S$.
This implies that cycles
$[C^{\fdot}]$ and $C(\varphi)\sta([\ce])$ have the same support along
the fiber over $q$ and that the multiplicities of their respective
components near the fiber over $q$ coincide.
Since $q$ is arbitrary, we must have
$C(\varphi)\sta([\ce])=[C^{\fdot}]$.
This proves the corollary.
\qed

\noindent
{\bf Remark.}
The proof shows that the cycle $[\ce]$ can be characterized as follows.
At each $q\in S$, there is a quotient vector bundle homomorphism
$$j: \vt(\e_2)\times_S\hat S\to \vt(T_2)\times\hat S,
$$
extending $\vt(\e_2)\times_S\{q\}\to \vt(T_2)$, such that
$j\sta([C^f])=r\sta[\ce]$,
where $r\mh\vt(\e_2)\times_S\hat S\to\vt(\e_2)$ is
the induced morphism and is flat.
Clearly, this criterion determines $[\ce]$ completely, if
it exists. The reason we need to use the relative Kuranishi
families is to ensure that $[\ce]$ does exist as a cycle.

Let $S_0\sub S$ be a closed subscheme. Then $\tdot_{Z/S}=\hhdot(\edot)$
induces canonically a relative tangent-obstruction complex
$\tdot_{Z\times_SS_0/S_0}
=\hhdot(\fdot)$, where $\fdot=\edot\otimes_{\o_S}\o_{S_0}$.
Let $(F,\Phi)$ be a Kuranishi family of $\tdot_{Z/S}=\hhdot(\edot)$.
Let $A_0=\Gamma(\o_{S_0})$, $M_0=M\otimes_AA_0$
and $F_i=E_i\otimes_AA_0$. The pair $F_0\in M_0\otimes_{A_0}F_2$ and
$\Phi_0\mh\spec M_0/(F_0)\to Z\times_SS_0$ defined by
$F_0=F\otimes_AA_0$ and $\Phi_0=\Phi|_{\spec M_0/(F_0)}$ is a relative
Kuranishi family of $\tdot_{Z\times_SS_0/S_0}=\hhdot(\fdot)$. We let
$\cc=C_{(\spec M_0/(F_0))/\spec M_0}$ and
set $C^{\fdot}=\cc\times_{\spec M_0}S_0$,
which is canonically embedded in $\vt(\f_2)$. Note that
$\vt(\f_2)=\vt(\e_2)\times_SS_0$.

\pro{Corollary 3.6}
$C^{\fdot}=\ce\times_SS_0$ as subschemes of $\vt(\f_2)$.
\endpro

\demo{Proof}
The proof is similar to the proof of Corollary 3.5. 
We continue to use the notations introduced there. Let $q\in S_0$ be any
closed point and let $\hat S_0$ be the formal completion of $S_0$ along $q$.
It follows that $\hat S_0=\hat S\times_SS_0$. Let $j_1$ be the map constructed
in the proof of Corollary 3.5 and let $j_2$ be the restriction of $j_1$ to
$\vt(\e_2)\times_S\hat S_0$.
Hence, $j_2\mh\vt(\f_2)\times_{S_0}\hat S_0\to\vt(T_2)\times \hat S_0$.
Then Lemma 3.3 shows that $j_1\upmo(C^f)=\ce\times_S\hat S$.
We claim that $j_2\upmo(C^f)=C^{\fdot}\times_{S_0}\hat S_0$.
Indeed, by the proof of Lemma 3.3, it suffices to check that if we
let $C_0$ be the formal completion of $M_0$ along the maximal ideal
$\mm_0$ of $q\in S_0$, then $C_0$ is flat over $B$ via the homomorphism
$B\to C_0$ induced by $p_2$. But this is obvious. This proves the claim.
Therefore,
$$\ce\times_S\hat S_0=j_1\upmo(C^f)\times_{\hat S}\hat S_0=j_2\upmo(C^f)
=C^{\fdot}\times_{S_0}\hat S_0.
$$
This proves that $C^{\fdot}=\ce\times_SS_0$.
\qed

In the remainder of this section, we will construct the virtual normal
cone and the virtual cycle of a perfect tangent-obstruction complex.
We will show in the end that this construction commutes
with the refined Gysin maps.

Let $Z$ be a quasi-projective scheme and
let $\tdot_Z$ be a tangent-obstruction complex of $Z$.
We assume that $\edot=[\e_1\to\e_2]$ is
a complex of locally free sheaves of $\o_Z$-modules
so that $\tdot_Z=\hhdot(\edot)$.
We cover $Z$ by affine open $S\lalp$ such that
$\Gamma_{S\lalp}(\e_i)$ are free $\Gamma(\o_{S\lalp})$-modules
for $i=1, 2$. It follows from  Lemma 3.2
that we have canonical cones $C\lalp^{\edot}\sub\vt(\e_2)\times_ZS\lalp$
of the tangent-obstruction complex $\tdot_{S\lalp}=\hhdot
(\edot\otimes_{\o_X}\o_{S\lalp})$.
By Lemma 3.3, $S\lalp\times_Z C^{\edot}_{\beta}=S_{\beta}\times_Z
C^{\edot}\lalp$ as subcones in $\vt(\e_2)\times_Z (S\lalp\cap S_{\beta})$.
Therefore $C\lalp^{\edot}$ patch together to form a global cone
scheme $C^{\edot}\sub\vt(\e_2)$.

We remark that a global resolution $\tdot_Z=\hhdot(\edot)$ allows
us to construct a global cone as a subscheme in $\vt(\e_2)$. However,
if we only have a locally free sheaf $\Cal V$ making $\ttwo_Z(\o_Z)$ its
quotient sheaf, then we can canonically construct a cone cycle as follows.
Since $\tdot_Z$ is perfect, we can find an open covering $S\lalp$ of $Z$
and complexes $\edot\lalp$ of sheaves of $\o_{S\lalp}$-modules
such that $\hhdot(\edot\lalp)=\tdot_{S\lalp}$ and that there are
quotient homomorphisms $\varphi\lalp\mh \Cal V\otimes_{\o_Z}\o_{S\lalp}
\to \e_{2,\alpha}$ such that
$$\CD
\Cal V\otimes_{\o_Z}\o_{S\lalp} @>>> \ttwo_{S\lalp}(\o_{S\lalp})\\
@V{\varphi\lalp}VV @|\\
\e_{2,\alpha} @>>> \ttwo_{S\lalp}(\o_{S\lalp})
\endCD
$$
is commutative. Because of Corollary 3.5, the flat pull backs
$C(\varphi\lalp)\sta([C^{\edot\lalp}])$ and $C(\varphi_{\beta})\sta([
C^{\edot_{\beta}}])$ coincide over $S\lalp\cap S_{\beta}$. Therefore
they patch together to form a  cycle $[C^{\Cal V}]\in Z\lsta\vt(\Cal V)$.
Because of Corollary 3.5 again, $[C^{\Cal V}]$ is unique.

Now we construct the virtual cycle of a perfect tangent-obstruction
complex $\tdot_Z$. We first present $\ttwo_Z(\o_Z)$ as a quotient sheaf
of a locally free sheaf $\Cal V$, 
which is possible since $Z$ is quasi-projective.
Let $i_V\mh Z\to\vt(\Cal V)$ be the zero section and let $i_V\sta\mh
A\lsta(\vt(\Cal V))\to A\lsta Z$ be the Gysin homomorphism.

\pro{Definition 3.7} 
Let the notation be as before. Then we define the
virtual cycle $[Z]\vir$ of $\tdot_Z$ to be
$$[Z]\vir= i_V\sta[C^{\Cal V}]\in A\lsta Z.
$$
\endpro

In order to show that $[Z]\vir$ is well-defined, we need to check that it
is independent of the choice of quotient homomorphism
$\Cal V\to\ttwo_Z(\o_Z)$. Assume that $\w$ is another locally free sheaf
of $\o_Z$-modules and $\w\to\ttwo_Z(\o_Z)$ is a quotient homomorphism.
We let $\m\pri$ be the pull back defined by the square
$$\CD
\m\pri @>>> \Cal V\\
@VVV @VVV\\
\w @>>> \ttwo_Z(\o_Z).\\
\endCD
$$
Then by making $\m\pri$ a quotient sheaf of a locally free sheaf,
say $\m$, we obtain $\phi_1\mh\m\to\Cal V$ and $\phi_2\mh
\m\to\w$.
It follows from Corollary 3.5 that
$$C(\phi_1)\sta([C^{\Cal V}])=[C^{\m}]=C(\phi_2)\sta([C^{\w}]).
$$
This implies that 
$i_V\sta[C^{\Cal V}]=i_M\sta[C^{\m}]=i_W\sta[C^{\w}]\in A\lsta Z$
as required. So $[Z]\vir$ is well-defined.

Refined Gysin maps play an important role in intersection theory.
Given a fiber product square of schemes
$$\CD
W_0 @>>> W\\
@VVi_0V @VViV\\
X_0 @>\xi>> \ X\,,
\endCD
\tag 3.1
$$
where $\xi$ is a regular embedding of codimension $d$, then the
refined Gysin map
$$\xi\sha: A\lsta W\lra A_{\ast-d} W_0
$$
sends $D\in A\lsta W$ to the intersection of $[C_{D\times_XX_0/D}]$
with the zero section of $i_0\sta N_{X_0/X}$.
In this section, we will show that the refined Gysin map
is compatible to our virtual cycle construction.

Let $W$ be a quasi-projective scheme over $X$ and let $X_0\sub X$ be
a regular embedding. We define $W_0$ by the Cartesian square (3.1).
We assume that $W$ (resp. $W_0$) admits a perfect tangent-obstruction
complex $\tdot_W$ (resp. $\tdot_{W_0}$). Let $\l$ be the sheaf of
normal bundle to $X_0$ in $X$. For any affine $S$, $\eta\mh S\to
W_0\sub W$ and $\f\in\mods$, there is a canonical sheaf
homomorphism
$$\tone_W(\eta)(\f)=\hom_{\o_S}(\eta\sta\Omega_W,\f)\lra
(i_0\circ\eta)\sta\l\otimes_{\o_S}\f,
$$
induced by $i\sta\i_{X_0\sub X}\to \i_{W_0\sub W}$,
that fits into the exact sequence
$$0\lra\tone_{W_0}(\eta)(\f)\lra\tone_W(\eta)(\f)\lra
(i_0\circ\eta)\sta\l\otimes_{\o_S}\f.
\tag 3.2
$$

\pro{Definition 3.8}
We say that $\tdot_W$ and $\tdot_{W_0}$ are compatible with the
Cartesian square (3.1) if (3.2) extends to a long exact sequence
$$\tone_W(\eta)(\f)\lra(i_0\circ\eta)\sta\l\otimes_{\o_S}\f\mapright{\delta}
\ttwo_{W_0}(\eta)(\f)\mapright{r}\ttwo_W(\eta)(\f)\lra 0
$$ of which the following holds. Let $S\to Y_0\to Y$ be a tuple of
$S$-schemes described in Definition 1.2 and let $\i=\i_{Y_0\sub Y}$. 
Let $\eta_0\mh Y_0\to W_0$
be any morphism and let
$o\in\ttwo_{W_0}(\eta_0)(\i)$
be the obstruction class to extending $\eta_0$ to $Y\to W_0$. Then
$r(o)\in\ttwo_{W}(\eta_0)(\i)$ is the obstruction class
to extending $\eta_0$ to $Y\to W$. Secondly, assume that $r(o)=0$.
Then we have an extension, say $\eta\mh Y\to W$. Let $\beta
\in(i_0\circ\eta_0)\sta\l\otimes_{\o_S}\i$ be the canonical
homomorphism $(i_0\circ\eta_0)\sta(\i_{X_0\sub X}/\i_{X_0\sub X}^2)
\to\i$ induced by $\eta$. Then $o=\delta(\beta)$.
\endpro

We make one technical assumption, which usually can be
checked explicitly in applications. 
Let $S\lalp$ be an open covering of $W_0$.
We assume that there are
complexes of locally free sheaves $\edot\lalp$ and $\fdot\lalp$
of $\o_{S\lalp}$-modules fitting into
the exact sequence
$$0\lra[0\to i_0\sta\l\otimes_{W_0}\o_{S\lalp}]\lra\edot\lalp
\lra \fdot\lalp\lra 0
$$
such that $\tdot_{S\lalp}=\hhdot(\edot\lalp)$, 
$\tdot_{S\lalp}=\hhdot(\fdot_{\alpha})$ and
the long exact sequence of sheaf cohomologies
$$\hh^j([0\to i_0\sta\l\otimes_{W_0}\o_{S\lalp}])\lra 
\hh^j(\edot\lalp)\lra\hh^j(\fdot\lalp)\lra
\hh^{j+1}([0\to i_0\sta\l\otimes_{W_0}\o_{S\lalp}])
$$
is the exact sequence given in the above definition. 
We assume that there are sheaves $\e_2$ and $\f_2$
of $\o_{W_0}$-modules so that $\e_{\alpha,2}$ and $\f_{\alpha,2}$
are restrictions of $\e_2$ and $\f_2$ to $S\lalp$ respectively.
Also, over $S\lalp\cap S_{\beta}$, there are isomorphisms
of $\edot\lalp$ and $\edot_{\beta}$ and isomorphisms of $\fdot\lalp$
and $\fdot_{\beta}$ so that the induced isomorphisms on their 
sheaf cohomologies are the identity maps. Finally, we
assume that $\f_2$
can be extended to a sheaf of $\o_W$-modules, say $\tilde\f_2$,
so that $\f_2\to\ttwo_{W_0}(\o_{W_0})$ extends to a quotient
homomorphism $\tilde\f_2\to\ttwo_W(\o_W)$.

\pro{Proposition 3.9}
Let $W_0\sub W$ be defined by the square (3.1) such that their
tangent-obstruction complexes $\tdot_{W_0}$ and $\tdot_W$ are
compatible. Assume further that
the technical conditions stated above are satisfied.
Let $[W_0]\vir$ and $[W]\vir$ be the virtual cycles of $\tdot_{W_0}$ and
$\tdot_W$ respectively. Then
$$\xi\sha[W]\vir=[W_0]\vir.
$$
\endpro

This identity is essentially a statement about
associativity of refined Gysin maps. As usual, we will
transform this problem to a problem about the commutativity of
Gysin maps and then apply the basic Lemma in [Vi] to
conclude the proof of the proposition. We now provide the details
of the proof, which will occupy the rest of this section.

We first introduce some notations. Let
$\hat W_0$ (resp. $\hat W$) be the
formal completion of $W_0\times W_0$ (resp. $W\times W$)
along its diagonal, considered as a scheme over $W_0$ (resp.
$W$) via the first projection of the product. We denote by $p_{W_0}
\mh \hat W_0\to W_0$
(resp. $p_W\mh \hat W\to W$)
the morphism induced by the second projection. We begin with a
locally closed affine subscheme $S\sub S\lalp$.
We fix the complexes $\edot\lalp$ and $\fdot\lalp$ given before.
We let $A=\Gamma(\o_S)$, $F_i=\Gamma(\f_{\alpha,i}\otimes_{\o_{W_0}}\o_S)$,
$E_i=\Gamma(\e_{\alpha,i}\otimes_{\o_{W_0}}\o_S)$ and $L=\Gamma(i_0\sta\l
\otimes_{\o_{W_0}}\o_S)$. By shrinking $S\lalp$ if
necessary, we can assume that all modules $F_i$, $E_i$ and
$L$ are free $A$-modules.
As before, we let $M=\hatsym(E_1\dual)$,
which is canonically isomorphic to $\hatsym(F_i\dual)$ using
$\e_{\alpha,1}\cong\f_{\alpha,1}$. 
We pick a relative Kuranishi family $(f,\varphi)$ of
$\tdot_{\hat W\times_WS/S}$, where $f\in M\otimes_A F_2$ and
$\varphi$ is an $S$-isomorphism $\spec M/(f)\to \hat W\times_WS$.
We now pick a relative Kuranishi family of $\tdot_{\hat W_0\times_{W_0}S
/S}$.

We first pick a splitting $\sigma\mh F_2\to E_2$ of the exact sequence
$L\to E_2\to F_2$. Let $g_1=(1_M\otimes\sigma)(f)$.
Note that $(g_1)=(f)$, hence $\spec M/(g_1)$ is isomorphic to
$\hat W\times_WS$. We denote this isomorphism by $\varphi$.
Let $l=\codim(X_0,X)$. Without loss of generality, we can assume that
near $i(S)\sub X$ the sheaf $\i_{X_0\sub X}$ is generated by $l$ sections, say
$s_1,\cdots,s_l$. We let $\bar s_1,\dots,\bar s_l\in M/(g_1)$ be
the pull backs of $s_1,\cdots,s_l$ via
$$\spec M/(g_1)\lra \hat W\times_WS\mapright{\text{pr}_1}
\hat W\mapright{p_W} W\mapright{i} X.
$$
Note that $(s_1,\cdots,s_l)$ form a basis of
$\i_{X_0\sub X}/\i_{X_0\sub X}^2$ near $i(S)$.
Then $(\bar s_1,\cdots,\bar s_l)$
defines a homomorphism
$\bar g_2\mh L\dual\to M/(g_1)$. Let $\phi\mh
\spec M/(g_1,\bar g_2)\to\hat W_0\times_{W_0}S$
be the morphism induced by $\varphi\mh\spec M/(g_1)\to\hat W\times_WS$.

\pro{Lemma 3.10}
Let $\tau\mh M\otimes_AL\to M\otimes_AE_2$ be the homomorphism
induced by $L\to E_2$. Then there is a lift $g_2\in M\otimes_AL$ of
$\bar g_2\in (M/(g_1))\otimes_AL$ such that
$$g:=g_1+\tau(g_2)\in M\otimes_AE_2
\qquad\and\qquad
\phi:\spec M/(g)\lra \hat W_0\times_{W_0}S
$$
form a relative Kuranishi family of
$\tdot_{\hat W_0\times_{W_0}S/S}=\hhdot(\edot\lalp\otimes_{\o_{W_0}}
\o_S)$.
\endpro

\demo{Proof}
Assume that we have found a lift $h\in M\otimes_AL$ of
$\bar g_2$ such that $g_1+\tau(h)\mod M_k$ and the above
$\phi$ 
satisfy the property in Lemma 2.5. Let $J=(g_1+\tau(h))\sub M$,
$J_{k-1}=J+M_k$ and $I_{k-1}=J\cdot M_1+M_{k+1}$.
We consider the epimorphism $M/J_{k-1}\to M/I_{k-1}$
and its kernel $J_{k-1}/I_{k-1}$. Now
let $O=\ttwo_W(\o_S)$, let $O_0=\ttwo_{W_0}(\o_S)$ and let
$r\mh O_0\to O$ be the homomorphism given in Definition 3.8.
Let $o$ (resp. $o_0$) be the obstruction class to extending
$$\varphi_k=\varphi|_{\spec M/J_{k-1}}: \spec M/J_{k-1}\lra
\hat W\times_WS
$$
to $\spec M/I_{k-1}\to \hat W\times_WS$
(resp. to $\spec M/I_{k-1}\to
\hat W_0\times_{W_0}S$).
Since the obstructions are compatible,
we have that $r(o_0)=o$. Now let $\bar g_1$ and $\bar\tau(h)$ be the
residue classes of $g_1$ and $\tau(h)$ in
$(J_{k-1}/I_{k-1})\otimes_AO_0$, respectively.
Since $f$ is a relative Kuranishi map, we have that 
$r(\bar g_1)=o$. On the
other hand, we know that $\varphi_k$ extends to
$\spec M/(I_{k-1}+(f))\to \hat W\times_WS$.
It follows that the residue class of $\tau(h)$ in
$\bl J_{k-1}/(I_{k-1}+(f))\br\otimes_A O_0$,
which also is
the residue class of $\bar g_2$, is the obstruction class to extending
$\varphi_k$ to
$$\spec M/(I_{k-1}+(f))\lra \hat W_0\times_{ W_0}S.
$$
Therefore $\bl\bar g_1+\bar\sigma(h)\br-o_0$ belongs to
$$\text{Ker}\bigl\{ (J_{k-1}/I_{k-1})\otimes_AO_0
\lra (J_{k-1}/I_{k-1})\otimes_AO\oplus
\bl J_{k-1}/(I_{k-1}+(f))\br\otimes_AO_0\bigr\}.
$$
This proves that there is an $\epsilon_k\in (f)\otimes_AL$ so that
$g_1+\tau(h+\epsilon_k)\mod M_{k+1}$
satisfies the property in Lemma 2.5. It follows from the proof of Lemma 2.5
that we can choose $\epsilon_k$ to be in $\bl(f)\cap M_k\br\otimes_AL$.
Hence an induction on $k$ shows that there is a lift
$g_2\in M\otimes_AL$ such that $g:=g_1+\tau(g_2)\in M\otimes_AE_2$
and $\phi=\varphi|_{\spec M/(g)}$ is a relative
Kuranishi family of $\tdot_{\hat W_0\times_{W_0}S/S}=
\hhdot(\edot\otimes_{\o_{W_0}}\o_S)$. 
\qed

Now let $Z=\spec M$ and let $Z(g)=\spec M/(g)\sub Z$. 
Then $Z$ is a scheme over $S$, thus a scheme over
$W_0$. Let $V_1=\vt(\l)\times_{W_0}Z$, $V=\vt(\e_2)\times_{W_0}Z$
and $V_2=\vt(\f_2)\times_{W_0}Z$. Then $V_1$ is a subbundle of $V$
and $V_2$ is the quotient vector bundle $V/V_1$. 
Let $C_{Z(g)/Z}$ be the normal cone to $Z(g)$ in $Z$. The cone
$C_{Z(g)/Z}$ is canonically embedded in $V\times_Z Z(g)$. 
We let 
$$D_1(S)=C_{Z/(g)/Z}\times_ZS\sub V\times_ZS.
$$
It follows from Lemma 3.1 that for the affine covering
$S\lalp$ of $W_0$, the collection
$\{ D_1(S\lalp)\}$ patches together to form a
cone $\bold D_1$ in $V$. By Definition 3.7, if we let
$\eta_1\mh W_0\to V$ be the zero section,
$$ [W_0]\vir=\eta\sta_1[\bold D_1].
$$

Next, we consider the subscheme $Z(f)=\spec M/(f)\sub Z$
and the normal cone $C_2=C_{Z(f)/Z}$,
which is naturally a subcone of $V_2\times_Z Z(f)$. We let
$C_2\to X$ be the morphism induced by 
$$ V_2\times_Z Z(f)\mapright{\text{pr}_2} Z(f)
=\hat W\times_WS\mapright{\text{pr}_1} \hat W\mapright{p_W}W\mapright{i}X.
$$
The normal cone $C_{C_2\times_XX_0/C_2}$ is canonically a subcone in
$$\bl V_2\times_Z(\hat W\times_WS)\br\times_Z\bl V_1
\times_Z(\hat W_0\times_{W_0}S)\br
=(V_1\times_Z V_2)\times_Z(\hat W_0\times_{W_0}S).
$$
We set
$$D_2(S)=C_{C_2\times_XX_0/C_2}\times_{\hat W_0\times_{W_0}S}S.
$$
For the same reason, the collection $\{D_2(S\lalp)\}$ patch together to form a 
cone $\bold D_2\sub \vt(\l)\times_{W_0}\vt(\f_2)$. We claim that
$$\xi\sha[W]\vir=\eta_2\sta[\bold D_2],
$$
where $\eta_2\mh W_0\to V_1\times_{W_0} V_2$ is the zero section.
By our technical assumption, $\f_2$ extends to $\tilde\f_2$ so that
$\f_2\to\ttwo\f_{W_0}(\o_{W_0})$ extends to
$\tilde\f_2\to\ttwo_W(\o_W)$. Let 
$\n\in Z\lsta \vt(\tilde\f_2)$ be the virtual cone cycle of $\tdot_W$
provided by Lemma 3.1 and Corollary 3.5. 
Then the normal cone cycle $[C_{\n\times_XX_0/\n}]$
is canonically a cone cycle in $\vt(\l)\times_{W_0}\vt(\f_2)$. It follows from 
[Vi, p643] that $\xi\sta[W]\vir=\eta_2\sta[C_{\n\times_XX_0/\n}]$.
However, by using Lemma 3.1 and Corollary 3.6,
we have that $[C_{\n\times_XX_0/\n}]=[\bold D_2]$. Therefore, 
$\xi\sha[W]\vir=\eta_2\sta[\bold D_2]$.

It remains to show that $\eta_1\sta[\bold D_1]\sim_{rat}\eta_2\sta[\bold D_2]$.
Our strategy is to transform it into a problem about commutativity of
Gysin maps and then apply work in [Vi].
Let $\rho_1\mh V_1\to V$ and $\rho_2\mh V\to V_2$ be the embedding and 
the quotient vector bundle morphisms,
and let $1_{V_1}\times\rho\mh V_1\to V_1\times_ZV$
be the product morphism. Let $\Gamma$ be the graph of the relative 
Kuranishi map $g\in M\otimes_AE_2$ and let $0_{V_i}$ be the
scheme of the 0-section of $V_i$. We set
$$Y=V_1\times_Z\Gamma,\quad
X_1=(1_{V_1}\times\rho_1)(V_1)\times
_{V_1\times_ZV}Y
\quad\and\quad
X_2=(0_{V_1}\times_ZV)\times_{V_1\times_ZV}Y.
$$
The scheme
$Y$ is a subscheme of $V_1\times_ZV$ and $X_1$ and $X_2$ are
subschemes of $Y$.
Clearly, $X_1\cong\hat W\times_WS$, $X_2\cong\Gamma$
and $X_1\times_{V_1\times_ZV}X_2\cong\Gamma\times_V0_V$.
It follows that the normal cone $C_{X_2/Y}$ is $V_1\times_Z\Gamma$, 
a cone over $X_2$. Now let $\Cal B_1(S)$ be the normal cone to 
$C_{X_2/Y}\times_YX_1$ in $C_{X_2/Y}$.
Since $C_{X_2/Y}\times_YX_1$ is $V_1\times_Z(\Gamma\times_V 0_V)$, 
the scheme $\Cal B_1(S)$ is the pull back of 
$C_{\Gamma\times_V0_V/\Gamma}\sub V$ under 
$$(V_1\times_ZV)\times_Z(\hat W_0\times_{W_0}S)\mapright{\text{pr}_2}
V\times_Z(\hat W_0\times_{W_0}S).
$$
Let $B_1(S)=\Cal B_1(S)\times_ZS$. Then
$B_1(S\lalp)$ patch together to form a cone 
$$\bold B_1\sub 
\vt(\l)\times_{W_0}\vt(\e_2).
$$
The cone $\bold B_1$ is the pull back of $\bold D_1\sub\vt(\e_2)$
via $\vt(\l)\times_{W_0}\vt(\e_2)\mapright{\text{pr}_2}\vt(\e_2)$. 
Hence if we let
$\eta_0$ be the zero section of $\vt(\l)\times_{W_0}\vt(\e_2)$, 
then $\eta_0\sta[\bold B_1]=\eta\sta_1[\bold D_1]$.

Next, we let $\Cal B_2(S)$ be the normal cone to $C_{X_1/Y}\times_YX_2$ in
$C_{X_1/Y}$. The cone $\Cal B_2(S)$ is canonically a subcone of 
$(V_1\times_ZV)\times_Z(\hat W_0\times_{W_0}S)$.
We claim that 
$$\CD
\Cal B_2(S) @>>>C_{C_2\times_XX_0/C_2}\\
@VVV @VVV\\
(V_1\times_ZV)\times_Z(\hat W_0\times_{W_0}S)@>{1_{V_1}\times\rho_2}>>
(V_1\times_ZV_2)\times_Z(\hat W_0\times_{W_0}S)\\
\endCD
$$
is a fiber square.
We first look at $C_{X_1/Y}$. Let $h\mh V_1\times_ZV\to V_1\times_ZV$
be the isomorphism defined by $h(a,b)=(a,b-\rho_1(a))$. Then
$h|_{(1_{V_1}\times\rho_1)(V_1)}$ is an isomorphism between
$(1_{V_1}\times\rho_1)(V_1)$ and $V_1\times_Z0_V$. Further, under 
$1_{V_1}\times\rho_2\mh 
V_1\times_ZV\to V_1\times_Z V_2$, $h(Y)$ is isomorphic to 
$V_1\times_Z(0_{V_2}\times_{V_2}\Gamma_f)$. 
Therefore, the subcone $C_{X_1/Y}$ of
$$C_{(1_{V_1}\times\rho_1)(V_1)/V_1\times_ZV}\times_{V_1\times_Z V_2}
X_1=
V\times_ZX_1
$$
fits into the Cartesian square
$$\CD
C_{X_1/Y}@>>> C_2\\
@VVV @VVV\\
V\times_ZX_1@>>>  V_2\times_Z(\hat W\times_WS)
\endCD
$$
(Note that $X_1$ is canonically isomorphic to $\hat W\times_WS$).
>From this description, we immediately see that $C_{X_1/Y}\times_YX_2$
fits into the Cartesian square
$$\CD
C_{X_1/Y}\times_YX_2@>>>C_{C_2\times_XX_0/C_2}\times_Z(\hat W_0\times_{W_0}S)\\
@VVV @VVV\\
(V_1\times_Z V)\times_Z(\hat W_0\times_{W_0}S)
@>>> (V_1\times_Z V_2)\times_Z(\hat W_0\times_{W_0}S).\\
\endCD
$$
This proves the claim. Finally, we let $B_2(S)=\Cal B_2(S)\times_ZS$.
For the same reason, $B_2(S\lalp)$ patch together to form a cone
$\bold B_2\sub\vt(\l)\times_{W_0}\vt(\e_2)$.
>From the local description, we see that $[\bold B_2]$
is the pull back of $[\bold D_2]\in Z\lsta\bl\vt(\l)\times_{W_0}\vt(\f_2)\br$.
Therefore $\xi\sha[W]\vir=\eta\sta_0[\bold B_2]$.

It remains to show that $[\bold B_1]\sim_{\text{rat}}[\bold B_2]$.
We will apply the basic Lemma in [Vi] to construct a cycle
$[\bold R]\in Z\lsta \bl\ao\times\vt(\l)\times_{W_0}\vt(\e_2)\br$
such that 
$$[\bold R]\cap [\piao\upmo(0)]-[\bold R]\cap[\piao\upmo(1)]=
[\bold B_1]-[\bold B_2].
\tag 3.3
$$
 The main conclusion of the basic 
Lemma in [Vi] is as follows. Let $Y$ be any reduced and 
equidimensional scheme and let $X_1,X_2\sub Y$ be closed subschemes.
Let $D_1$ be the normal cone to $C_{X_1/Y}\times_YX_2$ in $C_{X_1/Y}$
and let $D_2$ be the normal cone to $C_{X_2/Y}\times_YX_1$ in 
$C_{X_2/Y}$, both are canonically embedded in $C_{X_1/Y}\times
_YC_{X_2/Y}$. Then there is a cycle 
$[R]\in Z\lsta(\ao\times C_{X_1/Y}\times_YC_{X_2/Y})$
such that 
$$[R]\cap[\piao\upmo(0)]-[R]\cap[\piao\upmo(1)]=[D_1]-[D_2].
$$
Further, $R$ is canonical under \'etale base change.
The reason that we can not apply this result directly to our
choice of $X_1,X_2\sub Y$ is that the the ambient scheme
$Y$ in our situation, which is $V\times_Z\Gamma$,
may not be equidimensional. However, because of Lemma 3.3,
we will argue that the basic Lemma still apply.

We let $S=S\lalp$ be one of the open sets in the covering of $W_0$ and let
$T_1,\cdots,T_l$ be the irreducible components of $S$ 
(with reduced scheme structure). We fix one of these components,
denoted by $T$. $T\sub W_0$ is a
locally closed affine subscheme.  We then form cycles $\Cal B_1(T)$ and
$\Cal B_2(T)$. Because the corresponding $Y$ 
in constructing $\Cal B_i(T)$, which is
$(V\times_Z\Gamma)\times_Z(Z\times_ST)$, is reduced and 
equidimensional, the basic Lemma provides us a cycle 
$$\r(T)\in
Z\lsta(\ao\times(V_1\times_ZV)\times_ZT)\br
$$
such that
$$[\r(T)]\cap[\piao\upmo(0)]-[\r(T)]\cap[\piao\upmo(1)]=
[\Cal B_1(T)]-[\Cal B_2(T)].
$$

We need to show that the collection $\{\r(T\lalp)\}$ provides us a
global cycle $[\bold R]$ as required. For this, we need
a comparison Lemma similar to Lemma 3.3.
Let $q\in S$ be any closed point, let $Z_q$ be $Z\times_S\{q\}$,
let $\hat Z$ be the formal completion of $Z$ along $Z_q$ and let
$\hat q$ be the formal completion of $S$ along $q$.
Note that $\hat q$ is canonically a subscheme of $Z_q$.
For $i=1,2$ or $\emptyset$, we let $\hat V_i$
be the formal completion of $V_i\times_Z\hat Z$ along its zero section
$0_{V_i}\times_Z\hat Z$.

\pro{Sublemma}
There is a morphism $\varphi\mh \hat Z\to Z_q$ and there are
isomorphisms
$\phi_i\mh \hat V_i \to (\hat V_i\times_{\hat Z}Z_q)\times
_{Z_q}\hat Z$, where $i=1,2$ and $\emptyset$, of which
the followings hold.
\roster
\item
Let $\iota\mh\hat q\to \hat Z$ be the inclusion induced by
$\hat q\to S\to Z$. Then the restriction of $\varphi$ to
$\iota(\hat q)$ factor through the 
subscheme $\hat q\to  Z_q$ and the factored morphism 
$\iota(\hat q)\to\hat q$ is the identity map
between $\hat q=\iota(\hat q)$ and $\hat q\sub Z_q$.
\item
For $i=1,2$ and $\emptyset$, 
$\phi_i(0_{\hat V_i})=
(0_{\hat V_i}\times_{\hat Z} Z_q)\times_{Z_q}\hat Z$ and the restrictions of
$\phi_i$ to $\hat V_i\times_{\hat Z}Z_q$ are the identity morphisms of
$\hat V_i\times_{\hat Z}Z_q$.
\item
We have the commutative diagrams 
$$\CD
\hat V_1@>>> \hat V@>>> \hat V_2\\
@VV{\phi_1}V @VV{\phi}V @VV{\phi_2}V\\
(\hat V_1\times_{\hat Z}Z_q)\times_{Z_q}\hat Z @>>>
(\hat V\times_{\hat Z} Z_q)\times_{Z_q}\hat Z @>>>
(\hat V_2\times_{\hat Z} Z_q)\times_{Z_q}\hat Z,\\
\endCD
$$
where the lower sequence is induced by
$\hat V_1\mapright{\hat\rho_1}\hat V\mapright{\hat\rho_2}\hat V_2$.
\item
$\phi(\Gamma_g\times_Z\hat Z)=(\Gamma_g\times_Z Z_q)\times_{Z_q}\hat Z$
and $\phi_2(\Gamma_f\times_Z\hat Z)=(\Gamma_f\times_ZZ_q)
\times_{Z_q}\hat Z$.
\endroster
\endpro

\demo{Proof}
The proof is similar to the proof of Lemma 3.3. The only modification is to
make sure that the morphisms $\rho_1$ and $\rho_2$, and
the schemes $\Gamma_g$ and $\Gamma_f$ are compatible. 
This can be done easily following the argument used to construct
the relative Kuranishi families of $\tdot_{W_0}$ from $\tdot_W$.
We will omit the details.
\qed

Now if we view $q\in W_0$ as an affine subscheme, we obtain the schemes 
$\Cal B_1(q)$ and $\Cal B_2(q)$, and the cycle $\r(q)$. 
Assume that $q\in T_i$.
Let $\hat T_i$ be the formal completion of $T_i$ along
$q$. Then the flat morphism
$$(\hat V_1\times_{\hat Z}\hat V)\times_S\hat T_i\lra
(\hat V_1\times_{\hat Z}\hat V)\mapright{\phi_1\times\phi}
\bl(\hat V_1\times_{\hat Z}\hat V)\times_{\hat Z} Z_q\br
\times_{Z_q}\hat Z\mapright{\text{pr}_1}
\hat (V_1\times_{\hat Z}\hat V)\times_{\hat Z} Z_q
$$
induces isomorphisms between respective $X_1$, $X_2$ and $Y$ 
in the construction of $\Cal B_j(\hat T_i)$ and $\Cal B_j(q)$. Let
$$h_1: \ao\times(V_1\times_Z V)\times_S\hat T_i\to
\ao\times(V_1\times_Z V)\times_Z Z_q,
$$
where $h_1$ is induced by $\hat T_i\sub\hat Z\mapright{\varphi}
Z_q$, and let
$$h_2: \ao\times(V_1\times_Z V)\times_S\hat T_i\to
\ao\times(V_1\times_Z V)\times_S T_i
$$ 
be the obvious morphism. Note that both are flat.
We claim that
$$h_1\sta(\r(q))=h_2\sta(\r(T_i)).
\tag 3.4
$$
This is not exactly what was proved in [Vi],
since 
$$H: (\hat V_1\times_{\hat Z} \hat V)\times_S\hat T_i\to
(\hat V_1\times_{\hat Z} \hat V)\times_S\hat Z
\mapright{}
(\hat V_1\times_{\hat Z} \hat V)\times_{\hat Z} Z_q
$$
is not \'etale. However, it is clear that $H$ is of the form
$ \hat U\to U_1\times U_2\mapright{\text{pr}_1} U_1$,
where $U_1$ and $U_2$ are two reduced, irreducible formal
complete schemes each supported at a single closed point
and $\hat U$ is the formal completion of $U_1\times U_2$ along
its closed point. A step by step check of the proof of the basic Lemma 
in [Vi] shows that the isomorphism (3.4) does hold.

It is clear now how to construct the cycle $\r(S)\in
Z\lsta(\ao\times(V_1\times_ZV)\times_ZS)$. 
We let the support of $\r(S)$ be the union
of $\text{Supp}(\r(T_j))$. Because of the identity (3.4),
$\r(S)$ with reduced scheme structure is an 
equidimensional closed subscheme.
Now we assign multiplicity to each irreducible 
component of $\r(S)$. Let $p\in\r(S)$ be a general point of one of
its irreducible component, say $C$. Let $q\in S$ be the closed point
under $p$. We assume $q\in T_j$. The component $C$ corresponds to a
unique irreducible component $C\pri$ in $\r(T_j)$
and a unique irreducible component $C_0$ in $\r(q)$, by (3.4).
Let $m_C$ be the multiplicity of $C_0$ in $\r(q)$. Since $T_j$
is reduced and irreducible, $m_C$ is also the multiplicity of
$C\pri$ in $\r(T_j)$. We assign $m_C$ to be the multiplicity of
$C$ in $\r(S)$. Such an assignment is well-defined. 

\noindent
{\bf Remark.}
The cycle $\r(S)$ has the property that for any closed
point $q\in S$, the flat pull back $h_3\sta(\r(S))$
is isomorphic to the flat pull back $h_1\sta(\r(q))$,
where 
$$h_3: \bl\ao\times\vt(\l)\times_{W_0}\vt(\e_2)\br\times_{W_0}\hat q
\lra \ao\times\vt(\l)\times_{W_0}\vt(\e_2).
$$ 
Of course, if
such a $\r(S)$ exists, it is unique. The construction using
the relative Kuranishi families is to ensure that $\r(S)$ exists.

Finally, it follows from the \'etale base change property that 
$\r(S\lalp)$ patch together to form a cycle
$$[\bold R]\in Z\lsta\bl\ao\times_k\vt(\l)\times_{W_0}\vt(\e_2)\br.
$$
Because $\r(q)$ provides a rational equivalence of $[\Cal B_1(q)]$
and $[\Cal B_2(q)]$, $\bold R$ provides a rational equivalence of 
$[\bold B_1]$ and $[\bold B_2]$. Therefore,
$\xi\sta[W]\vir=[W_0]\vir$. This completes the proof of the 
proposition.

\head 
4. Gromov-Witten Invariants of smooth varieties
\endhead

Let $X$ be a smooth projective variety, $n$, $g$ integers, and
$\alpha \in A_1X/\simalg$. The GW-invariants are defined by
taking intersections on the moduli space of stable maps from
$n$-pointed genus $g$ curves to $X$ such that their image 
cycles are in
$\alpha$.  We denote this moduli space by $\mgna$.  When
$\mgna$ has the expected dimension, then the
GW-invariants can be defined as usual.  However,
this rarely happens. Thus we need to use the virtual moduli cycles
to define these invariants.

Let $S$ be an affine scheme and let $\eta\in\fgna(S)$ 
be an element represented by the morphism $f\mh\x\to X$, where
$\x$ is a curve over $S$ with marked sections $D\sub \x$
understood. Then the standard choice of the
tangent-obstruction complex of $\fgna$ is
$$\tdot\fgna(\eta)(\f)=
\bigl[ \ext_{\x/S}^{\bullet}\bl[f\sta\Omega_X\to\Omega_{\x/S}(D)],
\pi_S\sta\f\br\bigr],
$$
where $\f\in\mods$ and $\pi_S\mh\x\to S$ is the projection. 
We now show that there are complexes of locally free
sheaves over $\fgna$ so that their sheaf cohomologies are $\tdot\fgna$.

We fix a sufficiently ample invertible sheaf $\l$ on 
$X$ and then form the exact sequence
$$\CD
0 @>>>\w_2 @>>> \w_1 @>>> f^* \Omega_X  @>>> 0 \; ,
\endCD
$$
where $\w_1 \to f^*_S \Omega_X$ is the natural surjective
homomorphism 
$$
\pi_S\sta\pi_{S *}
\left( \omega_{\Cal{X}/S} (D)^{\otimes 5} \otimes
  f^*(\l\otimes \Omega_X)
\right)
\otimes
\left(\omega_{\Cal{X}/S} (D)^{\otimes 5} \otimes f^* \l
\right)\upmo
\lra f^* \Omega_X\, .
$$
We then form complexes 
$$\a^{\bullet}_{\eta}= [\w_2 \to 0]\quad\text{and}\quad
 \calb^{\bullet} _{\eta}=
[\w_1 \to \Omega_{\Cal{X}/S} (D)]\,  
$$
indexed at $-1$ and $0$ with $\w_1 \to
\Omega_{\Cal{X}/S}(D)$ the composite $\w_1 \to f^* \Omega_X \to
\Omega_{\Cal{X}/S} (D)$.  Let 
$\cc^{\bullet}_{\eta}=[f\sta\om_X\to\om_{\x/S}(D)]$. 
Then we have an exact sequence of complexes
$$
\CD
0 @>>> \a^{\bullet}_{\eta} @>>> \calb^{\bullet}_{\eta} @>>>
\cc^{\bullet}_{\eta} @>>> 0\,,
\endCD
$$ 
and hence a long exact sequence of sheaf cohomologies
$$
\ext^1_{\x/S} 
\bl\cc^{\bullet}_{\eta}, \Cal{O}_{\Cal{X}}  \br
\to \ext^1_{\x/S}
\bl \calb^{\bullet}_{\eta}, \Cal{O}_{\Cal{X}} \br
\to \ext^1_{\x/S}
\left(\a^{\bullet}_{\eta}, \Cal{O}_{\Cal{X}} 
\right)\to
\ext^2_{\x/S}\bl\cc^{\bullet}_{\eta},\o_{\x}\br
\, .
$$
Since $\l$ is sufficiently ample, $\ext^i_{\x/S}(\calb^{\bullet}_{\eta},
\Cal{O}_{\x})$ and $\ext^i _{\x/S}(\a^{\bullet}_{\eta},
\Cal{O}_{\Cal{X}})$ vanish for $i \neq 1$.  Hence
$$\e_{\eta,1} = \ext^1_{\x/S} (\calb^{\bullet}_{\eta}, \Cal{O}_{\Cal{X}})
\qquad\and\qquad 
\e_{\eta,2} =
\ext^1_{\x/S}(\a^{\bullet}_{\eta}, \Cal{O}_{\Cal{X}})
\tag 4.1
$$
are locally free and
the sheaf cohomology of 
$\ebul_{\eta}= [\e_{\eta,1}\to \e_{\eta,2}]$ is $\tdot\fgna(\eta)$.  
It is straightforward to check that the collection $\{\ebul_{\eta}\}$
satisfies the base change property in Definition 1.1, hence it forms a
complex of sheaves over $\fgna$.

To construct the virtual moduli cycle $[\mgna]\vir$,
we need to address one
technical issue, namely, $\mgna$ does not admit
universal families due to the presence of non-trivial
automorphisms.  An automorphism of a morphism
$f$ from $C$ to $X$ is an automorphism $\varphi\mh C \to C$ 
fixing its marked points such that $\varphi
\circ f = \varphi$.  Because $f$ is stable, $\aut (f)$ is
finite.
There are two approaches to get around this difficulty.
One is to realize the moduli space as a quotient 
by a reductive group, say $G$. The other is to use the
intersection theory on stacks developed in [Vi]. 
The former relies on constructing $G$-equivariant data
and then descending them to the quotient space. This can be
done directly
if the quotient is a good quotient. Otherwise, the \'etale slice of the
group action can be used to study the descent problem.
This approach allows one to work with Fulton's
operational cohomology theory of $\mgna$, rather
than the parallel theory on the moduli stack of $\FF_{\alpha,g,n}^X$.

>From [Al], there is a quasi-projective scheme
$\q$ and a reductive group $G$ acting on $\q$ 
such that $\mgna$ is the categorical
quotient of $\qgna$ by $G$.  Over $\qgna$,
there is a universal family
$$
\{F: D\sub \Cal{X}\lra X\}=
\xi\in\FF^X_{\alpha,g,n}(\qgna)
$$
acted on by $G$.  For any closed point $w \in \qgna$,
the stabilizer $G_w \subset G$ of $w$ is naturally the
automorphism group of $F_w : D_w\sub \Cal{X}_w \to X$.
Now by using this family we can construct the complex 
$\edot:=\edot_{\xi}$ in (4.1), after fixing a very ample invertible sheaf
$\l$. By our construction, 
$\tdot\fgna(\xi)=\hhdot(\edot)$.
Further, both $\edot$ and $\tdot\fgna(\xi)$ are
canonically $G$-linearized and the identity is
$G$-equivariant.

We now construct the virtual cycle $[\mgna]\vir$ with
the complex $\edot$ provided.
Let $z \in \mgna$ be any closed point and $w \in
\pi^{-1} (z)$, where $\pi : \qgna \to \mgna$ 
is the quotient projection.  Let $G_w$ be the stabilizer
of $w$.  By combining an argument in [Ko2] and the
construction of [Al], we can find a $G_w$-invariant slice $S
\subset \qgna$ containing $w$
such that $S / G_w$ is an
\'etale neighborhood of $z \in \mgna$. 
Let $\xi_S\in\fgna(S)$ be the object associated to the
restriction of $F\mh\x\to X$ to fibers over $S$.
Clearly, $\tdot\fgna(\xi_S)$ is a tangent-obstruction
complex of $S$ and $\edot_S=\edot\otimes_{\o_{\q}}\o_S$ 
is the complex whose sheaf cohomology if $\tdot\fgna(\xi_S)$. 
Therefore, by applying the construction in the previous sections, 
we obtain a canonical cone
$$
C^{\edot_S}\sub\vt(\e_{S,2})=\vt(\e_2)\times_{\q}S.
$$
Let $C^{\edot_S}/G_w \subset \vt(\e_2)\times_{\q}S/ G_w$ 
be their quotients. 
Assume that $T$ is another $G_w$-invariant slice passing through
$w\pri\in O_w$ such that $T / G_w$ is an \'etale neighborhood of
$z$, then we obtain the cones $C^{\edot_T}\subset \vt(\e_2)\times_{\q}T$
and their quotients by $G_w$.
Since $\edot$ is $G$-linearized, 
by Lemma 3.2, the pull-back of these schemes 
to $S / G_w\times_{\mgna} T / G_{w}$ from $S / G_w$ and $T / G_{w}$
are naturally isomorphic.  
Hence the collection $\vt(\e_2)\times_{\q}S/G_w$ descends to a scheme 
$\vt_{\m}(\e_2)$ over $\mgna$, and
the collection $C^{\edot_S} / G_w$ descends to a scheme
$C^{\edot}_{\m}$ which is a subscheme of $\vt_{\m}(\e_2)$.
Note that $\vt_{\m}(\e_2)$ is not necessarily a vector bundle. 
We will call $\vt_{\m}(\e_2)$ the $\QQ$-descent of the vector bundle 
$\vt(\e_2)$.

Similar to the ordinary case, we define the virtual moduli cycle 
$[\mgna]\vir$ to be
$$s\sta [C^{\edot}_{\m}] \in A_*\bl\mgna\br\otimes_{\ZZ}\QQ\,,
$$
where $s\mh \mgna\to \vt_{\m}(\e_2)$ is the zero section
and $s\sta$ is the Gysin map.
Note that $s\sta$ is well-defined. One way of
seeing it is by using the description
of Gysin maps in terms of Chern classes [Fu, \S 6.1]. This way, to
define $s\sta$, we suffice
to define the Chern classes of $\QQ$-descents of vector bundles (i.e.
 V-vector bundles), which are known to exist
with rational coefficients.

\pro{Lemma 4.1}
The cycle $[\mgna]\vir$ is independent of the choice of the
complex $\edot$ making $\tdot\fgna(\xi)=\hhdot(\edot)$.
\endpro

\demo{Proof}
We will apply Corollary 3.5 to prove the invariance. Clearly,
by the proof of Corollary 3.5 and the above construction, 
it suffices to show that if $\fdot=[\f_1\to\f_2]
$ is another complex of $G$-linearized locally free
sheaves such that
$\tdot\fgna(\xi)=\hhdot(\fdot)$ and that the identity is
$G$-equivariant, then there is a $G$-linearized
locally free sheaf of $\o_{\q}$-modules
$\k$ and surjective $G$-linearized sheaf homomorphisms
$\k\to\e_2$ and $\k\to\f_2$ such that
$$\CD
\k@>>> \e_2\\
@VVV @VVbV\\
\f_2 @>a>> \tdot\fgna(\xi)(\o_{\q})\\
\endCD
$$
is commutative.
We first let $\k_0$ be the pull-back
of $(a,b)$, where $a$ and $b$ are shown in the above square.
Then $\k_0$ is canonically $G$-linearized.
It remains to find a $G$-linearized locally free sheaf $\k$ so
that $\k_0$ is a $G$-quotient sheaf of $\k$.  
Let $L$ be an ample $G$-linearized line bundle on
$\qgna$.  Such an $L$ exists following [Al].  Let
$w \in \qgna$ be any closed point and $O_w$ its
$G$-orbit, which is closed. 
Using locally free sheaves $\e_2$ and $\f_2$, we
can find a $G$-equivariant surjective homomorphism
$$
\eta_w : \e_2 \oplus \f_2 \lra \k_0 \otimes \Cal{O}_{{O}_w}.
$$
Since $G_w$ is finite, for some power $L^{\otimes n_w}$,
the $G_w$-action on 
$L^{\otimes n_w} \otimes k (w)$ is trivial. Because
$\qgna$ is quasi-projective and $L$ is ample, for some
large $m$, the homomorphism
$$
\eta\pri_w : (\e_2 \oplus \f_2) \otimes \o(L^{n_w m}) 
\lra \k_0 \otimes \Cal{O}_{{O}_w} 
$$
induced by $\eta_w$, which is still surjective and
$G$-equivalent, lifts to a global homomorphism
$$
\eta : (\e_2 \oplus \f_2) \otimes \o(L^{n_w m})
\lra \k_0 \; .
$$
Then applying the Reynolds operator, we can assume that 
$\eta$\ is also $G$-equivariant and its restriction to
$O_w$ is $\eta_w\pri$. Since 
$\qgna$ is quasi-projective, a finite sum of sheaves of this
type gives us the desired $G$-equivariant surjective homomorphism
$\k \to \k_0$.  Therefore, by Corollary 3.5, if we let $\vt_{\m}(\k)$
be the $\QQ$-descent of the vector bundle $\vt(\k)$ and let
$\phi_1\mh \vt_{\m}(\k)\to\vt_{\m}(\e_2)$ and
$\phi_2\mh \vt_{\m}(\k)\to\vt_{\m}(\f_2)$ be the induced morphisms
between the $V$-vector bundles, then
$\phi_1\sta[C^{\edot}_{\m}]=\phi_2\sta[C^{\fdot}_{\m}]$.
Therefore, $[\mgna]\vir$ is independent of the choice of the complex
$\edot$. 
\qed

In the remainder of this paper, we will define the GW-invariants of
any smooth projective variety and prove some of its basic properties.
>From now on, unless otherwise is mentioned we will only consider 
homology theory with rational coefficients. We will denote the
(operational) cohomology and homology 
by $A\sta$ and $A\lsta$ respectively. When the varieties are over
complex numbers, sometimes we will use the singular homology theory,
which we will denote by $H\lsta$.
We now give the definition of 
GW-invariants of any smooth projective variety $X$.  
We fix the $\alpha\in
A_1X/\simalg$ and the integers $g$ and $n$ as before so
that $2g+n\geq 3$. Let
$$
[\mgna]\vir \in A_* (\mgna)\,
$$
be the virtual moduli cycle. By using the Riemann-Roch theorem,
it is a purely
$(3 - \dim X) (g - 1) + n + \alpha \cdot c_1 (X)$
dimensional cycle. Let $\pi\ua_n : \mgna \to \mgn$ be the stable contraction
morphism. The $k$-th marked points in curves naturally induces an
evaluation morphism $e_k : \mgna \to X$. 
We let $ev\mh\mgna\to X^n$ be the product $e_1\times\cdots\times e_n$.
Paired with the cycle 
$[\mgna]\vir$, we obtain a homomorphism
$$
\Psi\lgnax :
 A^* (X) ^{\times n} \times A^* (\mgn) 
  \lra A_* (\mgna) 
$$
defined by 
$$\Psi\lgnax(\beta,\gamma)=
\bl ev\sta(\beta)\cup
(\pi\ua_n)\sta(\gamma)\br\bl[\mgna]\vir\br\,.
$$
Composing $\Psi\lgnax$\ with the degree map
$A_*(\mgna)\to A_0(\mgna)\to\QQ$, we obtain the GW-invariants
$$\psi\lgnax:  A^* (X) ^{\times n} \times A^* (\mgn) 
  \lra \QQ\,.
$$
If we fix a polarization $H$ of $X$ and an integer $d$, we can 
define the GW-invariants
$$\psi^X_{d,g,n}\mh A\sta(X)^{\times n}\times A\sta(\mgn)\lra\QQ
$$
as follows. We let $\fdgn\mh\s\to (sets)^0$ be the moduli functor
of stable morphisms defined similar to $\fgna$ except 
that the condition $f_{ \ast}([C])\in\alpha$ is replaced by 
the condition that the degree of
$ c_1(H)(f_{\ast}([C]))$ is $d$. Because of
[Al], $\fdgn$ is coarsely represented by a projective scheme, denoted
$\mdgn$. The previous construction works for this moduli functor
without any change. Consequently we have the virtual moduli 
cycle $[\mdgn]\vir$ which in turn defines the GW-invariants
$\psi_{d,g,n}^X$. When $X$ is a smooth complex projective
variety, then we can use the ordinary homology theory to define
$$
\psi\lgnax: H^*(X) ^{\times n} \times H^*(\mgn) 
  \lra \QQ\,,
$$
where $\alpha\in H_2(X,\ZZ)$, 
by using the image of $[\mgna]\vir$ in $H_*(\mgna)$.
Here since $\alpha\in H_2(X,\ZZ)$, the moduli functor $\fgna$ 
parameterizes all stable morphisms $f\mh C\to X$ with 
$f\lsta([C])=\alpha$ understood as an identity in the 
singular homology $H_2(X,\ZZ)$. Note that if
$\alpha\not\in H_2(X,\ZZ)\cap H^{0,1}(X,\CC)^{\perp}$, then
$\psi_{\alpha,g,n}^X\equiv0$.

The GW-invariants satisfy some basic properties. One of them is
the invariance under deformations of $X$.
Let $\pi\mh\mh X_T\to T$ be a smooth family of relatively
projective varieties
over $T$. For $t\in T$ we let $X_t=\pi\upmo(t)$ and for 
$\beta\in A\sta (X_T)^{\times n}$ we let $\beta_t\in A\sta(X_t)^{\times n}$ 
be the pull back of $\beta$ under $X_t\to X_T$. 
We fix a relatively  ample line bundle of $X_T/T$.

\pro{Theorem 4.2}
Let $X_T/T$ be as before with $T$ an irreducible smooth curve.
Then for any $d\in\ZZ$ and cohomology classes 
$\beta\in A\sta(X_T)^{\times n}$ and $\gamma\in 
A\sta(\mgn)$, the values of the GW-invariants
$\phi_{d,g,n}^X(\beta_t,\gamma)$
are independent of $t\in T$.
\endpro

When $X_T/T$ is defined over $\CC$, we can
use knowledge of $H_2(X_t,\ZZ)$ to prove a finer version of
the invariance theorem.
Consider the analytic curve
$$\Lambda:= 
R_2\pi\lsta \ZZ_{X_T}\otimes_T R^{k_1}\pi\lsta\ZZ_{X_T}\times_T
\cdots\times_T R^{k_n}\pi\lsta\ZZ_{X_T}
$$ 
over $T$. Clearly any point $w\in\Lambda$ corresponds to 
$w=(\alpha_w, \beta_w)\in H_2(X_t,\ZZ)\times H\sta(X_t,\ZZ)^{\times n}$
for some $t\in T$. Hence for any $\gamma\in H\sta(\mgn)$
we can define 
$$\Psi_{\gamma}: \Lambda\lra \QQ
$$
that assigns $w$ to $\psi_{\alpha_w,g,n}^{X_t}(\beta_w;\gamma)$.

\pro{Theorem 4.2${}\pri$}
Let $X_T$ be defined over $\CC$ as before. Assume $T$ is a 
smooth connected curve, then for any $\gamma\in H\sta(\mgn)$
the function $\Phi_{\gamma}\mh\Lambda\to\QQ$ is locally constant.
\endpro

\demo{Proof of Theorem 4.2}
We first form the relative moduli functor.
For simplicity, we assume $T$ is affine.
Let $\frak{Sch}_T$ be the category of $T$-schemes and let
$\FF_{d,g,n}^{X_T/T}: \frak{Sch}_T\to (sets)^0$
be the functor that sends any $S\in
\frak{Sch}_T$ to the subset of $\FF_{d,g,n}^{X_T}(S)$
consisting of the isomorphism classes of
 $f\mh\x\to X_T$ such that
$$\CD
\x @>f>> X_T\\
@VVV @VVV\\
S @>>> T\\
\endCD
$$
is commutative.
$\FF_{d,g,n}^{X_T/T}$ is coarsely represented by a 
$T$-projective scheme $\m_{d,g,n}^{X_T/T}$ [Al].
Let $\m_{d,g,n}^{X_T/T}\to T$ be the obvious morphism.
Then we have the Cartesian square
$$\CD
\Cal M_{d,g,n}^{X_t} @>>> \m_{d,g,n}^{X_T/T}\\
@VVV @VVV\\
\{t\} @>{\eta}>> T.\\
\endCD
$$
Following the principle of conservation of number
[Fu, section 10.2], to prove the theorem it suffices to show that
$$\eta\sha[\m^{X_T/T}_{d,g,n}]\vir=[\m_{d,g,n}^{X_t}]\vir\,.
$$

We first determine the tangent-obstruction complex of
$\FF_{d,g,n}^{X_T/T}$. Let $S\mapto{c} T$ be any affine scheme over
$T$ and $\xi\in \FF_{d,g,n}^{X_T/T}(S)$ be represented by $f\mh\x\to X$
with marked points $D\sub\x$ understood.
>From the discussion in Section 1, the tangent of $\FF_{d,g,n}^{X_T/T}$
at $\xi$ is
$$\tone\FF_{d,g,n}^{X_T/T}(\xi)(\n)=
\ext^1_{\x/S}\bl\bbul_{\xi},\o_{\x}\otimes_\os\n \br\,,
$$
where $\n\in\mods$ and 
$\bbul_{\xi}=\bigl[f\sta\om_{X_T}\to \om_{\x/S}(D)\bigr]$.
We claim that the obstruction to deformations of $f$ lies in
the kernel of
$$\ext^2_{\x/S}\bl\bbul_{\xi},\o_{\x}\otimes_{\o_S}\n\br
\mapright{\beta}
\ext^2_{\x/S}\bl\cbul_{\xi},\o_{\x}\otimes_{\os}\n\br,
$$
where $\cbul_{\xi}=[f\sta\pi\sta\om_T\to 0]$ and $\pi\mh X_T\to T$.
Here the homomorphism $\beta$ is part of the long exact sequence
of cohomologies induced from the short exact sequence
$$ 0\lra \cbul_{\xi}\lra\bbul_{\xi}\lra\abul_{\xi}\lra 0
$$
induced by $\pi\sta\om_T\to\om_{X_T}\to\om_{X_T/T}$.
Here $\abul_{\xi}$ is the complex
$\bigl[ f\sta\om_{X_T/T}\to\om_{\x/S}(D)\bigr]$.
Indeed, let $S\to Y_0\to Y$ be a tuple of $S$-schemes described
in Definition 1.2 with $\i_{Y_0\sub Y}\cong\n$
and let $f_0\mh \x_0\to X_t$ with marked points 
$D_0\sub\x_0$ be a family of stable morphisms over $Y_0$. 
Let $D\sub\x$ be an extension of $D_0\sub \x_0$ and let $o$
be the obstruction to extending $f_0$ to $f$ over $Y$. Then
using the description of $o$ in Section 1 we see immediately that
$\beta(o)$ is the image of 
$$\{ \pi\sta\Omega_T\to\Omega_{\x_0/S}\}=0\in\Hom(\pi\sta\Omega_T,\Omega_{\x/S})
$$
in $\Ext^2_{\x}(\cbul_{\xi},\pi_S\sta\n)$,
which is zero.
We denote the kernel of $\beta$ by $\ttwo\FF_{d,g,n}^{X_T/T}(\xi)(\n)$.
$\tdot\FF_{d,g,n}^{X_T/T}$ is the tangent-obstruction complex of
$\FF_{d,g,n}^{X_T/T}$. Now assume $S\mapto{c} T$ factor
through $\{t\}\sub T$. Then for $\n\in\mods$ we have the exact sequence
$$\align
0\lra\tone\FF^{X_t}_{d,g,n}(\xi)(\n)&\lra
\tone\FF^{X_T/T}_{,g,n}(\xi)(\n)\lra
\n\lra\\
&\lra
\ttwo\FF^{X_t}_{d,g,n}(\xi)(\n)\lra
\ttwo\FF^{X_T/T}_{d,g,n}(\xi)(\n)\lra 0.
\tag 4.2
\endalign
$$ 
Also, from the description of the obstruction classes in 
section 1, we see immediately
that the obstruction classes are compatible in the sense of
Definition 3.8.

To prove the theorem, we need to choose complexes so that 
they satisfy the technical condition of Proposition 3.9. 
We fix a sufficiently ample invertible sheaf $\l$ on $X_t$.
We first let $\hbul$ (resp. $\kbul$) be the complex constructed in
(4.1) with $f\sta\Omega_X$ replaced by $f\sta\Omega_{X_t}$
(resp. $f\sta\Omega_{X_T}$). We let $\gbul$ be the
complex in (4.1) with $f\sta\Omega_X$ replaced by $f\sta\o_{X_t}$
and with $\Omega_{\x}(D)$ replaced by $0$.
Clear, they fit into the exact sequence
$0\to \hbul\to\kbul\to\gbul\to 0$.
Now let $\edot=\hbul\oplus[\o_S\mapto{\text{id}} \o_S]$,
let $\f^1$ be the kernel of $\k^1\to\g^1\to\g^2$ and $\f^2$
be the kernel of $\k^2\to\g^2$.
Note that $\hh^1(\gbul)=\o_S$.
We pick a homomorphism $\o_S\to \f_1$ so that 
$\o_S\to\f_1\to\hh^1(\gbul)$ is the 
identity homomorphism,
we pick a homomorphism $\o_S\to\h^2$ so that
the induced homomorphism $\o_S\to\hh^2(\hdot)$ is the middle
arrow in (4.2) and pick $\o_S\to \f^2$ be so that
$$0\lra [0\to\o_S]\lra \edot\lra\fdot\lra 0
$$
is exact. Hence its long exact sequence
of the sheaf cohomologies 
is exactly the sequence (4.2), and hence the tangent-obstruction complex 
of $\FF_{d,g,n}^{X_T/T}$ satisfies 
the technical condition of Proposition 3.9. 

Lastly, we need to modify the proof of Proposition 3.9 to
accommodate the fact that $\m_{d,g,n}^{X_T/T}$ has
no universal family, as we did in constructing the virtual
moduli cycle. We will omit the details here 
since it is a repetition of the previous argument.
Note that $\M_{d,g,n}^{X_T}$ is a categorical quotient
by a reductive group and the cones constructed in
the proof of Proposition 3.9 are all canonical under \'etale base change.
This completes the proof of Theorem 4.2.
\qed

\demo{Proof of Theorem 4.2 ${}\pri$}
We still assume $T$ is affine. 
Let $R$ be a connected component of $R_2\pi\lsta\ZZ_{X_T}$. $R$
is a smooth analytic curve \'etale over $T$. Let $X_R=X_T\times_T R$ and
let $\alpha\mh R\to R_2\pi\lsta\ZZ_{X_R}$ be the section induced
by the component $R$. Note that 
$\alpha(s)\in H_2(X_s,\ZZ)$ for $s\in R$. 
Note also that under $H_2(X_s,\ZZ)\to 
H_2(X_R,\ZZ)$ all $\alpha(s)$ have identical images.
We let $\alpha_0\in H_2(X_R,\ZZ)$ be their common images. 
Let $H_R$ be the pull back  of
the relatively ample line bundle on 
$X_T/T$. Then $\alpha(s)\cdot c_1(H_R)\in\ZZ$ is independent of $s$. 
We denote it by $d$. Clearly, for any $t\in T$
the disjoint union of $\M^{X_s}_{\alpha(s),g,n}$ for all $s\in R$
over $t$ is an open and closed subscheme of $\M^{X_t}_{d,g,n}$. 
Since $\M^{X_t}_{d,g,n}$ is projective, this is possible either
$R\to S$ is finite or the set of $s\in R$ of which $\M^{X_s}_{\alpha(s),g,n}
\ne\emptyset$ is discrete. We first look into the second situation.
>From the construction of virtual cycle, it is clear that if $V\sub
\M^{X_t}_{d,g,n}$ is a connected component, then using the induced
tangent-obstruction complex of $V$ we can construct the virtual
cycle $[V]\vir$. It follows from the proof of Theorem 4.2 that
for any connected component 
$V$ of $\M^{X_T/T}_{d,g,n}$ and immersion $\eta_t\mh\{t\}\to T$
we have the identity 
$[V_t]\vir=\eta_t\sha[V]\vir$, where $V_t=V\times_T\{t\}$. 
Since $\M^{X_T/T}_{d,g,n}$
consists of fibers over a discrete point set of $T$, $\eta_t\sha[V]\vir
\sim_{\text{rat}}0$ for all $t\in T$. Hence $[\M^{X_s}_{\alpha(s),g,n}]\vir=0$
for all $s\in R$. 
As to the first situation, since $R\to S$ is finite, $R$ is algebraic.
Hence Theorem 4.2 implies that for any $s\in R$,
$[\M_{\alpha(s),g,n}^{X_s}]\vir=\eta_s\sha[\M^{X_R/R}_{\alpha_0,g,n}]\vir$.
Theorem 4.2${}\pri$ then follows from the principle of
conservation of number. This proves Theorem 4.2${}\pri$.
\qed

The Gromov-Witten invariants are expected to satisfy a set of relations,
as explained in [KM, RT1, 2]. We state these relations in terms of the virtual
moduli cycles. We will provide their proofs except the composition law,
which will be proved in the next section.

We first recall the contraction transformation.
For $n\geq 1$ we let $\fgna\to\FF_{\alpha,g,n-1}^X$ be 
the transformation that sends any family
$f\mh\x\to X$ over $S$ in $\fgna(S)$ to the family $f\pri\mh\x\pri\to X$,
where $\x\pri$ is the curve over $S$ obtained by forgetting the $n$-th labeled
section of $\x$ and then stable contracting the resulting $(n-1)$-pointed
curve relative to $f$ and $f\pri$ is the unique morphism so that
$\x\to\x\pri\mapright{f\pri}X$ is $\x\mapright{f}X$.
We let $\pi_n\mh\mgna\to\mgnma$ be the induced morphism. 
Similarly, we let $p_n\mh\mgn\to \M_{g,n-1}$ be the morphism induced
by forgetting the last sections.

\pro{Theorem 4.3}
The virtual moduli cycle $[\mgna]\vir$ satisfies the following properties:
\roster
\item
$[\mgna]\vir\in A_k\mgna$ where $k=(3-\dim X)(g-1)+n+\alpha\cdot c_1(X)$.
\item
Let $\sigma\in S_n$ be any permutation of $n$ elements and let $\phi_{\sigma}
\mh\mgna\to\mgna$ be the morphism induced by permuting the $n$ marked
points of the domains of $f\in\mgna$. Then $\phi_{\sigma}$ is an isomorphism
and $\phi_{\sigma\ast}[\mgna]\vir=[\mgna]\vir$.
\item
The morphism $\pi_n\mh\mgna\to\M^X_{\alpha,g,n-1}$ is a
flat morphism of relative dimension 1. Further
$$\pi_n\sta[\mgnma]\vir=[\mgna]\vir.
$$
\item
Let $\beta\in A^1X$ and $e_n\mh \mgna\to X$ be the $n$-th evaluation morphism.
Then
$$\pi_{n\ast}(e_n\sta\beta\cdot[\mgna]\vir)=\deg(\alpha\cdot\beta)\cdot[\mgnma]\vir.
$$
\item
Composition law (See the statement and the proof in the next section).
\endroster
\endpro

We remark that what is known as the
fundamental class axiom is a direct consequence of (3) of the theorem.
Let $\xi_n\mh\mh\mgna\to X^n\times\mgn$ be the product of 
$ev\mh\mgna\to X^n$
and the projection $\pi_n^{\alpha}\mh\mgna\to\mgn$. Then another way
to describe the GW-invariants is by the homomorphism
$$ I\lagnx: A\sta(X)^{\times n}\lra A\lsta\mgn
$$
defined by 
$I\lagnx(\beta)=\pi_{2\ast}\bl \pi_1\sta(\beta)(\xi\lsta[\mgna]\vir)\br$,
where $\pi_1$ and $\pi_2$ are the first and the second projections of
$X^{n}\times\mgn$. The fundamental class axiom claims that 
for $n\geq 1$ and $2g+n\geq 4$, and for any 
$\beta\in A\sta(X)^{\times n-1}$, we have
$$I\lagnx(\beta\times 1_X)=p_n\sta I\lagnmx(\beta),
\tag 4.3
$$
where $1_X\in A^0(X)$ is the identity element.
We now show that (3) implies (4.3). Consider the commutative diagram
$$\CD
\mgna @>{\tilde\xi }>> X^{n-1}\times\mgn\\
@VV{\pi_n}V @VV{ \bold 1\times p_n}V\\
\mgnma @>{\xi_{n-1}}>> X^{n-1}\times\mgnm\\
\endCD
$$
where $\bold 1\mh X^{n-1}\to X^{n-1}$ is the identity map and
$\tilde\xi$ is the product of the first $n-1$ evaluation 
morphisms and the projection $\M^X_{\alpha,g,n}\to \M_{g,n}$.
Using the projection formula and the property of cohomology
classes [Fu, Definition 17.1], it is direct to check that
(4.3) follows from the identity
$$(\bold 1\times p_n)\sta \xi_{n-1\ast}([\mgnma]\vir)
={\tilde\xi}\lsta([\mgna]\vir).
\tag 4.4
$$
In light of (3) of the theorem, to prove (4.4) it
suffices to show that for any irreducible variety $Y\sub
\mgnma$,
$$\tilde\xi\lsta\pi_n\sta([Y])=(\bold 1\times p_n)\sta 
\xi_{n-1\ast}([Y])\in
Z\lsta(X^{n-1}\times\mgn).
\tag 4.5
$$
Note that the above square is not necessary a fiber square.
We now prove (4.5). 
Clearly, $\tilde\xi(\pi_n\upmo(Y))=(\bold 1\times p_n)\upmo(\xi_{n-1}(Y))$ 
as sets. Hence we suffice  to show that for any irreducible component
$W\sub \pi_n\upmo(Y)$ such that $\dim W=\dim\tilde\xi(W)$, 
the coefficient of $[\tilde\xi(W)]$ in $\tilde\xi\lsta\pi_n\sta([Y])$
is identical to its coefficient in 
$(\bold 1\times p_n)\sta \xi_{n-1\ast}([Y])$.
Let $W\sub\pi_n\upmo(Y)$ be any irreducible component
such that $\dim W=\dim\tilde\xi(W)$. Let
$w\in W$ be a general point associated to the stable map $f_0\mh C_0\to X$
with the marked points $x_1,\cdots,x_n\in C_0$. Let $E_0$ be the
irreducible component of $C_0$ that contains $x_n$. 
Let $(\tilde C_0,\tilde x_1,\cdots,\tilde x_{n-1})$ be the stable contraction 
of $(C_0,x_1,\cdots,x_{n-1})$. Let $\tilde E_0\sub\tilde C_0$ be the
image of $E_0$. Since 
$\dim W=\dim \tilde\xi (W)$, the map $E_0\to\tilde E_0$ is generically
one-to-one. Let $\tilde x_n\in\tilde E_0$ be the image of $x_n$.

Now let $w\pri=\pi_n(w)$, $z=\pi_n^{\alpha}(w)$ and $z\pri=p_n(z)$($=
\pi_{n-1}^{\alpha}(w\pri)$), and let $G$ be the automorphism group
of $(\tilde C_0,\tilde x_1,\cdots,\tilde x_{n-1})$.
We claim that there are $G$-schemes $U$, $U\pri$, $V$ and $V\pri$ and
a $G$-equivariant fiber square
$$\CD
U @>l>> V\\
@V{h_U}VV @V{h_V}VV\\
U\pri @>{l\pri}>> V\pri
\endCD
\tag 4.6
$$
such that their quotients by $G$ are \'etale neighborhoods of
$w\in\mgna$, $w\pri\in\mgnma$, $z\in\mgn$ and $z\pri\in\M_{g,n-1}$, 
respectively, 
and that the induced morphisms shown in the following  makes the square (4.6)
compatible to the middle square in
$$\CD
U @>p>> \mgna @>{\pi_n^{\alpha}}>> \mgn @<q<< V\\
@. @V{\pi_n}VV @V{p_n}VV @. \\
U\pri @>{p\pri}>> \mgnma @>{\pi_{n-1}^{\alpha}}>> \M_{g,n-1} @<{q\pri}<< V\pri.
\endCD
$$
Indeed, we can find a desired $G$-scheme $V\pri$ such that there is
a tautological family $\{D\pri\sub \cc\pri\}\in\FF_{g,n-1}(V\pri)$,
where $\FF_{g,n-1}$ is the moduli functor of 
stable $(n-1)$-pointed curves of genus $g$. Then $V$ can be
chosen as an open subset of the total space of $\cc\pri$ that
contains $\tilde x_n\in\tilde C_0$. By shrinking $V\pri$ and $V$
if necessary, we can assume
that $V\to V\pri$ is smooth with connected fibers. As to $U\pri$,
we can choose it so that in addition to $U\pri/G$ being an 
\'etale neighborhood of $w\pri$ there is a tautological family
$\{f\mh D\sub\cc\to X\}\in\FF_{\alpha,g,n-1}^X(U\pri)$
of which the following holds.
First, there is an isomorphism, denoted by $\varphi$,
between the stable contraction of $\{D\sub\cc\}$ with the pull back of 
$\{ D\pri\sub\cc\pri\}$ under $U\pri\to V\pri$;
Secondly, for any point $s\in V\pri$ we let $\{D_s\pri\sub\cc_s\pri\}$ be
the fiber of $\{ D\pri\sub\cc\pri\}$ over $s$ and let
$A_s$ be the set of isomorphism classes of pairs $(a,b)$,
where $a=\{\psi\mh D_0\sub C_0\to X\}\in\FF_{\alpha,g,n-1}^X(\spec k)$ 
such that the stable contraction of $\{D_0\sub C_0\}$
is isomorphic to $\{D_s\pri\sub\cc\pri_s\}$ and
$b$ is an isomorphism  between the contraction of
$\{ D_0\sub C_0\}$ and the curve $\{D\pri_s\sub\cc_s\pri\}$. 
Then the canonical map $l^{\prime-1}(s)\to A_s$ induced by the 
isomorphism $\varphi$ mentioned in the previous condition is an
isomorphism. Now let $U=U\pri\times_{V\pri}V$
and let $\{\tilde f\mh\tilde D\sub\tilde\cc\to X\}$ 
be the pull back of $\{f\mh D\sub\cc\to X\}$ under
$U\to U\pri$. Because $E_0\to\tilde E_0$ is
generically one-to-one, by shrinking
$V\pri$, $V$ and $U\pri$ if
necessary, there is a unique section $\tilde D_n\mh U \to
\tilde\cc$ such that $\{\tilde f\mh \tilde D\cup\tilde D_n\to X\}\in
\fgna(U)$ and the stable contraction of
$\{\tilde D\cup\tilde D_n\sub\tilde\cc\}$ is isomorphic
to the pull back of
the tautological family over $V\pri$. It is direct to check that
the induced map $U/G\to\mgna$ makes it an \'etale neighborhood of $w$.
Hence the choice of $U$'s and $V$'s satisfy the desired property.
With this choice of $U\pri$ and $V$, we can take
$U=U\pri\times_{V\pri}V$ that satisfies the desired property.
Now let $\tilde W$ be any irreducible component of 
$(\bold 1\times q)\upmo(\tilde\xi(W))$. 
Then with our choice of $U$, etc., it is clear that if we let
$\mu_{[\tilde W]}\bl 
(\bold 1\times q)\sta\tilde\xi\lsta\pi_n\sta([Y])\br$
be the coefficient of $[\tilde W]$ in 
$(\bold 1\times q)\sta\tilde\xi\lsta\pi_n\sta([Y])$, then it is equal to
$$\align
\mu_{[\tilde W]}\bl (\widetilde{ev}\times l)\lsta h_U\sta p^{\prime\ast}([Y])\br
&=\mu_{[\tilde W]}\bl (\bold 1\times h_V)\sta
(ev\pri\times l\pri)\lsta p^{\prime\ast}([Y])\br\\
&=\mu_{[\tilde W]}\bl(\bold 1\times q)\sta 
(\bold 1\times p_n)\sta\tilde\xi\lsta([Y])\br,
\endalign
$$
where $\widetilde{ev}\mh U\to X^{n-1}$
is the composite of $U\to\mgna$ with the product of the
first $n-1$ evaluation morphisms, and
$ev\pri\mh U\pri\to X^{n-1}$ is defined similarly. 
This proves the identity (4.5), and hence (4.6)

\demo{Proof of Theorem 4.3}
By the construction, the cycle
$[\mgna]\vir$ is an equidimensional cycle whose dimension
is the virtual dimension of $\mgna$. Using the Riemann-Roch theorem, one calculates
that it is exactly the $k$ given in the statement. This proves (1).
Also, it is clear that for any $\sigma\in S_n$, 
we have that $\phi_{\sigma}\mh\mgna\to
\mgna$ is an isomorphism of schemes. Since the tangent-obstruction complex
of $\fgna$ does not depend on the ordering of the marked sections,
the virtual moduli cycle will be invariant under $\phi_{\sigma\ast}$. This proves (2). Next, we prove statement (4) assuming property (3). 
Since $X$ is smooth,
$\beta\in A^1$ is the Chern class of a line bundle. By applying 
the projection formula to the flat morphism $\pi_n$, we obtain
$$\align
\pi_{n\ast}(e_n\sta\beta\cdot[\mgna]\vir)&
=\pi_{n\ast}(e_n\sta\beta\cdot\pi_n\sta[\mgnma]\vir)\\
&=\pi_{n\ast}(e_n\sta\beta)\cdot[\mgnma]\vir.
\endalign
$$
Because $\pi_n$ is flat of relative dimension 1, $\pi_{n\ast}(e_n\sta\beta)
=\deg(\alpha\cdot\beta)\cdot{\bold 1}$, where $\bold 1\in A^0\mgnma$ is the
identity element. This proves (4).

Now we prove property (3). Let $G$ be the reductive group and $\q$ be the
$G$-scheme mentioned before so that $\mgnma$ is the
categorical quotient of $\q$. Let $\Cal C$ be the universal curve over $\q$. Then
$\Cal C\to\q$ is flat of relative dimension 1. It follows from the universal property
of $\q$ that $G$ acts canonically on $\Cal C$ and $\mgna$ is the 
categorical quotient of the total scheme of $\Cal C$, denoted $\p$, by $G$. Let
$\pi_n\pri\mh \mgna\to\mgnma$ be the induced morphism. $\pi_n\pri$ is
the morphism described in the statement (3) of the theorem. 
Now we argue that $\pi_n$($=\pi_n\pri$) is flat of relative dimension 1.
It is obvious that $\pi_n$ has relative dimension 1.
Now let $w\in\q$ be any closed point
and let $G_w\sub G$ be the stabilizer of $w$. Then there is a
$G_w$-invariant slice $U\to\q$ such that $U/G_w\to\mgnma$ is
an \'etale neighborhood. It follows that $\cc\times_{\q}U/ G_w$ is
an \'etale neighborhood $\mgna$ and the projection 
$\cc\times_{\q}U/ G_w\to U/G_w$ is compatible to the projection
$\mgna\to\mgnma$. Therefore $\pi_n$ will be flat if we can show
that $\cc\times_{\q}U/ G_w\to U/G_w$ is flat, which follows
from the flatness of $\cc\to\q$ and that $G_w$ is a finite
group. This shows that $\pi_n$ is flat.

It remains to prove the identity in statement (3). 
Let $\tilde f\mh\tilde D\sub \tilde \cc\to X$ 
be the tautological family over $\p$,
where $\tilde D\sub\tilde \cc$ is a family of $n$-pointed curves over $\p$,
characterized by the following property.
There is a canonical morphism
$$\pi: \tilde \cc\lra \Cal C\times_{\q}\p
$$ 
such that the base change  $f\pri\mh \Cal C\times_{\q}\p\to X$ of $f$ with the
marked divisors $D\times_{\q}\p$ is the stable contraction of
$\tilde D_{<n}\sub\tilde \cc$, where $\tilde D_{<n}$ is the first
$(n-1)$-marked sections in $\tilde D$, relative to $\tilde f$;
The restriction of $\pi$ to $\tilde D_n$ is an isomorphism between $\tilde D_n$
and $\Cal C\times_{\q}\p$. 
Let $\xi\in\FF_{\alpha,g,n-1}^X(\q)$ (resp. $\tilde\xi\in\fgna(\p)$) be the
object corresponding to the family $f$ (resp. $\tilde f$).
Following the Remark after Lemma 4.1, to construct $[\mgnma]\vir$,
we suffice to find a $G$-linearized 
locally free sheaf $\v$ of $\o_{\q}$-modules such that
$$\ttwo\FF_{\alpha,g,n-1}^X(\xi)(\o_{\q})=
\ext^2_{\Cal C/\q}\bl[f\sta\Omega_X\to\Omega_{\Cal C/\q}(D)],\o_{\Cal C}\br
\tag 4.7
$$
is a $G$-linearized quotient of $\v$. As before, we pick a $G$-linearized
locally free sheaf $\w_1$ of $\o_{\Cal C}$-modules such that $\w_1\upmo$
is sufficiently ample along fibers of $\Cal C\to\q$. 
We then pick a $G$-equivariant quotient homomorphism 
$\w_1\upmo\to f\sta\Omega_X$ and let $\w_2=\ker\{\w_1\to f\sta\Omega_X\}$. 
$\w_2$ is also locally free. It follows that (4.7)
is a $G$-equivariant quotient sheaf of
$$\v:=\ext_{\Cal C/\q}^1\bl[\w_2\to0],\o_{\Cal C}\br=
\ext^0_{\Cal C/\q}(\w_2,\o_{\Cal C}).
$$ 
By the Remark after Lemma 4.1, there is a canonical cycle $[C^{\v}]\in
Z\lsta\vt(\v)$ so that the image of its $\QQ$-descent over $\mgnma$ under
the obvious Gysin map is $[\mgnma]\vir$. 
Now we pick a similar vector bundle over $\p$. We set $\tilde\w_1=p\sta\w_1$
and $\tilde\w_2=p\sta\w_2$, where $p$ is the composite
$\text{pr}_2\circ\pi\mh\tilde \Cal C\to\Cal C\times_{\q}\p\to\Cal C$.
Then $\w_2\to\w_1\to f\sta\Omega_X$ pulls back to $\tilde\w_2\to\tilde\w_1
\to {\tilde f}\sta\Omega_X$. Let
$$\tilde\v:=
\ext_{\tilde\Cal C/\p}^1\bl[\tilde\w_2\to0],\o_{\tilde\Cal C}\br
\lra\ttwo\fgna(\tilde\xi)(\o_{\p})=
\ext^2_{\tilde\Cal C/\p}\bl[{\tilde f}\sta\Omega_X\to
\Omega_{\tilde\Cal C/\p}(\tilde D)],
\o_{\tilde\Cal C}\br
$$
be the similar quotient homomorphism of sheaves. Because $\pi$ 
contracts at most one rational curve in each fiber of $\tilde\Cal C$ over $\p$, $\tilde\v$ is
canonically isomorphic to $p_n\sta\v$,
where $p_n\mh\p\to\q$ is the projection. We claim that there is a
canonical homomorphism $\phi$ making the following diagram commutative
$$\CD
\tilde v @>>> \ttwo\fgna(\tilde\xi)(\o_{\p})\\
@VV{\cong}V @VV{\phi}V\\
p_n\sta\v @>>> p_n\sta\ttwo\FF_{\alpha,g,n-1}^X(\xi)(\o_{\q}).\\
\endCD
\tag 4.8
$$
We consider the canonical exact sequences
$$\matrix
\ext_{\tilde \Cal C/\p}^1(\Omega_{\tilde\Cal C/\p}(\tilde D),\o_{\tilde \Cal C}) 
&\lra&
\ext_{\tilde \Cal C/\p}^1({\tilde f}\sta\Omega_X,\o_{\tilde \Cal C}) &\lra&
\ext_{\tilde \Cal C/\p}^2(\abul,\o_{\tilde \Cal C})
\\
\downarrow{\phi_2}&&\downarrow{\phi_1}&&\downarrow{\phi}\\
\pi_n\sta\ext^1_{\Cal C/\q}(\Omega_{\Cal C/\q}(D),\o_{\Cal C})&\lra&
\pi_n\sta\ext^1_{\Cal C/\q}(f\sta\Omega_X,\o_{\Cal C}) &\lra&
\pi_n\sta\ext^2_{\Cal C/\q}(\bbul,\o_{\Cal C})
\\
\endmatrix
$$
where $\abul=
[{\tilde f}\sta\Omega_X\to\Omega_{\tilde\Cal C/\p}(\tilde D)]$
and $\bbul=
[f\sta\Omega_X\to\Omega_{\Cal C/\q}(D)]$.
Clearly, there are canonical homomorphisms $\phi_1$ and $\phi_2$ as indicated
in the above diagram making the left square commutative. 
Because the two horizontal arrows on the right are
surjective,
there is a canonical homomorphism $\phi$ making the right square commutative. 
The commutativity of the diagram (4.8) follows immediately.

Now let $[C^{\v}]\in Z\lsta\vt(\v)$ (resp. $[C^{\tilde\v}]\in
Z\lsta\vt(\tilde \v)$) be the virtual normal cone cycle constructed in the
beginning of Section 4 associated to $\tdot\FF_{\alpha,g,n-1}^X(\xi)$
(resp. $\tdot\fgna(\tilde\xi)$), and the quotient homomorphism 
$\tilde\v\to \ttwo\fgna(\tilde\xi)(\o_{\p})$.
Let
$$\Phi: \vt(\tilde\v)=\vt(\v)\times_{\q}\p\mapright{\text{pr}_1}\vt(\v)
$$
be the projection. $\Phi$ is flat of relative dimension 1. 
Now assume that
$$\Phi\sta[C^{\v}]=[C^{\tilde\v}].
\tag 4.9
$$
Let $[C^{\v}_{\m}]$ and $[C^{\tilde\v}_{\m}]$ be the $\QQ$-descents of $
[C^{\v}]$ and $[C^{\tilde\v}]$ to $\mgnma$ and $\mgna$ respectively, 
and let $\vt_{\m}(\v)$ and
$\vt_{\m}(\tilde \v)$ be the $\QQ$-descents of $\vt(\v)$ and $\vt(\tilde\v)$ 
respectively. Then 
$$\CD
[C^{\tilde\v}_{\m}] @>>> \vt_{\m}(\tilde\v)\\
@VVV @VVV\\
[C^{\v}_{\m}] @>>> \vt_{\m}(\v)
\endCD
$$
is a pull back diagram. Therefore,
$$[\mgna]\vir=\eta_{\tilde V}\sta[C_{\m}^{\tilde\v}]=
\pi_n\sta\eta_V\sta[C^{\v}_{\m}]=[\mgnma]\vir,
$$
where $\eta_{\tilde V}$ and $\eta_V$ are the 0-sections of 
$\vt(\tilde\v)$ and $\vt(\v)$.

We now prove the identity (4.9). We first note that there is a largest 
subscheme $\Sigma\sub\p$ characterized by the property that the 
restriction of the contraction morphism 
$$\pi|_{\tilde\Cal C\times_{\p}\Sigma}: \tilde\Cal C\times_{\p}\Sigma\lra\Cal C\times_{\q}\Sigma
$$
contracts a $\bold P^1$-bundle over $\Sigma$. Since the image
of any of these $\bold P^1$ in $\Cal C$ must be either one of the marked points 
or one of the singular points of the curves in this family, 
the restriction of $p_n\mh\p\to\q$ to $\Sigma$ is
a finite morphism.
Let $w\in\p$ be any
closed point over $z\in\q$. Let $\tilde C$ and $C$ be the restriction of $\tilde \Cal C$
and $\Cal C$ to $w$ and $z$ respectively. We let $\tilde\varphi\mh\tilde C\to X$ and 
$\varphi\mh C\to X$ be the corresponding morphisms, and 
$\tilde D\sub\tilde C$ and $D\sub C$ be the corresponding marked points, respectively.
Then $(f,C,D)$ is the stable contraction of $(\tilde f,\tilde C,\tilde D_{<n})$ relative to 
$\tilde f$. We first consider the situation where $\tilde C$ is isomorphic to
$C$. Namely, $\tilde f\in \p-\Sigma$. Clearly
$$\phi_w:=\phi\otimes k(w):
\Ext^2_{\tilde C}\bl [{\tilde\varphi}\sta\Omega_X\to
\Omega_{\tilde C}(\tilde D)],\o_{\tilde C}\br\cong
\Ext_C^2\bl[\varphi\sta\Omega_X\to\Omega_C(D)],\o_C\br
$$
is an isomorphism. We claim that the obstruction theory to deformations of
$\tilde\varphi$ is identical to the obstruction theory to deformations of $\varphi$.
Indeed, let $B$ be any Artin ring with residue field $k$ and let $I\sub B$ be 
an ideal annihilated by the maximal ideal of $B$. Let $B_0=B/I$ and let
$\tilde\varphi_0\mh \tilde D_0\sub\tilde C_0\to X$ be a flat family over $\spec B_0$
whose restriction to the fiber over $\spec k$ is $\tilde \varphi$. Let
$\tilde o$ be the obstruction to extending $\tilde\varphi_0$ to families over
$\spec B$. Similarly, we let $\varphi_0\mh D_0\sub C_0\to X$ be the stable
contraction of $\tilde\varphi_0\mh\tilde D_{0,<n}\sub\tilde C_0\to X$ and let
$o$ be the obstruction 
to extending $\varphi_0$ to families over $\spec B$. By our
description of the obstruction theory of $\fgna$ in Section 1,
$\phi_w(\tilde o)=o$. 
Now let 
$$T_{i,w}=\Ext^i_{\tilde C}([{\tilde\varphi}\sta\Omega_X\to
\Omega_{\tilde C}(\tilde D)],\o_{\tilde C}\br
\quad\and\quad
T_{i,z}=\Ext^i_{C}([{\varphi}\sta\Omega_X\to
\Omega_{C}(D)],\o_{C}\br.
$$
Note that in our situation, $T_{2,z}\cong T_{2,w}$ and $T_{1,z}$ is 
canonically a quotient vector space of $T_{1,w}$ with
$\dim T_{1,w}=\dim T_{1,z}+1$. Let $h\mh T_{1,w}\to T_{1,z}$
be the projection.
Let $f_z\in\hatsym(T_{1,z}\dual)\otimes_k T_{2,z}$ be a Kuranishi map of
of $w$. Then $f_w: =h\sta(f_z)\in\hatsym(T_{1,w}\dual)\otimes_k T_{2,w}$ 
is a Kuranishi map of $w$. Now let
$$\hat z:=\spec \hatsym(T_{1,z}\dual)/(f_z)\sub
\hat T:=\spec \hatsym(T_{1,z}\dual)
$$ 
and let
$$\hat w:=\spec \hatsym(T_{1,w}\dual)/(f_w)\sub
\hat S:=\spec \hatsym(T_{1,w}\dual).
$$ 
Because the normal cones are canonical under flat base change [Vi], the normal 
cone
$$[C^{f_w}]:=[C_{\hat w/\hat S}]\in Z\lsta(\vt(T_{2,w})\times_k\hat w)
$$
is the pull back of
$$[C^{f_z}]:=[C_{\hat z/\hat T}]\in Z\lsta(\vt(T_{2,z})\times_k\hat z)
$$ 
under the canonical flat morphism $\hat w\to\hat z$. Hence, if we let
$U=\p-\Sigma$, then
$$R_U\sta(\Phi\sta([C^{\v}])=R_U\sta([C^{\tilde \v}])\in Z\lsta(\vt(\tilde\v)\times
_{\p}U),
$$
where $R_U\mh\vt(\tilde\v)\times_{\p}U\to\vt(\tilde\v)$ is the open immersion.
Since both $\Phi\sta([C^{\v}])$ and $[C^{\tilde\v}]$ are cycles of identical 
dimensions, to show $\Phi\sta[C^{\v}]=[C^{\tilde \v}]$, we suffice to show that
no irreducible components of $[C^{\tilde \v}]$ are contained in 
$\vt(\tilde\v)\times_{\p}\Sigma$.

Now we consider $w\in\Sigma$. Let $T_{1,w}$ and $T_{2,w}$ be as before,
and let $f_w\in\hatsym(T_{1,w}\dual)\otimes_kT_{2,w}$ be a Kuranishi map of
$w$. By the definition of Kuranishi families, there is an associated family
$$F:\x\lra X\qquad\and\qquad \d \sub\x
$$
of stable morphisms over $\hat w$. We let $\x\pri$ over $\hat w$ be the resulting
curve obtained by stable contracting $\d_{<n}\sub\x$ relative to $F$.
Let $\x\to\x\pri$ be the contraction morphism and let $\hat\Sigma\sub\hat w$
be $g\upmo(\Sigma)$, where $g\mh\hat w\to\p$ is a morphism such that $F$ is the
pull back of the tautological family $\tilde f$ over $\p$. 
It follows that $[C^{\tilde \v}]$ has no components 
supported on $\Sigma$ if and only
if for any $w\in\Sigma$ the cycle $[C_{\hat w/\hat S}]$
has no components supported over $\hat \Sigma$.

Now we prove this statement. Let $z=\pi_n(w)\in\q$, 
let $T_{1,z}$ and $T_{2,z}$ be
as before and let $f_z\in\hatsym(T_{1,z}\dual)\otimes_kT_{2,z}$ be a Kuranishi
map of $z$.
Similarly, we let
$$G:\y\lra X\qquad\and\qquad \e\sub\y
$$
be the associated family of stable morphisms over $\hat z$.
It follows that there is a canonical isomorphism
$\x\pri\times_{\hat w}\{w\}\cong\y\times_{\hat z}\{z\}$.
Now let $r\in\y\times_{\hat z}\{z\}$ be the image of the $\bold P^1\sub
\x\times_{\hat w}\{w\}$ that was contracted under $\x\to\x\pri$. 
Let $\hat r$ be the formal 
completion of the total scheme of $\y$ along $r$. It follows that $\hat r$ is
isomorphic to $\hat w$. Without loss of generality, we can assume that
$\y/\hat z$ can be extended to a family of nodal curves, say $\tilde\y$,
over $\hat T$. In case the total space
of $\tilde\y$ is smooth at $r$, we let $\hat {T\pri}=\hat T$ and let
 $\hat R$ be the formal completion of 
$\tilde\y$ along $r$. 
Otherwise, because $\tilde\y\to\hat T$
is a flat family of nodal curves, by embedding $\hat T$ in 
$\hat {T\pri}=\spec\hatsym(T_{1,z}\oplus k)$, 
we can assume that $\tilde\y/\hat T$ extends to an $\tilde {\y\pri}$ over
$\hat {T\pri}$ so that the total space of
$\tilde{\y\pri}$ is smooth at $r$. Then we let $\hat R$ be
the formal completion of $\tilde{\y\pri}$ along $r$. 
It follows that $\hat w\cong\hat r$ and $\dim\hat R
=\dim T_{1,w}$. Let 
$p\mh\hat R\to\hat {T\pri}$ be the induced projection.
Because $p$ is flat and $\hat r=\hat z\times_{\hat {T\pri}}\hat R$,
$$C_{\hat z/\hat {T\pri}}\times_{\hat z}\hat r\cong C_{\hat r/\hat R}.
$$
However, because $\hat r\cong\hat w$ and $\dim\hat R=\dim\hat S$, 
$\hat r\cong\hat w$ extends to an isomorphism $\hat R\cong\hat S$. 
Therefore, there is an isomorphism
$C_{\hat r/\hat R}\cong C_{\hat w/\hat S}$. 
Finally, because the restriction of the composite
$\hat w\cong\hat r\mapto{p}\hat z$ to $\hat\Sigma\sub\hat w$ is finite,
where $\hat \Sigma$ is the formal completion of
$\Sigma$ along $w$, $C_{\hat r/\hat R}$, and hence
$C_{\hat w/\hat S}$, has no components supported over $\hat\Sigma$.
This proves that $[C^{\tilde\v}]$ has no components supported over 
$\Sigma$ and hence $[C^{\tilde\v}]$ is isomorphic to the pull back of $[C^{\v}]$
as subcone cycles in $\vt(\tilde\v)$. 
This completes the proof of the theorem.
\qed

\head
5. Composition laws of GW-invariants
\endhead

The goal of this section is to prove the composition laws of the
GW invariants.

Before we state and prove the theorem, let us introduce some 
conventions which we will use. 
In [Kn], Knudsen described various clutching morphisms, of which two 
are basic to the composition laws. We fix a partition $g_1+g_2=g$ and 
a partition $n_1+n_2=n$ once and for all. 
We always assume that $g_1$, $g_2$, $n_1$ and $n_2$
are non-negative. Let $\FF_{g,n}$ be the moduli functor of
stable $n$-pointed genus $g$ curves and let $\mgn$ be its coarse
moduli scheme. Let $S$ be any scheme, let $\xi_1\in\FF_{g_1n_1+1}(S)$
be represented by the family $C_1$ and let $\xi_2\in\FF_{g_2,n_2+1}(S)$
be represented by $C_2$. We let $\xi\in\FF_{g,n}$ be the
object represented by the family $C$ obtained by
identifying the last marked section of $C_1$ with the first marked
section of $C_2$, and set its $n$-marked sections to be the
first $n_1$ sections of
$C_1$ followed by the last $n_2$ sections of $C_2$.
In this way, we obtain the so called clutching transformation
$$\tilde\tau\lggnn: \FF_{g_1,n_1+1}(S)\times\FF_{g_2,n_1+1}(S)\lra
\FF_{g,n}(S).
$$
If there is no confusion, we will write $C=\tilde\tau\lggnn(C_1,C_2)$.
We denote the induced morphism on their moduli schemes by
$$\tau\lggnn :\mgno\times\mgnt\lra\mgn\,
$$
and call it the clutching morphism. The other clutching morphism
is defined as follows. Given $\xi\in\FF_{g-1,n+2}(S)$ represented by the
curve $C$, we obtain a new curve by identifying the last
two marked sections of $C$ and keep the initial $n$ sections. The resulting
curve is in $\FF_{g,n}(S)$. We denote this transformation by
$\tilde\tau_{g-1,n+2}$ and denote the morphism between their
moduli schemes by
$\tau_{g-1,n+2}: \M_{g-1,n+2}\lra\mgn$.

\pro{Theorem 5.1}
Let $X$ be any smooth projective variety. Assume that $\tau_i$, $\tilde\tau_i
\in H\sta(X)$ are elements so that
$[\Delta]\dual=\sum_{i=1}^k \tau_i\otimes\tilde\tau_i$
is the Kunneth decomposition of the Poincare dual of the class $[\Delta]$,
where $\Delta\sub X\times X$ is the diagonal.
Then
 \roster
\item
For any $h_1\in H\lsta(\M_{g_1,n_1+1})$ and
$h_2\in H\lsta(\M_{g_2,n_2+1})$, we have
$$\align
&\psi_{\alpha,g,n}^X(\xi_1,\cdots,\xi_n,\tau_{g_{\bullet}n_{\bullet}\,\ast}
(h_1\times h_2) \dual)\\
=&\sum_{\aota}\sum_{i=1}^k
\psi_{\alpha_1,g_1,n_1+1}^X(\xi_1,\cdots,\xi_{n_1}, \tau_i, h_1\dual)\cdot
\psi_{\alpha_2,g_2,n_2+1}^X(\xi_{n_1+1},\cdots,\xi_n, \tilde\tau_i,h_2\dual)\,.\\
\endalign
$$
\item
For any $h\in H\lsta(\M_{g-1,n+2})$, we have
$$
\psi_{\alpha,g,n}^X(\xi_1,\cdots,\xi_n,\tau_{g-1,n+2\ast}(h)\dual)
=\sum_{i=1}^k
\psi_{\alpha,g-1,n+2}^X(\xi_1,\cdots,\xi_{n}, \tau_i, \tilde\tau_i,h\dual)\,.
$$
\endroster
\endpro

We now state the composition law at the level of cycles. The numerical version
above is a direct consequence of it.
Let ${\tilde\pi}_n^{\alpha}\mh\fgna\to\FF_{g,n}$ be the transformation
that send any map in $\fgna(S)$ to the curve obtained by first
forgetting the map and then stable contract the remainder $n$-pointed
curve. Let $\pi_n^{\alpha}\mh\mgna\to\m_{g,n}$ 
be the morphism between the 
respective moduli schemes.

\pro{Theorem 5.2} 
(1). Assume $n_1,n_2>0$. We form the fiber products
$$
\CD
Z\laa@>>> \prod_{i=1}^2 \mgnai\\
@VVV @V{(e_{\nopo}, e_1)}VV\\
X @>{\Delta}>> X\times X
\endCD
\ \text{and}\
\CD
W_{\alpha} @>>> \mgnaot\\
@VVV @V{\pi\ua_n}VV\\
\prod_{i=1}^2\mgni @>{\tau\lggnn}>> \ \,\mgn\,,
\endCD
$$
where $X\mapright{\Delta}X\times X$ is the diagonal
and $e_k\mh\mgna\to X$ is the $k$-th evaluation morphism.
Then there is a canonical morphism
$\Psi:\cup_{\aota} Z\laa \to W\la$ that is finite, 
unramified and dominant. Further,
$$\Psi\lsta\bl
\sum_{\aota}\Delta\sha
\bigl[\prod_{i=1}^2 \mgnai\bigr]\vir\br=
(\tau\lggnn)\sha[\mgnaot]\vir\,.
$$

\noindent
(2). Let $Z_1$ and $Z_2$ be defined by the following fiber diagrams
$$\CD
Z_1 @>>> \M^X_{\alpha,g-1,n+2}\\
@VVV @V{(e_{n+1},e_{n+2})}VV\\
X @>{\Delta}>> X\times X\\
\endCD
\ \and\
\CD
Z_2 @>>> \mgna\\
@VVV @V{\tau\ua_n}VV\\
\M_{g-1,n+2} @>{\tau_{g-1,n+2}}>> \,\mgn\,.\\
\endCD
$$
Then the canonical morphism $\Phi\mh Z_1\to Z_2$ is finite, 
unramified and dominant. Further,
$$\Phi\lsta\bl\Delta\sha[\M^X_{\alpha,g-1,n+2}]\vir\br
=(\tau_{g-1,n+2})\sha[\mgna]\vir\,.
$$
\endpro

We will refer them as the first composition law and the second composition law.

We will give a detailed proof for the first composition law in this paper.
The proof for the second is almost identical
The only difference with the proof of the second law is that when 
$n_1,n_2>0$ the 
clutching morphism $\tau\lggnn$ is a closed immersion while 
$\tau_{g-1,n+2}$ are only locally closed
embeddings. Some modifications are required for these cases, 
which we will mention at the end of this section.

We first observe that by property 4 of Theorem 4.3, the first composition
law for $\mgna$ can be obtained from that of $\m_{\alpha,g,n+1}^X$. 
Therefore to prove the first composition law it suffices to prove the
case when $n_1,n_2>0$, which we will assume from now on.

Before we explain the strategy of the proof, let us first recall the
notion of $\QQ$-schemes that is a straightforward generalization
of $\QQ$-varieties in [Mu].

\pro{Definition 5.3}
We define a $\QQ$-scheme to be a 
scheme $Z$ with the following data.
\roster
\item
A finite atlas of charts 
$Z\lb\mapto{\pi\lb} Z\lb/G\lb\mapto{p\lb} Z$, where $p\lb$ are 
\'etale, $G\lb$ is a finite group acting faithfully on a quasi-projective 
scheme $Z\lb$ and $Z=\cup(\text{Im}\, p\lb)$.
\item
For any pair of indices $\alpha$ and $\beta$, there is a chart 
$Z_{\alpha\beta}$ with the group $G_{\alpha\beta}=
G_{\alpha}\times G\lb$ 
such that there are equivariant finite \'etale
$Z_{\alpha\beta}\to Z_{\alpha}$, $Z_{\alpha\beta}\to Z\lb$ commuting with 
projection $Z_{\alpha}$, $Z\lb$, $Z_{\alpha\beta}\to X$ such that $\text{Im}(
p_{\alpha\beta})=\text{Im}(p_{\alpha})\cap\text{Im}( p\lb)$.
\item
For any triple $\alpha,\beta$ and $\gamma$, there is a chart 
$Z_{\alpha\beta\gamma}$ with the group $G_{\alpha\beta}=
G_{\alpha}\times G\lb\times G_{\gamma}$ 
such that there are equivariant finite \'etale morphisms from
$Z_{\alpha\beta\gamma}$ to $Z_{\alpha}$, $Z_{\beta}$ and
$Z_{\gamma}$ such that in addition to $\text{Im}(
p_{\alpha\beta\gamma})=\text{Im}(p_{\alpha})\cap\text{Im}( p\lb)\cap
\text{Im}( p_{\gamma})$, the diagram
$$\CD
Z_{\alpha\beta\gamma}@>>> Z_{\alpha\beta}\\
@VVV @VVV\\
Z_{\alpha\gamma} @>>> Z_{\alpha}\\
\endCD
$$
and the other two obtained by permuting the indices are all
commutative.
\endroster
\endpro

It is known that $\mgna$ is a projective $\QQ$-scheme.
Let $\beta\in\mgna$ be any point associated to the
morphism $f\mh C\to X$ and let $G\lb$ be its automorphism group.
Then we can find an affine open $Z\lb$ acted on by $G\lb$, a
$G\lb$-equivariant family $\xi\lb\in\fgna(Z\lb)$ of stable morphisms so 
that the classifying morphism $Z\lb/G\lb\to\mgna$ induced by
the family $\xi\lb$ is an \'etale neighborhood of $\beta\in\mgna$.
For pair $\beta$ and $\beta\pri$, we can take $Z_{\beta\beta\pri}=
\bold{Iso}_{Z\lb\times Z\lbp}(\pi_1\sta\xi\lb,
\pi_2\sta\xi\lbp)$ (see [DM, p84]). 
For triple $\beta,\beta\pri$ and $\beta^{\prime\prime}$,
we can take $Z_{\beta\beta\pri\beta^{\prime\prime}}$ as a
subscheme of $Z\lb\times Z\lbp\times Z_{\beta^{\prime\prime}}$
defined similarly.
One can introduce the notion of $\QQ$-sheaves of $\o_{\mgna}$-modules, or
$\QQ$-complex, in the obvious way. A $\QQ$-sheaf is a collection
of $G\lb$-equivariant sheaves $\f\lb$ on $Z\lb$ with isomorphism
$\f\lb\otimes_{\o_{Z\lb}}\o_{Z\lbb}\cong\f\lbp\otimes_{Z\lbp}\o_{Z\lbb}$,
satisfying the cocycle condition over the triple overlaps.
Given a $\QQ$-locally free sheaves $\e$ on $\mgna$, we can define the Chern 
class $c_i(\e)$ as a cohomological class with rational coefficients,
mimicking the similar definition over Deligne-Mumford stacks.
Lastly, if $L$ is a $\QQ$-line bundle
on $\mgna$ and $s$ is a section of $L$, by which we mean a collection 
of $s\lb\in H^0(Z\lb,L\lb)$ satisfying the obvious compatibility condition 
on double overlaps, we can define the localized first Chern 
class, denoted $c_1([L,s])$, using the normal cone construction,
mimicking the construction in [Fu, \S 14.1].
One key observation, which can be checked directly,
is the following. Assume that for some integer $k$,
$L^{\otimes k}$ is a line bundle on $\mgna$. Then $c_1(L^{\otimes k})
=k\cdot c_1([L,s])$.

The strategy to prove the first composition law is quite simple, 
at least conceptually. For any partition $\aota$,
we let $\tau\laa\mh Z\laa\to\mgna$ be the clutching 
morphism that sends pairs $f_1\mh C_1\to X$ and $f_2\mh C_2\to X$ to
$f\mh C\to X$, where $C=\tilde\tau\lggnn(C_1,C_2)$ and
$f$ is the obvious induced morphism. Since $n_1,n_2>0$,
$\tau\laa$ is a closed immersion.
We will first introduce a $\QQ$-line bundle $L\laa$ on $\mgna$ 
and a section $f\laa$ such that the $\QQ$-subscheme defined
by $f\laa=0$ is the image scheme of 
$W\laa\triangleq \tau\laa(Z\laa)$. 
We will construct a $\QQ$-line bundle $L_n$ on $\mgn$ and a section 
$f_n$ so that $f_n\upmo(0)$ is
the image scheme of the clutching morphism $\tau\lggnn$.
The significance of these line bundles and sections are given
in the following lemma.

\pro{Lemma 5.4}
There are isomorphism of $\QQ$-line bundles
$$\bigotimes_{\aota} L\laa\cong (\pi\ual_n)\sta L_n\quad
$$
so that under this isomorphism, we have
$$ \prod_{\aota}f\laa=(\pi\ual_n)\sta f_n.
$$
\endpro

It follows that $\tau\lggnn\sha[\mgna]\vir$ is
$$c_1((\pi\ual_n)\sta [L_n,f_n])[\mgna]\vir
=\sum_{\aota} c_1([L\laa,f\laa])[\mgna]\vir\,.
$$
In light of Proposition 3.9, the composition law will
follow from
$$c_1([L\laa,f\laa])=[Z\laa]\vir=\Delta\sha([\mgnao\times\mgnat]\vir),
$$
which will follow from establishing the compatibility of the
tangent-obstruction complexes of the corresponding functors.

We now give the details of the proof.
We first recall some known facts.
In this section, we will view $\mgna$ and 
$\mgn$ as $\QQ$-moduli schemes. We let 
$$\{(V\lb,V\lb/G\lb, \cc\lvb)\}
_{\beta\in\Lambda_n}
$$
be an atlas of $\mgn$,
where $\cc\lvb$ are the tautological families over $V\lb$,
and let 
$$\{(W\lb,W\lb/K\lb, F\lwb,\cc\lwb)\}_{\beta\in\Lambda_{\alpha}}
$$
be an atlas of $\mgna$, where $F\lwb\mh\cc\lwb\to X$ 
are the tautological families over $W\lb$.
Assume that for some $w\in W\lb$ and $u\in V\lbp$, $\cc_{V\lbp}|_u$
is isomorphic to the stable contraction (i.e. ${\tilde\pi}\ua_n$) of 
$\cc\lwb|_w$,
then there is an \'etale  neighborhood  $(\tilde w,\tilde W\lb)\mapto{\pi}
(w, W\lb)$ such that there is a morphism $\tilde W\lb\mapright{\phi}
U\lbp$ induced by the transformation ${\tilde\pi}\ua_n$,
where the family over $\tilde W\lb$ is the pull back from $W\lb$.
Similarly, let $\beta_1$ and $\beta_2$ be any two indices and
let  $\psi_1\mh \tilde W_{\beta_1}
\to V_{\beta_1\pri}$ and $\psi_2\mh \tilde W_{\beta_2}\to V_{\beta_2\pri}$
be the so defined morphisms. Then
there is a canonical lifting $\psi_{12}\mh \tilde W_{\beta_1\beta_2}\to 
V_{\beta_1\pri\beta_2\pri}$. This follows immediately from the way 
the double overlaps are constructed.
Another property we need is the following. Let $S$ be any scheme
and $\x$ be a family in $\FF_{g,n}(S)$. Assume that for some $s\in S$ and
$u\in V\lb$, the restriction $\x|_s$ is isomorphic to the
restriction $\cc\lvb|_u$. Then there is an
\'etale neighborhood $(s_0, S_0)\mapto{\pi}(s,S)$ such that there is a
morphism $\varphi\mh (s_0,S_0)\to V\lb$ such that 
$\varphi\sta\cc\lvb\cong \pi\sta\x$ and that their restrictions to $s_0$ is exactly 
the isomorphism $\x|_s\cong \cc\lvb|_u$ given.

We first introduce $\QQ$-line bundles on $\mgnm$.
We fix an integer $k\geq 0$ and let $m=n+k$. For any 
$$K=\{h_1,\cdots,
h_{k_1}\}\sub\Sigma\triangleq\{n+1,\cdots,m\},
$$
We let $K\pri=\{h\pri_1,\cdots,h\pri_{k_2}\}$ be the complement of 
$K$. Here we assume that both $h_i$ and $h\pri_j$ 
are strictly increasing. We call $K$ and $K\pri$ a partition
of $\Sigma$. Given any $K\sub\Sigma$,
there is an obvious clutching transformation
$$\tilde \tau\lkk: \FF_{g_1,m_1+1}\times\FF_{g_2,m_2+1}
\lra \FF_{g,n},
$$ 
where $m_i=n_i+k_i$.
Namely, given $C_1\in \FF_{g_1,m_1+1}(S)$ and $C_2\in
\FF_{g_2,m_2+1}(S)$, we let $\tilde\tau\lkk(C_1,C_2)$ be the family
obtained by identifying the last section
in $C_1$ with the first section in $C_2$, and set the $m$ marked 
sections to be the union of all but the last sections of $C_1$ and all
but the first sections of $C_2$. They are ordered as follows.
The first $n$ sections are
$s_1^1,\cdots,s_{n_1}^1,s^2_2,\cdots,s_{n_2+1}^2$,
where $\{s_i^1\}$ and $\{s_j^2\}$ are sections of $C_1$ and $C_2$
respectively. We then place 
$s_{n_1+i}^1$ in $h_i$-th place and place $s_{n_2+j+1}^2$ in $h\pri_j$-th
place. We denote by $d\mh S\to C$ the section of nodal points 
along which the gluing
is taking place. We call such an nodal section decomposable
nodal section, and say $C$ is decomposable into families in
$\FF_{g_1,m_1+1}$ and $\FF_{g_2,m_2+1}$ along 
partition $K\cup K\pri$ (or along the nodal section $d$).

Since $n_1, n_2>0$, $\tau\lkk$ is a closed immersion [Kn]. 
The image scheme of
$\tau\lkk$ is a $\QQ$-Cartier divisor. Let 
$\{(U\lb,U\lb/G\lb, C\lub)\}_{\beta\in\Lambda_{m}}$ be an atlas of $\mgnm$.
For each $\beta\in\Lambda_{m}$, there is an $f\lkkb\in\o\lub$
such that $T\lkkb=\{f\lkkb=0\}$ is exactly the subscheme such that
the restriction of $C\lub$ to it is decomposable into
families in $\FF_{g_1,m_1+1}$ and $\FF_{g_2,m_2+1}$ 
along partition $K\cup K\pri$ (See [DM, p83]). If $T\lkkb=\emptyset$,
we will set $f\lkkb\equiv1$. Note that in doing
so, we might need to shrink $U\lb$ if necessary (In the
following, we shall feel free to shrink $U\lb$ and make necessary
adjustment whenever necessary). By the theory of deformation
of nodal points, $d f\lkkb|_{T\lkkb}$ is nowhere vanishing.  Furthermore, 
over $U\lbb$, there is an $f\lkkbb\in\o\sta\lubb$ such that 
$$ f\lkkb\circ\pi\lub =f\lkkbb\cdot (f\lkkbp\circ\pi\lubp)\,,
$$
where as usual we denote by $\pi\lub\mh\ubb\to\ub$ the projection.
Clearly, $\{f\lkkbb\}$ satisfies the cocycle condition
over the triple overlaps. Hence $\{f\lkkbb\}$ defines a $\QQ$-line bundle, 
denoted by $L\lkk$, and
$\{f\lkkb\}$ defines a section $f\lkk$ of $L\lkk$.
By our choice of $f\lkkb$, the image $\QQ$-scheme
$\text{Im}(\tau\lkk)$ is defined in each $\ub$ by the vanishing of $f\lkkb$.

We now turn our attention to the contraction transformation
${\tilde\pi}_n^{m}: \FF_{g,m}\lra\FF_{g,n}$
that sends any curve to the curve obtained by
 first forgetting the last $k$ sections of the curve
and then stable contracting the resulting family. 
We denote by $\pi\unm_n$ the morphism between the 
corresponding moduli schemes.
For each $\beta\in\Lambda _{m}$, there is an \'etale  covering 
$\pi_j\mh A_j\to U\lb$ of $U\lb$ and $\beta_j\pri\in\Lambda_n$
such that there is a morphism
$\pi\unm_{n,\beta,j}\mh A_j\to V_{\beta_j\pri}$ induced by the 
transformation ${\tilde\pi}\unm_n$.
For simplicity, after replacing $U\lb$ by an \'etale  covering of itself 
and rearranging the index, we
can assume that $\Lambda_n=\Lambda_{m}\triangleq\Lambda$ 
and that the 
maps $\pi\unm_{n,\beta,j}$ just mentioned are from $\ub$ to $V\lb$,
denoted $\pi\unm_{n,\beta}$ 
(We allow some $U\lb$ or $V\lb$ to be empty sets).
Clearly, $\pi\unm\lnb$ is
$G\lb$-equivariant, under the obvious group homomorphism $G\lb\to H\lb$.
Let $f\lnb\in\o_{V\lb}$ be a section whose vanishing locus defines the closed 
subscheme $T\lnb\sub V\lb$, where $T\lnb$ is the largest subscheme
 over which the family $\cc_{V\lb}$
decomposes into families in $\FF_{g_1,n_1+1}$ and $\FF_{g_2,n_2+1}$.
Let 
$$f_{n,\bbpri}=(f\lnb\circ\pi_{V\lb})/( f_{n,\bp}\circ\pi_{V\lbp}),
$$
where
$\pi_{V\lb}\mh V\lbb\to V\lb$ are the projections.
$\{f_{n,\bbpri}\}$ defines a $\QQ$-line bundle $L_n$ on $\mgn$ and $\{f\lnb\}$
is a global $\QQ$-section of $L_n$ whose vanishing locus defines the
image $\QQ$-scheme of the clutching morphism
$$\tau\lggnn: \mgno\times\mgnt\lra\mgn\,.
$$
By [DM, p83], there is a nowhere vanishing $g\lb\in\o\lub\sta$ such that
$$\prod_{K\sub\Sigma}f\lkkb=g\lb\cdot(f\lnb\circ \pi\unm\lnb)\,.
$$
By replacing one $f\lkkb$ by $f\lkkb/g\lb$, we can assume that
$g\lb\equiv 1$. Then 
$$\bigotimes_{K\sub\Sigma} L\lkk=(\pi\unm_n)\sta L_n\quad
\and \quad 
\prod_{K\sub\Sigma}f\lkk=(\pi\unm_n)\sta(f_n) \,.
$$

In the next part, for any partition $\aota$
 we will construct the $\QQ$-line bundle $L\laa$ on $\mgna$
and a section $f\laa$ of $L\laa$ such that
$f\laa=0$ defines the image $\QQ$-scheme of the 
clutching morphism $\tau\laa$.
We fix a sufficiently large $m$.
Let $\beta\in\Lambda_{\alpha}$ and let   
$(W\lb,K\lb, F\lwb,\cc\lwb)$ be the chart in the atlas of
$\mgna$. 
Without loss of generality, 
we can assume that there are sections $s_{n+1,\beta},\cdots,s_{m,\beta}
\mh \wb\to\cc_{\wb}$ such that $\cc\lwb$ with these extra sections
is a family of $m$-pointed stable curves. We
fix such a choice of new sections for each $\beta$ once and for all.
We denote the resulting $m$-pointed curve by
$\bar\cc\lwb$.  As we did before, after rearrangement we can assume that
$\Lambda_{\alpha}=\Lambda$ and that for each $\beta\in\Lambda$ there is
a morphism $\pi\ua\lnmb\mh\wb\to\ub$ such that the pull back
of $\cc\lub$ is
isomorphic to $\bar\cc\lwb$. We fix such $\pi\ua\lnmb$ once and
for all. 

Now let $\beta\in\Lambda$ be an index.
We assume that for some partition $\aota$, $p\lb\upmo(\text{Im}(\tau\laa))
\ne\emptyset$, where $p\lb\mh W\lb\to\mgna$ is the classifying map. 
Let $w$ be a point in this set.
Then there is a unique partition $K\cup K\pri=\Sigma$ 
such that the curve associated to $w$ 
with the extra sections is in the image of the clutching
transformation $\tau\lkk$. Without loss of generality,
we can assume that the partitions $K\cup K\pri$
associated to $w\in p\lb\upmo(\text{Im}(\tau\laa))$ are independent 
of $w$. Then the subscheme 
$$T\laab=\{f\lkkb\circ\pi\ua\lnmb=0\}\sub W\lb
\tag 5.1
$$
 is the closed subscheme such that the restriction of the
family $F\lwb$ to $T\laab$ is in the image of the 
clutching transformation 
$$\tilde\tau\laa: \FF_{\alpha_1,g_1,n_1+1}^X\times
\FF_{\alpha_2,g_2,n_2+1}^X\lra\fgna
$$
(Here we might need to shrink $W\lb$ so that
any element in $\{f\lkkb\circ\pi\ua\lnmb=0\}$ is over the image of
$\tau\laa$).
We let $f\laab=f\lkkb\circ\pi\ua\lnmb$. 
In case $p\lb\upmo(\text{Im}(\tau\laa))=\emptyset$, we set 
$f\laab\equiv1$. Clearly, 
should there be a unique $f\laabb\in\o\sta_{W\lbb}$ such that
$$f\laab \circ \pi\lwb=f\laabb\cdot ( f\laabp\circ\pi_{W\lbp})
$$
for all pairs $\{\beta,\beta\pri\}$, where $\pi_{W\lb}\mh W\lbb\to W\lb$ 
is the projection,
then $\{f\laabb\}$ would satisfy the desired cocycle condition and then 
define a
$\QQ$-line bundle, and $\{f\laab\}$ would define us a global section
whose vanishing locus is the image of $\tau\laa$, as desired.
However, for our moduli space $\mgna$, it may happen
that some $f\laab\equiv0$. Thus the 
above direct argument needs to be modified. In the following,
we will first thicken the $\wb$'s so that the desired
$f\laabb$ are still well-defined.

We fix a two-term complex $\edot$ of locally free sheaves 
such that its sheaf cohomology 
$\Hbul(\edot)=\tdot\fgna$. Such complexes were constructed in Section 4. 
Over each $W\lb$, $\tdot\fgna$ is represented by a tangent-obstruction complex
$\tdot_{W\lb}$ of $\wb$, and $\edot$ is represented by a two-term complex of 
locally free sheaves of $\o_{\w\lb}$-modules, denoted
$\edot\lb=[\e\lbo\to\e_{\beta,2}]$, such that $\Hbul(
\edot\lb)=\tdot_{W\lb}$. Let $E\lb$ be the vector bundle on $\wb$ such that
$\o\lwb(E\lb)=\e\lbo$. By abuse of notation, we will also use
$E\lb$ to denote the total space of $E\lb$.
Because the way $\edot$ is constructed, there is a canonical homomorphism
of locally free sheaves
$$\varphi\lb: \e\lbo\lra\ext^1_{\cc\lwb/\wb}(\Omega_{
\cc\lwb/W\lb}(D\lwb),\o_{\cc\lwb}\br\,,
$$
where $D\lwb$ is the divisor of $n$-marked sections of $\cc\lwb$. 
By the deformation theory, if we let $J\sub\o_{E\lb}$ be the
ideal sheaf of $\wb\sub E\lb$, embedded via the zero section,
and let $F\lb\sub E\lb$ be the subscheme defined by the
ideal $J^2$, then there is a family of $n$-pointed curves $\cc\lfb$
over $F\lb$ such that the homomorphism $\varphi\lb$ above is the 
Kodaira-Spencer map of the family $\cc\lfb$ along the 
normal bundle to $\wb$ in $F\lb$. Here by the normal bundle to
$\wb$ in $F\lb$ we mean the subbundle of $TF\lb|_{W\lb}$
that is the kernel of $TF\lb|_{W\lb}\to TW\lb$.
Since the restriction of $\cc\lfb$ to $\wb\sub F\lb$ is $\cc\lwb$ and 
$\bar\cc\lwb$ is smooth over $\wb$ near its sections, we can
extend the last $k$ sections of $\bar\cc\lwb$ to sections in
$\cc\lfb$ (over $F\lb$)
and place them in the same order as $\bar\cc_{\w\lb}$ has.
We denote the resulting $m$-pointed
curve by $\bar\cc\lfb$, which is a stable curve over $F\lb$. 
Let 
$$\bar\pi^{\alpha}_{m,\beta}: F\lb\to\ub
$$
 be the morphism induced by the family $\bar\cc_{F\lb}$. 
Then $\bar\pi^{\alpha}_{m,\beta}$ is an extension of 
$\pi\ua\lnmb\mh\wb\to\ub$.

Now we go back to the partition $K\cup K\pri=\Sigma$
and the ${T\laab}$ defined in (5.1).
We let $\bar f\laab$ be $ f\lkkb\circ\bar\pi\ua_{m,\beta}\in\o\lfb$.
Then $\bar T=\{\bar f\laab=0\}$ is exactly the subscheme
of $F\lb$ such that the restriction of $\bar \cc_{F\lb}$ to $\bar T$ 
belongs to the image of the clutching transformation
$$\tilde\pi_K: \FF^X_{\alpha_1,g_1,m_1+1}\times\FF^X_{
\alpha_2,g_2,m_2+1}
\lra\fgna.
$$
We fix such a $\bar\cc\lfb$, and then
$\bar\pi\ua\lnmb$ and $\bar f\laab$ for each $\beta\in\Lambda$ once and
for all. Recall that $\pi_{W\lb}\mh W_{\bbpri}\to W\lb$ is the projection.
Since there is a canonical isomorphism 
$$(\pi\lwb)\sta E\lb\cong (\pi\lwbp)\sta E_{\beta\pri},
$$
$F\lb\times_{W\lb}W\lbb$ 
is canonically isomorphic to 
$F\lbp\times_{W\lbp}W \lbb$. We denote this scheme by $F\lbb$. We
let $\pi\lfb\mh F\lbb\to F\lb$ be the projection.

\pro{Lemma 5.5}
There are $\phi\in\o_{F\lbb}$ such that 
$$\bar f\laab\circ\pi_{F\lb}=\phi\cdot(
\bar f\laabp\circ\pi_{F\lbp}).
$$
Further, if we let $\iota\mh W\lbb\to F\lbb$ be the inclusion, then
$\phi\circ\iota$ is unique and is nowhere vanishing. We denote $\phi\circ\iota$ by 
$f\laabb$. $\{f\laabb\}$ defines a $\QQ$-line bundle on $\mgna$ and
$\{f\laab\}$ defines a section of this line bundle.
\endpro

\demo{Proof}
The existence of $\phi$ follows from the fact that the subscheme
$\{\bar f\laab\circ\pi_{F\lb}=0\}$ is identical to the subscheme
$\{\bar f\laabp\circ\pi_{F\lbp}=0\}$. This is true because both define the
subschemes over which the two families of curves,
which are isomorphic after discarding the last $k$ sections,
can be decomposed along the same nodal sections (See [DM, p83]).
To show that $\phi\circ\iota$ is unique, it suffices to show that at each 
$w\in W\lbb$ there is a tangent
vector $v\in T_w F\lbb$ contained in the kernel of $T_w F\lbb\to T_w W\lbb$
such that $d(\bar f\laab)(v)\ne0$. Here $T_w F\lbb\to T_w W\lbb$
is induced by the vector bundle projection $E\lb\to\wb$. This is true because
$$\e\lbo\otimes k(w)\lra \Ext^1_{\cc_w}(\Omega_{\cc_w}(D_w),\o_{\cc_w})
\lra\ext^1_{\hat r}(\Omega_{\hat r},\o_{\hat r})
$$
is surjective. Here $\cc_w$ is the fiber of $\cc\lub$ over $w$, 
$D_w\sub\cc_w$ is the
divisor of $n$-marked points, $r\in\cc_w$ is the nodal point where the clutching
$\tau\laa$ is taking place and $\hat r$ is the formal completion of
$\cc_w$ along $r$. Because $f\laabb$ are unique, the collection
$\{f\laabb\}$ satisfies the cocycle condition on triple overlaps.
Therefore it defines a $\QQ$-line bundle on $\mgna$, denoted by $L\laa$.
For the same reason, $\{f\laab\}$ form a global section, denoted by 
$f\laa$, of $L\laa$ such that the image scheme $\text{Im}(
\tau\laa)\sub\mgna$ is defined over each $W\lb$ by $\{f\laab=0\}$.
This proves the lemma.
\qed

Now, we let $\pi\ua_n\mh\mgna\to\mgn$ be the stable contraction
morphism. Locally, $\pi\ua_n$ is represented by maps
$$\pi\ua\lnb=\pi\unm\lnb\circ\pi\lnmb\ua:\wb\lra V\lb\,.
$$
Let $\pi_{n,\bbpri}\ua\mh\wbb\to V\lbb$ be the lifting of the pair $(\pi\ua\lnb,
\pi\ua_{n,\bp})$, which exists and is unique. 
Let $f\lab\in\o\lwb$ be $f_{n,\beta}\circ\pi\ua\lnb$ and let $f_{\alpha,\bbpri}
=f_{n,\bbpri}\circ\pi\ua_{n,\bbpri}\in\o\lwbb\sta$. 
Then $\{f_{\alpha,\bbpri}\}$ defines a 
$\QQ$-line bundle, denoted by $L_{\alpha}$, which is isomorphic to
the pull back $(\pi\ua_n)\sta L_n$. For the same reason, $\{f\lab\}$ is
a section $f\la$ of $L\la$, which is the pull back of $f_n$ of
$L_n$. 

\pro{Lemma 5.6}
There is a canonical isomorphism of $\QQ$-line bundles
$\otimes L\laa\cong L\la$
such that under this isomorphism, we have 
$\Pi  f\laa=f\la$, where the product is over all possible $\aota$.
\endpro

\demo{Proof}
First note that given any $\aota$ and $w\in p\lb\upmo(\text{Im}(\tau\laa))$,
there is a unique partition $K\cup K\pri=\Sigma$ such that
the decomposition of $\bar\cc_w$ along partition $K\cup K\pri$ is the inverse 
to the clutching transformation $\tilde\tau\laa$. 
We denote this correspondence by
$\alpha_1\alpha_2\mapsto K_w(\alpha_1\alpha_2)$. By making each 
$W\lb$ small enough, we can assume that the following two conditions hold
automatically for all $\beta\in\Lambda$. First, for $w,w\pri\in
p\lb\upmo(\text{Im}(\tau\laa))$ we have
$K_w(\alpha_1\alpha_2)=K_{w\pri}(\alpha_1\alpha_2)$; Second,
the map
$$
\Bigl\{(\alpha_1,\alpha_2)\vert
\matrix \aota \ \and \\ p\lb\upmo(
\text{Im}(\tau\laa))\ne\emptyset \endmatrix
\Bigr\}\lra
\Bigl\{ (K,K\pri) \vert 
\matrix 
\cc_{U\lb}|_u \ \text{is decomposable along}\\
\ K\sub\Sigma\ \text{for some}\ u\in\ub
\endmatrix
\Bigr\}
$$
is injective and onto. The first is possible because $\tau\laa$ is a closed embedding.
We define $K\lb(\alpha_1\alpha_2)=K_w(\alpha_1\alpha_2)$ for some
$w\in p\lb\upmo(\text{Im}(\tau\laa))$. Hence for each $\beta$,
$$\prod_{K\sub\Sigma} f\lkkb=f\lnb\circ\pi\lnb\unm\,.
\tag 5.2
$$
Here if $K\cup K\pri$ is not in the image of the above correspondence, 
then $f\lkkb\equiv1$.
To show that $\otimes L\laa=L\la$, we need to check 
$\Pi f\laabb=f\labb$, which amounts to prove the identity
$$\prod_{\aota} { (\pi\ua\lnmb\circ\pi\lwb)\sta f_{K\lb(\aa),\beta}\over
(\pi\ua_{m,\bp}\circ\pi\lwbp)\sta f_{K\lbp(\aa),\beta\pri} }
={(\pi\ua\lnb\circ\pi\lwb)\sta f\lnb\over 
(\pi\ua_{n,\bp}\circ\pi\lwbp)\sta f_{n,\bp}}\,,
\tag 5.3
$$
on $\wbb$.
This is obvious from the previous identity if all quotients are
well-defined and unique. To prove (5.3),
we first embed $F\lbb\supset W\lbb$ into a smooth affine
scheme, say $F\lbb\sub R$. Without loss of generality, we can assume
$W\lb$, $U\lb$ and $W\lbb$ are affine.
Let $\phi\mh R\to\ub$ be a morphism  extending the 
morphism 
$$\psi\lb: F\lbb @>{\pi\lfb}>>  F\lb @>{\bar\pi\ua_{m,\beta}}>>\ub\,.
$$
We claim that there is a morphism $\phi\pri\mh
R\to\ubp$ extending $\psi\lbp\mh F\lbb\to\ubp$ such that 
$$p\lb\circ\pi\lnb\unm\circ\phi=
p\lbp\circ\pi_{n,\bp}\unm\circ\phi\pri\,,
\tag 5.4
$$
where $p\lb\mh V\lb\to\mgn$ is the obvious morphism.
Indeed, the morphism $\phi\mh R\to\ub$ provides us a family of
$m$-pointed stable curves over $R$ via pull back. Let it be $\x_R$ with
the sections $\{s_i\}$. Let $\psi\lbp\sta\cc\lubp$ be the pull back 
of $\cc_{U\lbp}$ via $\psi\lbp\mh F\lbb \to U\lbp$. 
Then after restricting $\x_R$ to $F\lbb$ and discarding its last $k$ sections,
the resulting family is isomorphic to $\psi\lbp\sta\cc\lubp$ as $n$-pointed
curve. We choose $k$ sections $\tilde s_{n+1},\cdots
\tilde s_{m}$ of $\x_R$ extending the last $k$ sections of 
$\psi\lbp\sta\cc\lbp$, which is
possible because $\x_R$ is smooth over $R$ along sections 
of $\psi\lbp\sta\cc\lbp$.
We denote $\x_R$ with these
new sections $(\cdots,s_n,\tilde s_{n+1},\cdots)$ by 
$\tilde \x_R$. Then $\tilde\x_R$ induces a morphism $\phi\pri
\mh R\to\ubp$ that has the desired property.

Now we show that $\Pi f\laabb=f\labb$. Because of (5.4) above
and because $\phi\pri$ is induced by the family of curves,
we have a unique lifting $\bar\phi\mh R\to V\lbb$. 
Then because $R$ is smooth and
because none of $f_{K\lb(\aa),\beta}\circ\phi$ and 
$f_{K\lbp(\aa),\beta\pri}\circ\phi\pri$ are zero, there is a unique
rational $h\laabb$ such that 
$$f_{K\lb(\aa),\beta}\circ\phi=h\laabb\cdot(
f_{K\lbp(\aa),\beta\pri}\circ\phi\pri)\,.
$$
Clearly, $h\laabb$ is regular and non-vanishing 
near $F_{\beta\beta\pri}\sub R$. It follows that
$$\prod_{\aota} ( f_{K\lb(\aa),\beta}\circ\phi)
=\bl\prod_{\aota} h\laabb\br\bl\prod_{\aota}
(f_{K\lbp(\aa),\beta\pri}\circ\phi\pri)\br\,.
\tag 5.5
$$ 
Now let $\tau\lbb\mh W\lbb\to R$ be the inclusion.
Of course, $h\laabb\circ\tau\lbb\equiv f\laabb$. Thus 
to show that $\otimes L\laa=(\pi\ua_n)\sta L_n$, it suffices to
show that 
$$(\prod h\laabb)\circ\tau\lbb= f_{n,\beta\beta\pri}\circ\pi\ua_{n,\bbpri},
$$
where $\pi\ua_{n,\bbpri}\mh W\lbb\to V\lbb$ is the projection.
Indeed, because $\Pi_{K\sub\Sigma} f\lkkb =f\lnb\circ
\pi\unm\lnb$,
$$g\lb:= f_{n,\beta}\circ
\pi\unm\lnb\circ\phi=\prod_{\aota} f_{K\lb(\aa),\beta}\circ\phi
$$
Therefore (5.5) holds with $\Pi h\laabb$ replaced by
$g\lb/g\lbp$.
Since all terms in (5.5) are non-trivial and since $R$ is smooth,
$g\lb/g\lbp$ is regular near $\wbb$ and its composite with
$\tau\lbb$ is identical to the composite of
$\Pi h\laabb$ with $\tau\lbb$. Finally, because this term is also the pull back
via $\bar\phi\mh R\to V\lbb$ of 
$f\lnb\circ\pi_{V\lb}/f_{n,\bp}\circ\pi_{V\lbp}$, where
$\pi_{V\lb}\mh V\lbb\to V\lb$ is the projection,
its composite with $\tau\lbb$ is the pull back of $f_{n,\bbpri}$.
This proves that
$\Pi f\laabb=f_{\alpha,\bbpri}$ and consequently 
$\otimes L\laa=(\pi\ua_n)\sta L_n$.
It remains to show that $\Pi f\laab=f\lab$. This is true because 
$\Pi_{K\sub\Sigma} f\lkkb=f_{n,\beta}\circ\pi\unm\lnb$.
This completes the proof of Lemma 5.5.
\qed

As was explained before, to complete the proof of the
first composition law, we remain to 
investigate the tangent-obstruction complex of $Z\laa$.
We first introduce the functor $\FF\laa$ so that its coarse 
moduli scheme is $Z\laa$. For any scheme $S$, we let
$\FF\laa(S)$ be the subset of $\fgna(S)$ consisting of
families $f\mh\x\to X$ such that there are distinguished 
sections of nodal points $d\mh S\to\x$ such that $f$ are
decomposable to pairs of families in $\FF^X_{\alpha_1,g_1,n_1+1}(S) 
\times\FF_{\alpha_2,g_2n_2+1}^X(S)$ along the nodal sections $d$.
Clearly, $\FF\laa$ is coarsely
represented by a scheme $V\laa$ which is canonically isomorphic
to the $Z\laa$ defined before. By forgetting the distinguished sections,
the resulting transformation $\FF\laa\to\fgna$ defines a morphism
$\Phi\mh Z\laa\to \mgna$. $\Phi$ is a closed
immersion since $n_1, n_2>0$ and induces an
isomorphism $Z\laa\cong W\laa$.
In the following, we will
not distinguish $Z\laa$ from $W\laa$ unless otherwise is mentioned.

We now give the tangent-obstruction complex
of $\FF\laa$. Let $S$ be any affine 
scheme and let $\xi\in\FF\laa(S)$ be represented by
$f\mh\x\to X$ with the marked divisor $D\sub\x$ and
the distinguished nodal section $d\mh S\to \x$.
Let $\bbul(\xi)$ be the complex 
$[f\sta\Omega_X\to\Omega^0_{\x/S}(D)]$ indexed at $-1$ and 0, 
where $\Omega^0_{\x/S}$ is the quotient sheaf of $\Omega_{\x/S}$
by its torsion supported at the distinguished nodal section
and $f\sta\Omega_X\to\Omega^0_{\x/S}$ is induced by
$f\sta\Omega_X\to\Omega_{\x/S}$. Because the 
first order deformations
of nodal curves with distinguished nodal sections
are $\Ext^1_{\x}(\Omega_{\x/S}^0,\o_{\x})$, from the description
of the tangent-obstruction complex $\tdot\fgna$, we see
immediately that
$$\tone\FF\laa(\xi)(\f)=\ext^1_{\x/S}(\bbul(\xi),\pi_S\sta\f),
$$
where $\pi_S\mh \x\to S$ is the projection, and that there is
an obstruction theory to deformations of $\FF\laa$ taking values
in 
$$\ttwo\FF\laa(\xi)(\f)=\ext^2_{\x/S}(\bdot(\xi),\pi_S\sta\f).
$$
Now let $\hat T$ be the formal completion of $\x$ along
$d(S)$. It follows that we have the long exact sequence
$$\align
0\lra\tone\FF\laa(\xi)(\f)\lra\tone&\fgna(\xi)(\f)\mapright{c}
\ext_{\hat T/S}^1(\Omega_{\hat T/S},\pi_S\sta\f)\lra\\
&\mapright{b}\ttwo\FF\laa(\xi)(\f)
\mapright{}\ttwo\fgna(\xi)(\f)\lra 0.
\endalign
$$
Note that $\ext_{\hat T/S}^1(\Omega_{\hat T/S},\pi_S\sta\f)\cong
\f$.
It follows from [DM] that
$\ext_{\hat T/S}^1(\Omega_{\hat T/S},\o_{\hat T})$ is 
isomorphic to the pull back of the $\QQ$-invertible
sheaf $\o_{\mgna}(L\laa)$ and the homomorphism $b$ is induced
by the differential of the defining equation $f\laa$.
Now let $S\to Y_0\to Y$ be the triple described in Definition
1.2 and let $\xi_0$ is a family in $\FF\laa(Y_0)$ extending $\xi$.
Let $o$ (resp. $\tilde o$) be the obstruction class to extending $\xi_0$
to $\FF\laa(Y)$ (resp. to $\fgna(Y)$).
A straightforward analysis of the definitions of $o$ and $\tilde o$ 
shows that $\tilde o=c(o)$. 
Next, assume that $\tilde o=0$. Namely $\xi_0$
extends to a family $\tilde f\mh\tilde\x\to X$ over $Y$.
Let $h\in\Ext_{\hat T}^1(\Omega_{\hat T/S},\pi_S\sta\i_{Y_0\sub Y})$
be the section associated to the family $\tilde\x\to Y$. Then
it is direct to check that
$o$ is the image of $h$ under $b$ (See [DM]).
This proves that the obstruction theory of $\FF\laa$ and $\fgna$
are compatible with respect to the defining equation $f\laa=0$. Finally,
following the construction of the complex $\edot$ in the
beginning of the Section 4, we can find complexes
$\edot$ and $\fdot$ such that $\hhdot(\edot)=\tdot\FF\laa$ and
$\hhdot(\fdot)=\tdot\fgna$ that satisfy the technical condition of
Proposition 4.9. Therefore, by the argument at the
end of the proof of Theorem 4.2, we can apply Proposition 4.9 to 
conclude
$$c_1([L\laa,f\laa])[\mgna]\vir=[Z\laa]\vir.
\tag 5.4
$$

In the following discussions, we will abbreviate $\FF_{\alpha_i,g_i,n_i+1}^X$
to $\FF_{\alpha_i}$.
It remains to show that $\tdot\FF\laa$ is compatible to the
tangent-obstruction complex of $\tdot(\FF_{\alpha_1}\times
\FF_{\alpha_2})$ with respect to the fiber product
$$\CD
\FF\laa(S) @>>> 
\FF_{\alpha_1}(S)\times\FF_{\alpha_2}(S)\\
@VVV @VVV\\
\bold{Mor}(S,X) @>{\Delta}>>\bold{Mor}(S,X\times X).
\endCD
\tag 5.5
$$
Here, given $\xi\in\FF\laa(S)$ represented by the map $f$
and the distinguished section $d$,  
the first vertical arrow will send it to $f\circ d\mh S\to X$.
Similarly, for $(\xi_1,\xi_2)\in\FF_{\alpha_1}\times\FF
_{\alpha_2}(S)$ represented by the pair of
maps $f_1$ and $f_2$, the second arrow will send it to
$(e_{n_1+1},e_1)\mh S\to X\times X$, where $e_{n_1+1}$ is the
$(n_1+1)$-th evaluation map of $f_1$ and $e_1$ is the first 
evaluation map of $f_2$.
Let $\xi\in\FF\laa(S)$ be as before. By decomposing $f$ along
its distinguished section $d$, we obtain a pair of family
$\{f_i\mh D_i\sub\x_i\to X\}=\xi_i\in\FF_{\alpha_i}(S)$
for $i=1,2$. We let 
$$\rdot_i(\xi)=[f_i\sta\Omega_X(-d_i(S))\to
\Omega_{\x_i/S}(D_i-d_i(S))],
$$
where $d_1$ is the last marked section of $D_1\sub\x_1$
and $d_2$ is the first marked section of $D_2\sub\x_2$.
Let $\iota_i\mh\x_i\to\x$ be the immersion. Then we have the exact
sequence
$$0\lra\iota_{1\ast}\rdot_1(\xi)\oplus\iota_{2\ast}\rdot_2(\xi)\lra
\bdot(\xi)\lra [f\sta\Omega_X\otimes_{\o_{\x}}\o_{d(S)}\to0]\lra0
$$
and its induced long exact sequence
$$\align
0\lra&\ext^1_{\x/S}(\bdot(\xi),\pi_S\sta\f)\lra
\oplus_{i=1}^2\ext^1_{\x/S}(\iota_{i\ast}\rdot_i(\xi),\pi_S\sta\f)\lra
\\&
\mapright{h_1} \ext^1_{\x/S}(f\sta\Omega_X\otimes_{\o_{\x}}\o_{d(S)},
\pi_S\sta\f)
\mapright{h_2} \ext^2_{\x/S}(\bdot(\xi),\pi_S\sta\f)\lra
\tag 5.6
\\&\qquad\qquad
\lra\oplus_{i=1}^2\ext^2_{\x/S}(\iota_{i\ast}\rdot_i(\xi),\pi_S\sta\f)\lra 0.
\endalign
$$
It is direct to check that there are canonical isomorphisms
$$\ext^j_{\x/S}(\iota_{i\ast}\rdot_i(\xi),\pi_S\sta\f)\cong
\ext^j_{\x_i/S}([f_i\sta\Omega_X\to\Omega_{\x_i/S}(D_i)],\pi_S\sta\f)
$$
and 
$$\ext^1_{\x/S}(f\sta\Omega_X\otimes_{\o_{\x}}\o_{d(S)},\pi_S\sta\f)\cong
(f\circ d)\sta\Omega_X\dual\otimes_{\o_S}\f.
$$
Therefore, (5.6) is the long exact sequence of cohomologies 
mentioned in Definition 3.8. Now we show that the tangent-obstruction
complexes of $\FF\laa$ and of 
$\FF_{\alpha_1}\times\FF_{\alpha_2}$
are compatible with respect to the defining diagram (5.5).
Let $\xi\in\FF\laa(S)$ be as before and let $\xi_i\in
\FF_{\alpha_i}(S)$ be represented by the family $f_i$.
We first show that the canonical homomorphism
$$\tone(\FF_{\alpha_1}\times\FF_{\alpha_2})(\xi_1,\xi_2)(\f)\lra
(e_{n_1+1},e_1)\sta\o(T_{X\times X})\otimes_{\o_S}\f
\lra (f\circ d)\sta\Omega_X\dual\otimes_{\o_S}\f
$$
is the homomorphism $h_1$ in the exact sequence (5.6).
Let $S$ be affine, let $\f\in\mods$ and let 
$v_1\in\tone\FF_{\alpha_1}(\xi_1)(\f)$
be represented by a flat extension $\tilde f_1$ of $f_1$. 
Let
$$
\CD
@.@. f_1\sta\om_X @= f_1\sta\om_X\\
@.@. @VVV @VVV\\
0@>>> \pi_S\sta\f @>>> \bb @>>> \om_{\x_1/S}(D_1) @>>> 0
\endCD
$$
be the associated diagram. We first set
$\Cal B\pri$ be the sheaf defined by the push-forward diagram
$$
\CD
\pi_S\sta\f(-d_{1}(S)) @>>> \Cal B(-d_{1}(S))\\
@VVV @VVV\\
\pi_S\sta\f @>>>\ \, \Cal B\pri\,.
\endCD
$$
Then $f_1\sta\om_X(-d_1(S)) \to\bb\pri$
lifts to $f_1\sta\om_X\to\bb\pri$. Let $\bdot_1(\xi_1)=
[ f_1\sta\om_X\to\om_{\x_1/S}(D)]$ and let
$p(\xi_1)$ be the homomorphism
$$\tone\FF_{\alpha_1}(\xi_1)(\f)=
\ext^1_{\x_1/S}\bl \bdot_1(\xi_1),\pi_S\sta\f\br
@>{p(\xi_1)}>>
e_{n_1+1}\sta\om_X\dual\otimes_{\os}\f
$$
that sends $v_1$ to the composite
$$e_{n_1+1}\sta\om_X\lra e_{n_1+1}\sta\bb\pri\lra
\coker\{\bb(-d_1(S)) \to\bb\pri\}\equiv \f\,.
$$
One checks directly that $p(\xi_1)$ assigns the flat extension
$\tilde f_1$ to the tangent direction of 
$\tilde f\circ\tilde e_{n_1+1}$, where $\tilde e_{n_1+1}$ is the 
$(n_1+1)$-th marked section of $\tilde f$.
Similarly, we let $p(\xi_2)$ be
$$p(\xi_2): \tone\FF_{\alpha_2}(\xi_2)(\f)=
\ext^1_{\x_2/S}\bl \bdot_2(\xi_2),\pi_S\sta\f\br
\lra e_1\sta\om_X\dual\otimes_{\os}\f
$$
that is defined with $f_1$, etc.
replaced by $f_2$, etc. respectively.
Then $p(\xi_1)-p(\xi_2)$ is the homomorphism $h_1$ in the exact
sequence (5.7).
Therefore, the induced homomorphism $(e_{n_1+1},e_1)$
on the tangent spaces
$$
\tone (\FF_{\alpha_1}\times \FF_{\alpha_2})
(\xi_1,\xi_2)(\f) 
\lra
(e_{n_1+1},e_1)\sta\om_{X\times X}\dual
\otimes_{\os}\f
$$
coincides with $(p(\xi_1),p(\xi_2))$, and consequently, its
composite with 
$$(e_{n_1+1},e_1)\sta\om_{X\times X}\dual 
\otimes_{\os}\f
\lra
(e_{n_1+1},e_1)\sta \n_{\Delta(X)/X\times X}
\otimes_{\os}\f
$$
is $p(\xi_1)-q(\xi_2)$, after identifying the normal bundle 
$N_{\Delta(X)/X\times X}$
with $T_{X}$. Similarly, it is direct to check that the
obstruction classes to extending a given family to families in
$\FF\laa$ and to families in $\FF_{\alpha_1}\times
\FF_{\alpha_2}$ are compatible in the sense of
Definition 3.8. 
Finally, similar to the case studied before, we can construct
complexes of locally free sheaves whose sheaf cohomologies
are the tangent-obstruction complexes of 
$\FF_{\alpha_1}\times\FF_{\alpha_2}$ required by Proposition 3.9.
Therefore, by applying Proposition 3.9, we have
$$\Del\sha\bl[\mgnao\times\mgnat]\vir\br=[Z\laa]\vir\,.
\tag 5.7
$$
Finally, we choose a sufficiently large $l$ such that $L\laa^{\otimes l}$ 
and $L_n^{\otimes l}$ are conventional line bundles. Then
$$
\align
&l\cdot(\tau\lggnn)\sha[\mgna]\vir=c_1((\pi\ua_n)\sta 
L_n^{\otimes l})[\mgna]\vir
\\=&\sum_{\aota} c_1(L\laa^{\otimes l})[\mgna]\vir
=\sum_{\aota} l\cdot c_1([L\laa,f\laa])[\mgna]\vir.
\endalign
$$
By (5.4), the terms in the last summation are
$l\cdot [Z\laa]\vir$. Combined with (5.7), we obtain
$$(\tau\lggnn)\sha[\mgna]\vir=
\sum_{\aota}\Delta\sha([\prod_{i=1}^2\mgnai]\vir).
$$
This proves the first composition law.

In the end, we will indicate the necessary change needed to prove the 
second composition law. 
Let $Z_1$ and $Z_2$ be the $\QQ$-schemes and $\Phi$ be the 
morphism defined in the statement of the theorem. For convenience, we 
will consider $Z_1/\ZZ_2$ and $Z_2/\ZZ_2$, where $\ZZ_2$ acts on $Z_1$
and $Z_2$ by interchanging the last two marked points of the curves
in $\m_{\alpha,g-1,n+2}^X$ and $\m_{g-1,n+2}$. Let 
$Z_i\pri= Z_i/\ZZ_2$.
Clearly, $\Phi$ factor through $\Psi\mh Z_2\pri\to Z_1\pri$.
$\Psi$ is a local embedding in the sense that
it is finite and unramified. Let $(L_n,f_n)$ be the
$\QQ$-line bundle and its section on $\mgn$ such that $f_n=0$ defines the
image $\QQ$-scheme $\m_{g-1,n+2}\to\mgn$. We pick a 
$w\in Z_2$ and let $\{z_1,\cdots,z_k\}=\Psi\upmo(w)$.
Now let $W\to\mgna$ be a chart of $\mgna$ containing $w$ with
the tautological family $\xi$. We let $U_i\to Z_1\pri$ be
the charts
of $Z_1\pri$ containing $z_i$ with the tautological family
$\eta_i$. Recall that each family $\eta_i$ has a distinguished
section of nodal points. Without loss of generality, we can assume that
there are morphisms $\varphi_i\mh U_i\to W$
such that $\varphi_i\sta(\xi)=\eta_i$. Now by using the 
technique of adding extra sections, we can find an \'etale
covering $\tilde W\to W$, $k$ sections $g_1,\cdots,g_k\in
\o_{\tilde W}$ and \'etale covering $\tilde U_i\to U_i$
of which the following holds.
Firstly, after fixing a trivialization of $L_n$ over a
chart $V$ of $\mgn$, the product $\Pi^k g_i$ is the pull back of
$f_n$ under the obvious map $\tilde W\to V$;
Secondly, there are morphisms $\tilde\varphi_i\mh\tilde V_i\to \tilde W$ 
making the following diagram 
$$\CD
\tilde U_i @>{\tilde\varphi_i}>> \tilde W\\
@VVV @VVV\\
U_1 @>{\varphi_i}>> W\\
\endCD
$$
commutative such that $\tilde\varphi_i$ are embeddings and the
image schemes $\tilde\varphi_i(\tilde U_i)=\{g_i=0\}$.
Using the distinguished section of nodal points in the
family $\eta_i$, one can construct a $\QQ$-invertible sheaf 
$\l$ on $Z_1\pri$ such that over the chart $\tilde U_i$, it is 
the locally free sheaf defined similar to the far right
term in (5.2). For convenience. let us assume 
that $[\mgna]\vir$ is a cycle $R$ supported on an equidimensional
scheme with multiplicity. Now let $Y_i=\{g_i=0\}\sub\tilde W$ and
let $R_{\tilde W}$ be the pull back of $R$ under $\tilde W
\to\mgna$. Consider the normal cone cycle 
$$[C_{R_{\tilde W}\times_{\tilde W}Y_i/R_{\tilde W}}].
$$
Using the isomorphism $\tilde\varphi_i\mh\tilde U_i\to Y_i$,
one can pull back the cycles 
$[C_{R_{\tilde W}\times_{\tilde W}Y_i/R_{\tilde W}}]$
and patch them together to form a global
cycle in the total space of
$$\vt_{Z_1\pri}(\l)\times_{\mgna}R.
$$
We denote this cycle by $\bold D$ and by $\zeta$
the zero section of $\vt_{Z_1\pri}(\l)\times_{\mgna} R$.
Then by studying the tangent-obstruction complex of 
$Z_1\pri$ induced by the defining equation $g_i=0$ and
that of $Z_1\pri$ induced by the defining square of
$Z_1$ in the statement of the theorem, on concludes that
$$2\Psi\lsta(\zeta\sta[\bold D])=\iota\lsta
\Phi\lsta(\Delta\sha[\m_{\alpha,g-1,n+2}^X]\vir),
$$
where $\iota\mh Z_2\to Z_2\pri$ is the projection.
However, it is clear that
$$\iota\lsta(\tau_{g-1,n+2})\sha[\mgna]\vir=2c_1(
[(\tau_n^{\alpha})\sta(L_n),(\tau_n^{\alpha})\sta(f_n)])(R).
$$
Therefore, the second composition law will follow from
$$ c_1( [(\tau_n^{\alpha})\sta(L_n),(\tau_n^{\alpha})\sta(f_n)])(R)
=\Psi\lsta(\zeta\sta[\bold D]).
$$
But this can can be checked directly. This proves
the second composition law.

\enddocument

\bye